\documentclass[11pt]{amsart}
\usepackage{amssymb, amsmath, amsthm}
\usepackage[latin1]{inputenc}
\usepackage{graphicx}
\usepackage[final]{hyperref}
\usepackage[a4paper, centering]{geometry}
\usepackage{color}
\usepackage{graphicx} 
\usepackage{wrapfig}
\usepackage{tikz,amsmath,amsthm,bbm}
\usepackage{xypic}
\usepackage[T1]{fontenc}
\usepackage{subfigure}

\usepackage[english]{babel}

\usepackage{mathtools}

\numberwithin{equation}{section}
\newtheorem{theorem}{Theorem}[section]
\newtheorem{proposition}[theorem]{Proposition}
\newtheorem{lemma}[theorem]{Lemma}
\newtheorem{corollary}[theorem]{Corollary}

\theoremstyle{definition}
\newtheorem{definition}[theorem]{Definition}

\theoremstyle{remark}
\newtheorem{remark}[theorem]{Remark}

\def\A{\widehat {\mathcal A}}

\def\X{\mathcal B}
\def\K{\mathcal K}
\def\R{\mathbb{R}}

\def\N{\mathbb{N}}

\def \o{\overline}

\font\maius=cmcsc10 scaled1200
\def \div {\mathop {\rm div}\nolimits}
\def \esssup{\mathop {\hbox{\rm - ess sup}}\limits}
\def\e{\varepsilon}
\def\1{{{\bf 1}
\kern-0,28em \rm l}}  
\def \res{\mathop{\hbox{\vrule height 7pt width .5pt depth 0pt
\vrule height .5pt width 6pt depth 0pt}}\nolimits}
\def \S{{\mathcal I} ^ * }

\begin{document}

\title[Duality for non-convex variational problems]%
{A duality theory for non-convex problems \\ in the Calculus of Variations}%
\author[G.~Bouchitt\'e, I.~Fragal\`a]{Guy Bouchitt\'e,  Ilaria Fragal\`a}
\address[Guy Bouchitt\'e]{UFR des Sciences et Techniques\\ Universit\'e de Toulon et du Var, BP 132\\ 83957 La Garde Cedex (France)
}
\email{bouchitte@univ-tln.fr}

\address[Ilaria Fragal\`a]{
Dipartimento di Matematica \\ Politecnico  di Milano \\ 
Piazza Leonardo da Vinci, 32 \\
20133 Milano (Italy)
}
\email{ilaria.fragala@polimi.it}

\keywords{}
\subjclass[2010]{  }

\date{\today}

\begin{abstract}  We present a new duality theory for non-convex variational problems, under possibly mixed Dirichlet and Neumann boundary conditions. The dual problem reads nicely as a linear programming problem, and  our main result states that there is no duality gap. Further, we provide necessary and sufficient optimality conditions, and we show that our duality principle can be reformulated as a min-max result which is quite useful for numerical implementations. As an example, we illustrate the application of our method to a celebrated free boundary problem.  The results were announced in \cite{BoFr}. 
\end{abstract} 

\maketitle

\centerline{\maius Plan of the paper}

1. Introduction

2. Setting of the primal problem

3. The duality principle 

\hskip .5 cm 3.1. Heuristic genesis

\hskip .5 cm 3.2. The admissible fields

\hskip .5 cm 3.3. The dual problem

4. Convexification recipe

\hskip .5 cm 4.1. Construction of the convex extension of the primal energy

\hskip .5 cm 4.2. Integral representation of $H$

\hskip .5 cm 4.3. Generalized coarea formula

\hskip .5 cm 4.4. Proof of Theorem \ref{l:IJ}

5. Optimality conditions and min-max formulation

6. Application to a free boundary problem 

\hskip .5 cm 6.1. Description of the problem

\hskip .5 cm 6.2. Numerical algorithms 

\hskip .5 cm 6.3. Some simulations in case $N=1$

\hskip .5 cm 6.4. Some simulations in case $N=2$

7. Completion of the proofs

\section{Introduction}\label{secintro}

A central issue of Convex Analysis is the development of a duality theory: this allows to associate with an initial convex variational problem a dual problem which has the same extremal value and in many cases is easier to solve; moreover,  solutions to both primal and dual problem can be nicely characterized through necessary and sufficient optimality conditions. This is by now a very classical road, which in the last decades has found applications in 
different areas, such as mechanics, optimal control, economics, mass transportation, and many more. 
In fact, one of the reasons is that the duality approach enables to set up very stable and efficient approximation schemes. 
We refer the reader to the reference monograph \cite{EkTe} for the theoretical framework (see also \cite{Bo}), and to \cite{BV, N} for more recent surveys including applications and numerical algorithms.  

Unfortunately, such theory  completely breaks down as soon as some nonconvexity appears in the optimization  problem under study. In particular, this drawback is often met in Calculus of Variations, where even very classical problems involve non-convex energy costs. As no systematical tool is available to characterize a global optimum, a dramatic consequence is that all currently available numerical methods loose their efficiency, because they are not able to rule out local minimizers and detect the global ones.

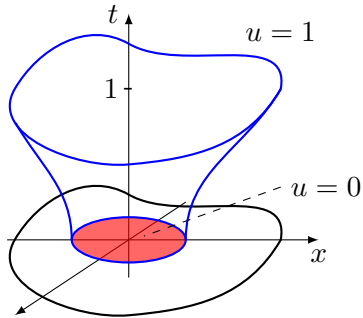
\begin{figure}[ht]
\centerline{
\begin{tikzpicture}
[>=stealth,scale=1]\draw [-latex] (0,-0.5)  -- (0,3)node [left] {$t$};
\draw [-latex](-1.6,0)-- (2.5,0) node [below] {$x$};
\draw [-latex](0.75,0.5)-- (-1.5,-1);
\draw [thick] (-1.5,0) to [out=60,in=150] (0,0.6) to [out=-30,in=80] (2,0) to [out=-120,in=5] (0,-1) to [out=-175, in=-120] (-1.5,0);
\draw [thick, blue] (-1.5,2) to [out=60,in=150] (0,2.6) to [out=-30,in=80] (2,2) to [out=-120,in=5] (0,1) to [out=-175, in=-120] (-1.5,2);
\fill [opacity=0.6, red] (0,0) ellipse (0.75cm and 0.3cm);
\draw (0,0) ellipse (0.75cm and 0.3cm)[blue, thick];
\draw [thick, blue] (-1.49,1.55) to [out=-60,in=90] (-0.75,0);
\draw [thick, blue] (2,2) to [out=-120,in=90] (0.75,0);
\draw (0,2) node [left] {$1$};
\draw (0,2)node{-} (2,2.5)node[above]{$u=1$};
\draw [dashed](2,0.7)node[right]{$u=0$} -- (0.2,0.05);
\end{tikzpicture}
}
\caption{The free boundary problem \eqref{pAC}} \label{figcaff}
\end{figure}

To have in mind a prototype situation,  let us mention for instance the free boundary problem studied in the seminal paper \cite{AlCa}:
\begin{equation}\label{pAC}
\inf \left\{\int_\Omega  \frac{1}{2} |\nabla u| ^ 2  \,dx  + \lambda
\big | \{ u >0 \} | \ :\ u \in H ^ 1  (\Omega)\, , \ u = 1 \text{ on } \partial \Omega \right\} \,, 
\end{equation}
the free boundary being the frontier of the positivity set $\{ u >0 \}$ (see Figure \ref{figcaff}). 
A huge literature about free boundaries stemmed from the existence and regularity results proved in \cite{AlCa}
(without any attempt of completeness, see  for instance \cite{ACF, CS, CF, CJK, KT, O}).  However,   these papers are mainly focused on the study of local minimizers, through the Euler-Lagrange equation and the related free boundary condition, intended in the variational or in the viscosity sense.   
To the best of our knowledge, a systematic way to evidence {\it global} minimizers for problem \eqref{pAC} is still missing. 

In this work we present a new duality theory  for non-convex variational problems, which aims at filling the lack depicted so far. In this respect, the papers \cite{AlBoDa, AtTh, Ga, PCBC1, PCBC2} should  be mentioned among the few attempts outside the convex framework. 

We consider very general minimization problems of the form
\begin{equation}\label{I}
\mathcal I := \inf \left\{\int_\Omega  f(u, \nabla u) \,dx  + \int_{\Gamma_1} \gamma(u) \, d \mathcal H^{N-1}\ :\ u\in W^{1,p}(\Omega) \ ,\ u = u_0 \hbox{ on } \Gamma_0 \right\}\,,
\end{equation}
where $\Omega$ is an open bounded domain of $\R ^N$ with a Lipschitz boundary and $(\Gamma_0, \Gamma_1)$ is a 
partition of $\partial \Omega$: $\Gamma_0$ and $\Gamma _1$ correspond respectively to the Dirichlet 
part (the datum $u _0$ is a given function in $W ^ {1,p} (\Omega)$), and 
to the Neumann part of the boundary.

The bulk integrand $f: \R \times \R ^N \to \R$ is assumed to  lower semicontinuous in both variables, and convex in $z$, 
but the key point is that it may have a {\it non-convex} dependence in $t$. 

The boundary integrand 
$\gamma$ is assumed to be Lipschitz, 
and  suitable $p$-growth conditions are imposed on $f$ and $\gamma$ to ensure the existence of a minimizer in $W ^ {1,p} (\Omega)$ (for some $p>1$).

Clearly, problem \eqref{pAC} falls into this general framework, by taking $\Gamma_0 = \partial \Omega$,  $u_0 \equiv 1$, and $f (t, z)= \frac{1}{2}  |z| ^ 2 + \chi _{(0, + \infty)}(t)$, 
where $\chi _{(0, + \infty)}$ is the characteristic function of $(0, + \infty)$. 

As a further example, one can take
$f (t, z):= \varepsilon |z| ^ 2 + W (t) - \lambda  t$, 
$W$ being a two-wells potential, $\varepsilon$  a small positive parameter, and $\lambda$ a Lagrange multiplier. 
In this case, if $\Gamma _1 = \partial \Omega$, 
problem (\ref{I}) describes the configuration of a  Cahn-Hilliard fluid in presence of a wetting term $\gamma$ on the whole of the boundary. 

For general minimization problems of the form \eqref{I}, the dual problem we propose is formulated as follows
\begin{equation}\label{I*}
{\S}:= \sup \Big \{ \int _{G_{u_0}} \sigma \, \cdot \,  \nu _{u_0} \, d \mathcal H ^N + \int_{\Gamma _1} \gamma (u_0) \, d \mathcal H ^{N-1} 
 \ :\ \sigma \in \X  \Big \}\,.
 \end{equation}
 and any optimal $\sigma$ is called a {\it calibration}, in analogy to the case of classical principle of calibration for minimal surfaces 
(see \cite{Fe, Mo, AlBoDa} and references therein). 

The class $\mathcal B$ of admissible competitors is a family of bounded divergence free vector fields $\sigma$, defined on $\Omega \times \R$,  which have a given normal trace on $\Gamma _ 1 \times \R$ and satisfy suitable convex pointwise constraints. The first integral appearing in \eqref{I*}  denotes the flux of $\sigma$ across the graph of the function $u _0$, and it is well-defined as admissible fields turn out to
admit a normal trace on any set with finite perimeter. We refer to Section \ref{secdual2} for all the details, including the precise statement of the convex constraints satisfied by the admissible fields, and its comparison with the classical dual problem in the convex case. 
 
\begin{figure}[ht]
\begin{center}
\begin{tikzpicture}
    \node[anchor=south west,inner sep=0] at (0,0) {\includegraphics[scale=0.35]{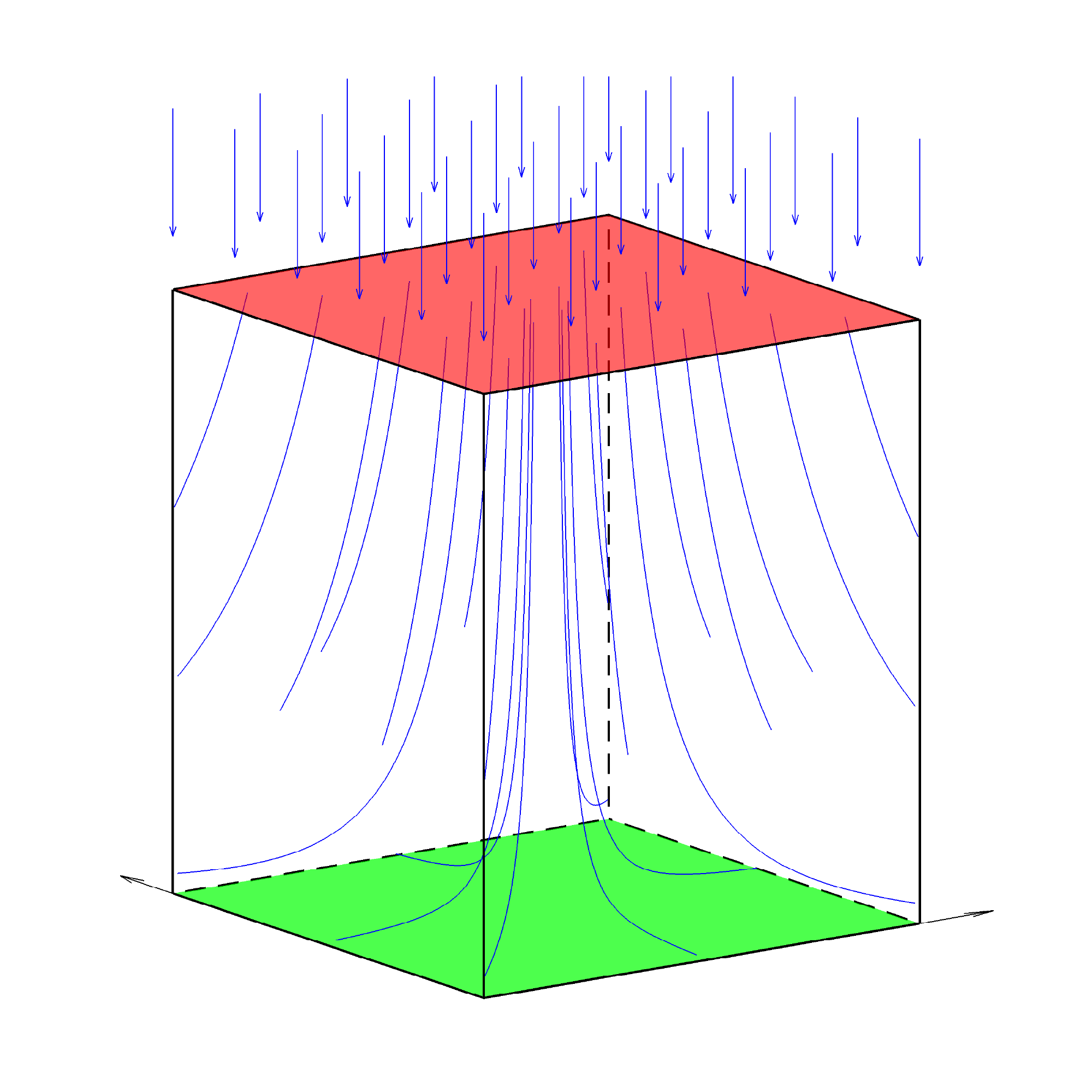}};
    \draw (3.2,4)node[right]{$\frac1{2} |\sigma^x|^2 \le \lambda + \sigma^t$};
    \draw[black!60!green](4.2,1.2)node{$\sigma^t(x,0) \ge 0$};
    \draw(6.8,1.2)node{$x_1$};
    \draw(0.55,1.4)node{$x_2$};
    \draw(0.6, 5.5)node{$t=1$};
\end{tikzpicture}
\end{center}
\caption{The optimal flow problem  \eqref{pAC*}  }\label{flow}
\end{figure}
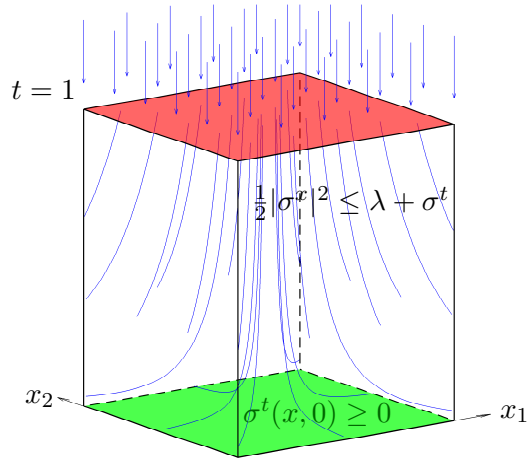

Here let us just give the complete formulation in case of problem  \eqref{pAC}, when the dual problem reads:
\begin{equation}\label{pAC*}
\mathcal I^* = \sup \Big \{ - \int _\Omega \sigma ^ t (x, 1) \, dx  \ :\ \sigma \in \mathcal B \Big \}\,. 
\end{equation}
Notice that in this case the integral on $\Gamma _1$ is missing (since $\Gamma _0 = \partial \Omega$), whereas the integral on $\Omega$ represents the flux term across the graph of the boundary datum $u _0 \equiv 1$. Namely, $\sigma ^ t$ denotes the vertical component of an element $\sigma = (\sigma ^ x , \sigma ^ t)$ lying in the admissible class $\mathcal B$, which for the problem under consideration is given by all  bounded divergence free vector field on $\Omega \times \R$ satisfying the constraints
$$\sigma ^ t (x, t) + \lambda \geq \frac{1}{2} |\sigma ^ x ( x, t) | ^ 2   \quad \text{a.e. on } \Omega \times \R \, , \qquad \sigma ^ t (x, 0) \geq 0  \text{ a.e. on } \Omega\,.
$$
Thus problem \eqref{pAC*} has a nice fluid mechanic interpretation:  it consists in maximizing the downflow 
through the top face $\Omega \times \{1\}$  
of an incompressible fluid constrained into the cylinder $\Omega \times \R$, whose speed $\sigma$ satisfies the conditions above, preventing in particular the fluid to pass across the bottom face (see Figure \ref{flow}, in which $\Omega = (0, 1) ^ 2 \subset \R ^2$).

Our main result establishes that, in the general setting sketched above and fixed more precisely in Section \ref{secdual}, there is no duality gap: the infimum $\mathcal I$ in \eqref{I} and the supremum $\mathcal I ^*$ in \eqref{I*} coincide. 
The result is stated, along with several comments, in Section \ref{secdual2} (see Theorem \ref{t:nogap}), after providing a heuristic 
description of the underlying idea, and giving all the required details about the class of admissible fields. 

The proof is quite delicate and to it is devoted most part of the paper. Here we limit ourselves to give just few hints. 
The approach we adopt  is based on the idea of reformulating the primal problem \eqref{I}  in $(N+1)$ space dimensions. More precisely, in the same spirit of what done in the paper \cite{AlBoDa} for the Mumford-Shah functional (see also \cite{Ch}), the starting point is to identify any admissible function $u: \Omega \to \R$ with the characteristic function $\1 _u$ of its subgraph. 
Then the building block of our method is a convexification recipe, which is carried over in Section \ref{recipe}. Roughly speaking, it consists in embedding the class $\mathcal A$ of competitors for the primal problem \eqref{I} into an enlarged class $\widehat {\mathcal A}$ of functions $v$ defined on $\Omega \times \R$ (via the identification $u \mapsto \1 _u$), and in constructing a {\it convex} functional $\widehat E$, which extends the primal energy 
$E (u):= \int _\Omega f (u, \nabla u) \, dx+ \int _{\Gamma _1} \gamma (u) \, d \mathcal H ^ {N-1}$ to the class $\widehat {\mathcal A}$. The key intermediate result (see Theorem \ref{l:IJ}) states that the infimum of the convex functional $\widehat E$  over the class  $\widehat {\mathcal A}$ coincides with $\mathcal I$, and that the solutions to the two problems are closely related to each other.  To establish such result, we exploit as a crucial ingredient a new very general coarea type formula (see Theorem \ref{theo}).

The completion of the proof of Theorem \ref{t:nogap}  
is postponed in Section \ref{secproofs} (since this last part is not needed for the comprehension of the contents of Sections \ref{secminmax} and \ref{secexamples}).    It is obtained essentially by using convex duality in $(N+1)$ space dimensions, in synergy with several ad-hoc arguments, driven from convex analysis and geometric measure theory, needed to handle the involved functions and fields.

The companion results of our duality theory are presented in Section \ref{secminmax}: in Theorem \ref{t:calib} we show that solutions to the primal and to the dual problem can be characterized through an equality holding on the graph of an optimal function $\overline u$, and in Corollary \ref{p:september} we give a practical way to check such condition in concrete situations; in Theorem \ref{minmax} we reformulate our duality principle under the form of a min-max result, and a variant which is conceived especially for numerical purposes is added in Proposition \ref{c:minmax}.   

In Section \ref{secexamples} we exemplify the application of our method to problem \eqref{pAC}.

To conclude, let us stress that this paper aims to give a breakthrough by settling the bases of the non-convex duality theory, 
but of course it cannot contain the many developments which are expected and which will be studied in forthcoming works. 

In particular, the existence of a 
solution to the dual problem, that we call a {\it calibration},  is a major issue.
In the forthcoming paper \cite{BoFrPh}, by using rearrangement techniques for integrals with  non-constant densities, we are going to provide an existence result for problems with linear growth (for which a variant of Theorem \ref{t:nogap} can be established).
Moreover, the numerical results given in Section \ref{secexamples} will be detailed and expanded in \cite{BoGaPh}. 

As further open problems and possible generalizations, let us mention that 
our duality principle may be easily extended to the case when $f$  and $\gamma$ depend also on the spatial variable $x$. On the other hand,  possible adaptations of the same idea to variational integrals involving the Hessian of $u$ are not straightforward and deserve further investigation. 
Finally, our results open the innovative perspective of studying the stability of minimizers 
of non-convex functionals by computing their shape derivatives (in fact, our duality result should allow to extend successfully to the non-convex setting the approach recently proposed in \cite{BoFrLu, BoFrLu2}).

\smallskip
{\bf Acknowledgments.} We acknowledge the financial support the University of Toulon, Politecnico di Milano, the University of Pavia, and the Italian institutions MIUR and INDAM, which helped the accomplishment of this work through PRIN and GNAMPA projects. 
We are very grateful to C\'edric Galusinski and Minh Phan for handling the numerical simulations presented in Section \ref{secexamples}.

\section{Setting of the primal problem}\label{secdual}

Let $\Omega$ be an open bounded domain of $\R ^N$, and let $(\Gamma_0, \Gamma_1)$ be a partition of  $\partial \Omega$. 

We consider as primal problem the non-convex infimum problem 
\begin{equation}\label{defI}\qquad (\mathcal P)    \qquad \qquad \mathcal I:= \inf  \Big \{ E (u) :\ u \in \mathcal A \Big \}  \,, \qquad \qquad \qquad \qquad \qquad
\end{equation}

where the energy cost is of the form
\begin{equation}\label{f:E}
E (u):= \int _{\Omega} f (u, \nabla u ) \, dx + \int _{\Gamma _1} \gamma (u) \, d \mathcal H ^ {N-1} \,,
\end{equation}

and the class of admissible functions is given by 
\begin{equation}\label{f:A}
\mathcal A := \Big \{ u \in W^ {1, p} (\Omega) \ :\ u = u _0 \hbox{ on } \Gamma_0 \Big \}\, ,
\end{equation}

being $u_0$ a fixed element in $W ^ {1, p} (\Omega)$.

\medskip
We work  under  the setting of hypotheses listed hereafter. 

\medskip
{\bf Standing assumptions}:

\medskip
$\bullet$  The boundary $\partial \Omega$ is Lipschitz with unit outer normal $\nu _\Omega$.


\medskip
$\bullet$ The integrand $f = f (t, z)$
is a function $f : \R \times \R ^N \to ( - \infty, + \infty]$ sastisfying:
\begin{eqnarray}
& \forall t \in \R \, , \ z \mapsto f(t, z) \text{ is convex; } \label{H2} \\
&  (t, z) \mapsto f(t, z)\text{ is lower semicontinuous on }\R \times \R ^N; 
\label{H1} \\
&  \forall  (t, z) \in \R \times \R ^N,  \     f (t, z) \geq  \alpha |z |^ p - r (t) \, ,   \label{H3}
\end{eqnarray}

where $p\in(1, + \infty)$, 
$\alpha$ is a positive constant, and $r = r (t)$ is a Borel function such that   
\begin{equation}\label{hypr}
0 \leq r(t) \le C  \ \ \text{ for some } C>0 \,.
\end{equation}


\medskip
$\bullet$  There exists a {\it Lebesgue negligible} set of $D\subset \R$ such that, for every $z \in \R ^N$,  the map $t \mapsto f(t, z)$  is upper semicontinuous on $\R\setminus D$, namely
\begin{equation}\label{H1bis}
f (t, z) \geq \limsup_{s \to t} f (s, z)  \qquad \forall z \in \R ^N \, , \ \forall t \in \R \setminus D \,.
\end{equation}

\medskip
$\bullet$ $\gamma: \R \to \R $ is a  Lipschitz  function such that  $\gamma (0) = 0$ and 
\begin{equation}\label{f:condinf}
\begin{cases}
\displaystyle \inf _{t \in \R}    \gamma (t) > - \infty   & \text{ if } \Gamma _0 \neq \emptyset
\\ \noalign{\bigskip} 
\displaystyle \liminf _{|t|  \to + \infty } \frac{ \gamma (t)}{|t|} >0    & \text{ if } \Gamma _0 = \emptyset\,,
\end{cases}
\end{equation}

\medskip
$\bullet$ The set $\Big \{ u \in \mathcal A \ :\ E (u) < +\infty \Big \}$ is not empty.

\medskip
\begin{remark}\label{r:stand}{\rm (i) We emphasize that the function $f$ is {\it not} assumed to be convex in $t$.

(ii)  
We point out that, by taking $\Gamma_1=\partial\Omega$ and $\gamma=0$, we can handle homogeneous Neumann boundary conditions. Notice also that the  condition $\gamma (0) = 0$ is not restrictive up to adding a constant.  }

(iii)  Allowing   a nonempty discontinuity set $D$ of vanishing Lebesgue measure for the map $t \mapsto f (t, z)$ (according to \eqref{H1bis}) is quite important in order to make our duality method applicable in case of free boundary problems, {\it cf.} Section \ref{secexamples}.

(iv) The boundedness of $r$ is a technical condition which  will be exploited mainly in the proof of Lemma \ref{l:B0}.

(v) One of the main roles of the growth conditions \eqref{H3} and \eqref{f:condinf} imposed respectively on $f$ and $\gamma$ is to  ensure the well-posedness of the primal problem, as stated in the next result. 

(vi) We stress that, for the validity of Proposition \ref{l:existence}, it is important to have chosen $p>1$ in \eqref{H3}, since for $p=1$ the primal problem  may fail to admit a solution. 
The main reason
is that in such case the energy $E$  is no longer lower semicontinuous
(whereas coercivity still holds, as it is easy to see by inspection of the proof below). Thus one needs to relax the energy $E$ in $BV (\Omega)$ (see \cite{Da}), which is made extremely delicate by 
the presence of the boundary integral in \eqref{f:E}, in particular when $\partial \Omega$ exhibits corners (see \cite{suquet}). 

On the other hand, with minor modifications in the proof,  our duality Theorem \ref{t:nogap}  remains true also in the case $p=1$ (provided $\Gamma_0 = \partial \Omega$), and this is precisely the setting in which it seems easier to obtain the existence of a solution  for the dual problem. 
An existence result for the dual problem in the framework of nonconvex functionals with linear growth under Dirichlet boundary conditions will be the topic of a forthcoming paper. 
\end{remark}

\begin{proposition}{\rm (well-posedness of the primal problem)}  \label{l:existence} 
The infimum $\mathcal I$ in \eqref{defI}  is finite and attained. 
\end{proposition}
 \proof   Since we assumed that the class $\mathcal A$ of admissible competitors contains some element $u$ of finite energy, we may apply the direct method of the Calculus of Variations. Thus we are reduced to showing that, under the standing assumptions,  the energy $E$ defined in \eqref{f:E} is both lower semicontinuous and coercive  respect to the weak topology of $W ^ {1,p} (\Omega)$. 
 
 The weak lower semicontinuity  of the first addendum of the functional $E$ follows  
 well-known results of weak-strong convergence (see for instance \cite[Chapter 4]{Bu}), 
 which can be applied in particular thanks to the growth conditions \eqref{H3}. 
 
The weak lower semicontinuity of the second addendum follows as a consequence of the compact embedding of $W ^ {1,p} (\Omega)$ into $L ^ p (\partial \Omega)$, by applying Fatou's lemma. 

We then focus attention on the coercivity property. We claim that there exists positive constants $C_1$, $C_2$ such that
\begin{equation}\label{f:coerc}
E (u) \geq C _1 \| u \| _{W ^ {1,p} (\Omega)} - C _2 \,.
\end{equation}
In case $\Gamma _1 = \emptyset$, the coercivity follows immediately from the lower bound in \eqref{H3}, taking into account that $r(t)$ 
satisfies \eqref{hypr}.
 
In case $\Gamma _1 \neq \emptyset$,  we further distinguish the cases $\Gamma _0 \neq \emptyset$ and $\Gamma _0 = \emptyset$. 
If $\Gamma _0 \neq \emptyset$,  the coercivity follows again from the lower bound in \eqref{H3}, taking into account that $r (t)$ satisfies \eqref{hypr}.
If $\Gamma _0  = \emptyset$, the lower bound in \eqref{H3} tells us merely that $u _n$ are bounded in $W ^ {1,p} (\Omega)$ modulo 
constants, but
by invoking the second condition in \eqref{f:condinf}, we obtain that the  boundary traces of $u _n$ are bounded in $L ^ 1 (\partial \Omega)$, and hence the constants are bounded.

\qed

\section{The duality principle}\label{secdual2}

In this section we present our new duality principle:

\medskip
\begin{itemize}
\item[--] in Subsection \ref{genesis} we provide an intuitive presentation of the underlying idea; 

\smallskip
\item[--] in Subsection \ref{af} we introduce the class of admissible fields in the dual problem;

\smallskip
\item[--] in Subsection \ref{dp} we state the result (see Theorem \ref{t:nogap}), along with some basic remarks.  
\end{itemize}

\subsection{Heuristic genesis.}\label{genesis} The original idea, already exploited in \cite{AlBoDa} for free-discontinuity problems,  relies on geometric measure theory and stems from the so-called calibration method for minimal surfaces  (see \cite{Mo,Fe}).
It consists in considering a suitable convex set $\K$ of vector fields $\sigma=(\sigma^x,\sigma^t): \Omega\times \R \to \R^{N+1} $ satisfying the following requirement: 
\begin{equation} \label{calib=}
\int_\Omega f (u, \nabla u) \, dx  \ =\ \sup_{\sigma\in\K} \, \int_{G_u}  \sigma \cdot \nu_u \, d \mathcal H^N  \qquad \forall \, u\in W^{1, p}(\Omega)\,. 
\end{equation}
The integral at the right hand side of \eqref{calib=}, which  is well defined provided $\sigma$ is regular enough,  represents the flux of
 $\sigma$ across the graph $G_u$ of $u$, seen as a $N$-dimensional rectifiable subset of $\R^{N+1}$  and endowed with 
oriented unit normal 
\begin{equation}\label{f:nu}
\nu_u= \frac{ (\nabla u , -1) } {\sqrt{1 + |\nabla u | ^ 2}}\,.
\end{equation}

\begin{figure}
\includegraphics[scale=0.4]{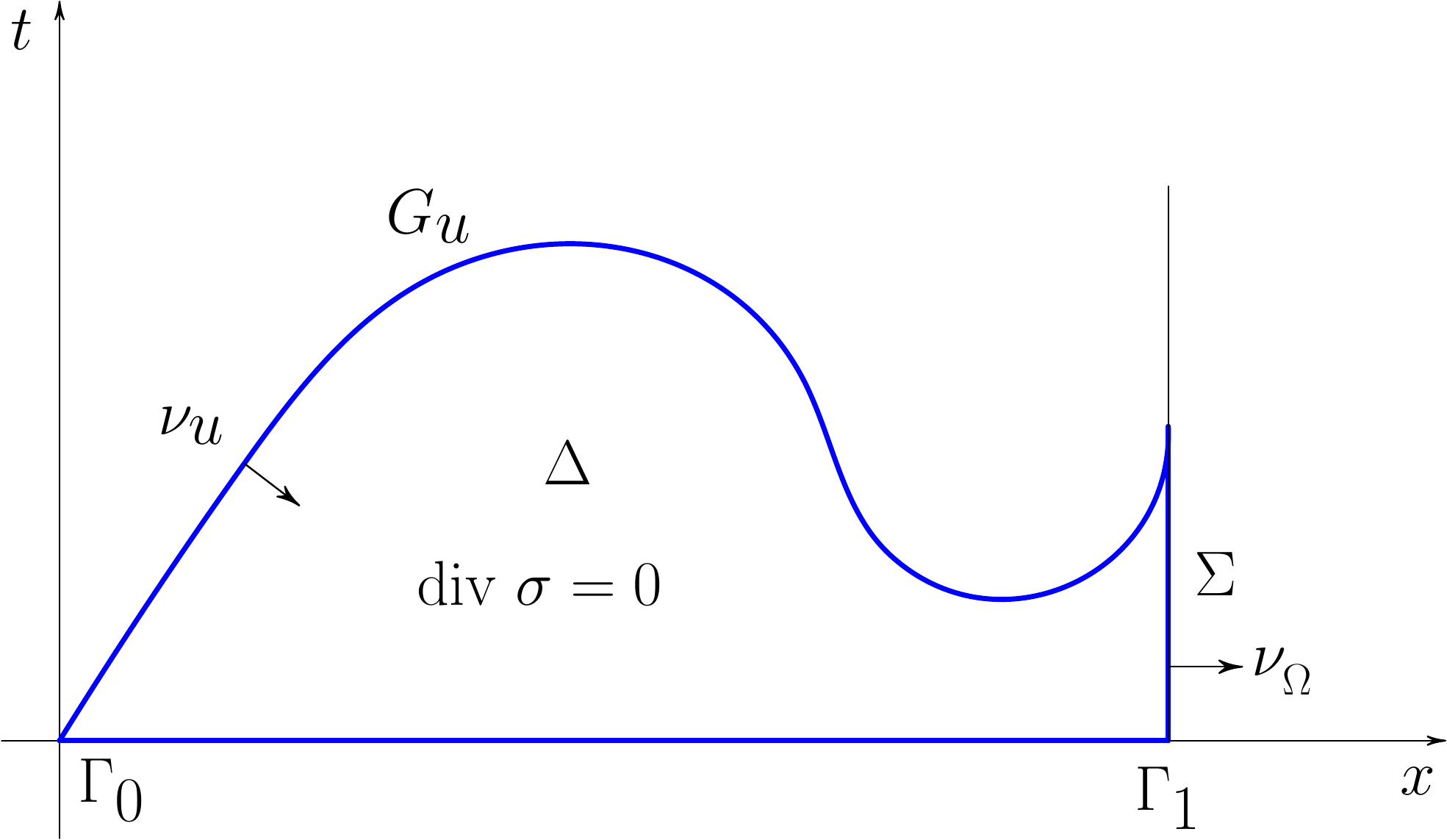}
\label{figdelta}
\end{figure}

Given  a function $u$ in $W^{1, p} (\Omega)$  such that $u = u _0$ on $\Gamma_0$, we denote by $\Delta$ the subset of  $\Omega \times \R$ lying between $G_{u_0}$ and $G_u$, and by $\Sigma \subset \Gamma _ 1 \times \R$  the ``lateral part'' of $\partial \Delta$, namely the set of points $(x, t)$ with $x \in \Gamma_1$, and $t$ between $u _0 (x) $ and $u (x)$.  In case $N=1$, taking $u _0 = 0$ and $u \geq 0$,  the region $\Delta$ is represented in Figure \ref{figdelta}.   

Let now  $\sigma$ be a  smooth element belonging to a class $\K$ verifying \eqref{calib=}, and assume that $\sigma$ satisfies the additional conditions
\begin{equation}\label{ac}
\div \sigma = 0 \ \hbox { in } \Omega \times \R \qquad \hbox{ and } \qquad \sigma \cdot \nu _\Omega = - \gamma '  \ \hbox{ on } \Gamma_1  \times \R\,. 
\end{equation}
By applying the divergence theorem on the region $\Delta$, we obtain:
$$\begin{array}{ll}\displaystyle
\int_{G_u}  \sigma \cdot \nu_u \, d\mathcal H^N   - \int_{G_{u_0}}  \sigma \cdot \nu_{u_0} \, d\mathcal H^N &\displaystyle= \int _{\Sigma} {\rm sign} (u - u _0)\   \sigma \cdot \nu _\Omega  \, d \mathcal H ^{N}  
\\  \noalign{\bigskip} &  \displaystyle = - \int _{\Sigma} {\rm sign} (u(x) - u _0(x))\   \gamma '(t)  \, d \mathcal H ^{N-1}(x) \, dt 
\\  \noalign{\bigskip} &  \displaystyle =  \int _{\Gamma_1}\big (  \gamma (u_0) - \gamma (u ) \big ) \, d \mathcal H ^{N-1} \,.    
\end{array}
$$

In view of \eqref{calib=}, and recalling the definition \eqref{f:E} of the energy $E$, we deduce that 
$$E (u) \geq  \int_{G_{u_0}}  \sigma \cdot \nu_{u_0} \, d\mathcal H^N +  \int _{\Gamma_1} \gamma (u_0)  \, d \mathcal H ^{N-1} \,.$$

It is then natural to optimize the above inequality by considering the   linear programming problem
\begin{equation}\label{lp}
\sup \Big \{ \int_{G_{u_0}}  \sigma \cdot \nu_{u_0} \, d\mathcal H^N +  \int _{\Gamma_1} \gamma (u_0)  \, d \mathcal H ^{N-1} \ : \  
\sigma \in \K  
\text{ satisfying }   (\ref{ac})
\Big \} \,. 
\end{equation}

Clearly from the above discussion the supremum in \eqref{lp} turns out to be bounded from above by  the infimum $\mathcal I$ of the primal problem. 
We have thus found a linear programming problem which is a good candidate for being the dual problem. To elect it as such, we have to complete the plan, by choosing 
 $\K$ so that the equality \eqref{calib=} holds and
the  supremum in \eqref{lp} equals $\mathcal I$. 

Let us now focus our attention on the  construction of the class $\K$,  by giving some heuristic arguments (the rigorous definition is postponed to Section \ref{af} below).

Assume that  $\sigma= (\sigma ^ x, \sigma ^t) \in \mathcal C ^ 1 (\Omega \times \R ; \R ^ {N+1})$ satisfies  the pointwise inequality
\begin{equation}\label{stima}
\sigma ^ t (x, t) \geq f ^*_z (t, \sigma ^x(x, t)) \qquad \forall  (x, t) \in \Omega \times \R \,, 
\end{equation}

where $f ^*_z$ denotes the Fenchel conjugate of $f$ with respect to $z$:
$$f ^ * _z (t, z^*):= \sup _{z \in \R ^N } \big [ z \cdot z ^* - f (t, z)  \big ] \,.$$
By using \eqref{stima} on the graph of $u$ and the Fenchel inequality, we obtain 
$$
\begin{array}{ll}
\displaystyle \int_\Omega  f(u, \nabla u) \,dx & \displaystyle \geq   \int _\Omega \big [ f ^ * _z (u(x), \sigma ^x(x, u(x))) + f (x, \nabla u(x)) - \sigma ^ t (x, u(x))  \big ] \, dx 
 \\  \noalign{\medskip}
& \displaystyle \geq  \int_\Omega  \big[ \sigma^x(x,u(x))\cdot \nabla u- \sigma^t(x,u(x))\big ] \,dx 
 =  
\int_{G_u}  \sigma \cdot \nu_u \, d \mathcal H^{N}   
\,.    
\end{array}
$$
The above inequality turns out to optimal: actually, as it will be shown later, 
if $\K$ is chosen as the class of fields in  $\mathcal C ^ 1 ( \Omega \times \R ; \R ^ {N+1})$ satisfying \eqref{stima}, not only the equality
\eqref{calib=} holds true, but in addition the supremum in \eqref{lp} equals $\mathcal I$.  

However, the class of competitors we are going to choose in our dual problem has also to be large enough in order to allow the existence of optimal fields. In this respect, it will be clear from the examples considered in Section \ref{secexamples} that one cannot expect optimal fields to be $\mathcal C ^1$ regular, and not even to be continuous fields which satisfy the inequality \eqref{stima} pointwise at {\it every} $(x, t)$ in $\Omega \times \R$. 

We are thus led to relax condition \eqref{stima} and to  work with fields which are less regular, but still admit a mathematically meaninfgul notion of flux and normal trace.

\subsection{The admissible fields}\label{af}

We consider the space 
\begin{equation}\label{X1}
X_1 (\Omega \times \R) := \Big \{ \sigma \in L ^ \infty (\Omega \times \R; \R ^ {N+1}) \ :\  \div \sigma \in L ^ 1 (\Omega \times \R) \Big \}\,, 
\end{equation}
where the divergence is intended in distributional sense. 

For any $\sigma \in  X_1 (\Omega \times \R)$, a notion of weak normal trace can be defined as follows.  
Given an open set $A \subset \Omega \times \R$ with  Lipschitz boundary and unit outer normal $\nu _A$, 
there exists a unique function $\sigma\cdot \nu _A \in L^\infty(\partial \Omega)$ such that
\begin{equation}\label{gen_trace}
\int_{\partial A}(\sigma\cdot \nu _A) \,\varphi\,d{\mathcal  H^{N}}=\int_A \big (\sigma\cdot \nabla \varphi +\varphi\,\div\sigma\big ) \,d x \qquad \forall \varphi \in \mathcal C ^\infty _0 (\Omega \times \R)\ .
\end{equation}
The same assertion remains true when $A$ is merely a Lebesgue measurable set with finite perimeter, provided $\partial A$ is intended as  the reduced boundary of $A$, and $\nu _A$ as the measure theoretic unit normal vector defined $\mathcal H ^N$-a.e.\ on $\partial A$. 

In particular, for any field $\sigma \in  X_1 (\Omega \times \R)$ and any function $u \in W ^ {1, p} (\Omega)$, the flux integral
$$\int_{G_u} \sigma \cdot \nu _u \, d \mathcal H ^N$$ 
is well-defined according to \eqref{gen_trace} (precisely,  by taking 
as a set $A$ the subgraph of $u$, we have $\nu _A = - \nu _u$, with $\nu _u$ given by \eqref{f:nu}). 
 
For later use, let us notice that, as \eqref{gen_trace}  can be extended to all $\varphi\in L^\infty(A) \cap W^{1,1}(A)$,
a duality argument easily yields the following equality
\begin{equation}\label{gen_claim}
\begin{array}{ll}
& \displaystyle \Big \{ (-\div q, q\cdot \nu_A) : q\in X_1(A) \Big \}  = \\ 
\noalign{\medskip}
&
\displaystyle \Big \{ (f, g) \in L ^ 1 (A) \times L ^ \infty (\partial A) : \int_{A} \!f   dx + \int _{\partial A} \! g  \, d \mathcal H ^N \, = 0 \Big \} \,.
\end{array}
\end{equation} 
 
We refer to \cite{A, CF} for more details  on these topics (see also Section \ref{proof1}, where we shall need to exploit a generalized version of the Gauss-Green Theorem involving BV functions).

%
%
\begin{definition}\label{d:B} (i) We set $\K$ the class of fields $\sigma = (\sigma ^ x, \sigma ^ t) \in   X_1 (\Omega \times \R)$ such that 
\begin{eqnarray}
& \sigma ^ t (x, t) \geq f ^*_z (t, \sigma ^x(x, t)) \ \text{ for } \mathcal L ^ {N+1}\text{-a.e.}\ (x, t) \in \Omega \times \R  &\label{s3}
\\
\noalign{\smallskip}
& \sigma ^ t (x, t) \geq - f (t, 0) \quad \forall\, t \in D  \text{ and for } \mathcal L ^ {N}\text{-a.e.}\ x \in \Omega\, ,   & \label{s4}
\end{eqnarray}
where $D$ is the Lebesgue negligible set introduced in the standing assumption \eqref{H1bis}. 

\medskip
(ii) We denote by $\X$ the class of fields $\sigma \in \mathcal K$ satisfying the following two conditions:
\begin{eqnarray} 
& \div \sigma = 0   \ \hbox{ in } \Omega  \times \R &\label{f:div}
\\
\noalign{\smallskip}
&  \sigma ^x \cdot \nu _\Omega = - \gamma '  \ \hbox{ on } \Gamma_1  \times \R\,, \label{f:bc}
\end{eqnarray}
where $\sigma^x \cdot \nu _\Omega$ is meant as the weak normal trace of $\sigma$ on $\partial (\Omega \times \R)$
(as $\nu_{\Omega\times \R}  = (\nu _\Omega, 0)$).

\end{definition}

\begin{remark} \label{r:june}
Few comments are in order about condition \eqref{s4}, which did not appear in our previous heuristic discussion. First we observe that, for every fixed $t \in D$, $\sigma ^ t (\cdot, t) $ makes sense as the weak normal trace of $\sigma$ on $\Omega \times \{ t \}$ according to \eqref{gen_trace}.  
The role of 
condition \eqref{s4} is to make the class of admissible fields sensitive to the possible discontinuities of the integrand $f$. 
In this respect, the almost everywhere inequality \eqref{s3}  alone would be too weak, since it is independent from the behaviour of $f$ on sets of vanishing measure:  for instance, condition \eqref{s3} reads exactly the same  in the two cases when  $f (t, z) = \frac{1}{2} z ^ 2 +  \chi _{\{t \neq 0 \}}$ or  $f (t, z) = \frac{1}{2} z ^ 2+1$ (since the discontinuity set $\{t = 0 \}$ is $ \mathcal L ^ {N+1}$-negligible).
Finally, let us mention that the inequality \eqref{s4} is actually satisfied also on the complement of $D$. Namely we shall see later on that, for any $\sigma \in \X$, there holds
$\sigma ^ t (x, t) \geq f (t, 0)$ for {\it every} $t \in \R$ and $\mathcal L ^N$-a.e.\ $x \in \Omega$ ({\it cf.} Remark \ref{r:easyine2}). 
\end{remark}

\begin{remark} Conditions \eqref{s3}-\eqref{s4} can be rephrased as
\begin{eqnarray}
& \sigma (x, t) \in K (t) \quad \text{ for } \mathcal L ^ {N+1} \hbox{-a.e. } (x, t) \in \Omega \times \R\,,  & \label{s3'}
\\ \noalign{\smallskip}
& \sigma ^ t (x, t) \in \Pi_{N+1} [K(t)] \quad \forall \, t \in D \text{ and for }\mathcal L ^ {N} \hbox{-a.e. } x \in \Omega\, ,  & \label{s4'}
\end{eqnarray}
where $K(t)$ is the convex subset  of $\R ^ {N+1}$ given for every $t \in \R$ by
\begin{equation}\label{f:C}
K(t):= \big \{ q = ( q ^ x, q ^ t) \ :\ q ^ t \geq f ^*_z (t, q ^x) \big \}\,,  
\end{equation}
and  $\Pi_{N+1} [ \cdot ]$ denotes the projection on the last component of $\R ^ {N+1}$ (that is, on the space spanned by $e _{N+1} = (0, 1)$).  In particular, the equivalence between \eqref{s4} and \eqref{s4'} follows from the identity $\inf\limits_{z^*} f ^*_z (t, z^*) = -f (t, 0)$. 
\end{remark}
\medskip

\subsection{The dual problem}\label{dp}

\medskip Recall that $u _0 \in W^ {1, p} (\Omega) $ is the prescribed trace on the Dirichlet part $\Gamma_0$  of the boundary 
({\it cf.}\ \eqref{f:A}), that  $\gamma$ is the energy density on the Neumann  part $\Gamma_1$ of the boundary ({\it cf.}\ \eqref{f:E}), and that, 
 for every field $\sigma$ belonging to the class $\X$ 
introduced in Definition \ref{d:B}, the flux across the graph of $u _0$ is well-defined
as explained in Section \ref{af}. 
We set 
\begin{equation}\label{defI*}
({\mathcal P}^*)
 \qquad  {\S}:= \sup \Big \{ \int _{G_{u_0}} \sigma \, \cdot \,  \nu _{u_0} \, d \mathcal H ^N + \int_{\Gamma _1} \gamma (u_0) \, d \mathcal H ^{N-1} 
 \ :\ \sigma \in \X  \Big \}\,.
 \end{equation}
\medskip

The core of our duality theory is the following

\begin{theorem}[duality principle]\label{t:nogap}
The extrema of the primal and dual problems defined respectively in \eqref{defI} and \eqref{defI*} coincide: 
\begin{equation}\label{dualprinciple}
{\mathcal I} = {\S}\,.
\end{equation}
\end{theorem}

Several comments are listed in the next remarks. 

\begin{remark}{\rm 
In the pure Neumann case when $\Gamma _0 = \emptyset$ (so that the boundary datum $u_0$ is not defined), definition \eqref{defI*} must be intended as if $u _0 = 0$, namely 
 $\mathcal I ^*$  can be reformulated as  ({\it cf.}\ \cite{BoFr})
$$\mathcal I ^ * = \sup \Big \{ -\int _\Omega \sigma ^ t (x, 0) \, dx \ :\ \sigma \in \X   \Big \}\,. $$
}
 \end{remark}

\begin{remark}\label{r:bounded}{\rm In many cases, when the boundary datum $u _0$ is a bounded function, there exist a priori lower or upper bounds for the minimizers of the primal problem $({\mathcal P})$,
so that the infimum value ${\mathcal I}$ is unchanged if we impose $u$ to take values in a suitable closed interval $[m, M]$ of the real line. We are thus
led to consider the variant of the primal problem \eqref{defI} where the class of admissible functions is changed into 
\begin{equation}\label{AmM}
\mathcal A (m, M) := \Big \{ u \in W ^ {1, p} (\Omega; [m, M]) \ :\ u = u _0 \text{ on } \Gamma _0 \Big \}
\end{equation}
In this case, our duality result continues to hold (with a simpler proof, see Proposition \ref{c:minmax}), provided the admissible fields  in the dual problem $({\mathcal P}^*)$ 
are taken in the class $\mathcal B (m, M)$ of elements $\sigma \in X _ 1 (\Omega \times (m, M))$  satisfying:
\begin{eqnarray}
& \sigma ^ t (x, t) \geq f ^*_z (t, \sigma ^x(x, t)) \ \text{ for } \mathcal L ^ {N+1}\text{-a.e.}\ (x, t) \in \Omega \times (m , M)  &\label{s31}
\\
\noalign{\smallskip}
& \sigma ^ t (x, t) \geq - f (t, 0) \quad \forall\, t \in D \cup \{m, M \} \text{ and for } \mathcal L ^ {N}\text{-a.e.}\ x \in \Omega & \label{s41}
\\
\noalign{\smallskip}
& \div \sigma = 0   \ \hbox{ in } \Omega  \times (m , M) &\label{f:div1}
\\
\noalign{\smallskip}
&  \sigma ^x \cdot \nu _\Omega = - \gamma '  \ \hbox{ on } \Gamma_1  \times (m , M)\,.\label{f:bc1}
\end{eqnarray}
This reduction of the dual problem to a bounded set will be of course crucial in the implementation of efficient  algorithms for the numerical approximation of its solutions.}
\end{remark}

\begin{remark}{\rm In general the solution to the dual problem $({\mathcal P}^*)$ is not unique (see Section \ref{secexamples}). However,  if
the 
infimum of $({\mathcal P})$ is reached in $\mathcal A (m, M)$ 
and the supremum of $({\mathcal P}^*)$ is reached  in $\mathcal B (m, M)$, then a unique solution to $({\mathcal P}^*)$ can be selected  by considering the Tikhonov regularization 
$$({\mathcal P_\e}^*) \qquad \qquad  \sup \Big \{\int _{G_{u_0}} \sigma \, \cdot \,  \nu _{u_0} \, d \mathcal H ^N + \int_{\Gamma _1} \gamma (u_0)  - \ \e  \int_{\Omega\times (m , M)} |\sigma|^{2} \, dx \ :\ \sigma \in \mathcal B (m, M)   \Big \}\,. \qquad \qquad\qquad \qquad$$
As $\e \to 0$, we are led to the solution of minimal $L ^{2}$-norm. }
\end{remark}

\begin{remark}\label{r:convexdual}{\rm
In case the integrand $f$ is convex in $(t, z)$, the inequality $\mathcal I ^* \geq \mathcal I$ (which is the most delicate part in the proof of Theorem \ref{t:nogap}) is a straightforward consequence of classical duality theory. 
To see this, consider  vector fields   of the form $\sigma (x, t)=  (\eta (x), a (x) - t \div \eta (x))$. For such fields, the inequality  
$\sigma ^ t(x, t)\geq f ^ *_z ( t , \sigma ^ x(x, t) ) $ is satisfied if and only if
$$a (x) \geq \sup _t \big \{ t \div \eta (x) + f^ * _z (t, \eta)  \big \} = \sup_{ (z, t)} \big \{ t \div \eta + 
z \cdot \eta - f (t, z) \big \}= f ^ * (\div \eta, \eta) \,.$$
We deduce that $ \X $ contains the class  $\Theta$ given by fields of the form  $\sigma (x, t)=  (\eta (x), a (x) - t \div \eta (x))$, with 
$\eta \in \mathcal C ^1 (\overline \Omega; \R ^N)$,  $\eta \cdot \nu _\Omega = - \gamma'$  on  $\Gamma_1$, and 
$a \in \mathcal C ^0  (\Omega)$, $a (x)\geq f ^ * (\div \eta, \eta)$ in $\Omega \times \R$. Therefore, 
\begin{equation}\label{i:convex} \begin{array}{ll}{\S} & \displaystyle \geq\sup \Big \{ \int _\Omega  - \sigma ^ t (x, 0) \, dx \ :\ \sigma \in \Theta    \Big \}  \\
\noalign{\medskip}
& \displaystyle = \sup \Big \{ \int _\Omega - f ^ * (\div \eta, \eta) \, dx \ :\ \eta \in \mathcal C ^1 (\overline \Omega; \R ^N) \, ,\  \eta \cdot \nu _\Omega = - \gamma' \text{ on } \Gamma_1  \Big \}\,.
\end{array}
\end{equation}
The variational problem in the last line is the classical dual problem of $(\mathcal P )$, and its supremum coincides with $\mathcal I$  by standard convex duality (see for instance \cite{Bo, EkTe}).}
\end{remark}

\begin{remark}\label{r:value} {In case $N=1$, when  the variational problem $({\mathcal P})$ is settled on an interval $(0, h)$ of the real line, 
every competitor $\sigma$ in the dual problem is a bounded divergence free vector field on $(0, h ) \times \R$, so that it can been written under the form $ \sigma = (\partial _ t w, - \partial _ x w)$, for some  function $w\in {\rm Lip} ((0, h) \times \R )$. 
For instance, in the pure Dirichlet case $\Gamma _0 = \{0, h \}$ with boundary conditions 
$u (0 ) = u (h) = c$, when the primal problem reads
\begin{equation}\label{1d}
\mathcal I = \inf \Big \{ \int _0 ^ h f (u, u' ) \, dt \ :\ u \in H ^ 1 (0, h)\, , \ u (0) = u (h) = c \Big \}\,, 
\end{equation}
the dual problem \eqref{defI*} written in terms of rotated gradients becomes:
\begin{equation}\label{MK}
\begin{array}{ll} \S= \sup \Big \{ w (h, c) -w (0, c) \, :\, & w \in {\rm Lip} ((0, h) \times \R ) \, , \\ \noalign{\smallskip}   &\!\!-\partial _x w \geq f ^* _z (t, \partial _ t w) \ \mathcal L ^ 2 \hbox{-a.e. on } (0, h) \times \R\, , \\ \noalign{\smallskip}    &\!\! -\partial _x w \geq - f (t, 0) \ \forall t \in D \, ,\ \mathcal L ^ 1 \hbox{-a.e. on } (0, h) \Big \}\,.\end{array} 
\end{equation}

Notice that problem \eqref{MK} looks like the dual formulation of Monge-Kantorowich transport problem, with marginals equal to the Dirac masses at $(0, c)$ and $(h, c)$, and a modified gradient constraint with respect to the usual one $|\nabla w| \leq 1$.

Inspired by dynamic programming and optimal control, a natural candidate to solve \eqref{MK}  is  the {\it value function} 
\begin{equation}\label{vf}
V(x, t):= \inf \Big \{ \int_0 ^ x f (u, u') \, ds   \ :\ u \in H ^ 1 (0, h)  \, , \ u (0) = c \, , \ u (x) = t \Big \}\,,
\end{equation}
or equivalently a candidate calibration is the rotated gradient $(- \partial _ t V, \partial _ x V)$. 

Indeed, if $V$ is admissible in \eqref{MK}, it is automatically optimal. Namely, 
$\S \geq V (h, c)- V (0, c)  = \mathcal I$ 
and by Theorem \ref{t:nogap} the first inequality holds necessarily as an equality. 

Thus the key point is to check the admissibility of $V$ in problem \eqref{MK}.  
By using Bellman's optimality principle (see for instance \cite[Theorem 1.2.2]{CaSi}), 
it is easy to check that 
$V$ satisfies the constraints asked  in \eqref{MK} at every differentiability point. 
Unfortunately, it misses to satisfy the last important requirement of being Lipschitz regular close to $s=0$.
In Section \ref{secexamples} we shall be back to this phenomenon in connection with a relevant example of free boundary problem. 
}\end{remark}

\medskip

\section{Convexification recipe}\label{recipe}

The synopsis of this section is the following:

\medskip
\begin{itemize}
\item[--] in Subection \ref{construction} we introduce a convex functional $\widehat E$, defined in one more space dimension,
of the form $H +\ell$, with $H$ and $\ell$ conceived respectively with the aim of extending  the bulk and the surface parts of  the primal non-convex energy $E$; then we state the main result of the section (Theorem \ref{l:IJ}), which makes the link between the primal problem \eqref{defI} and a minimization problem for $\widehat E$. 
\smallskip

\item[--] in Subsection \ref{representation} we provide an integral representation result for  $H$; 
\smallskip

\item[--] in Subsection \ref{coarea} we state a generalized coarea formula, which turns out to be satisfied in particular by $H$
(as it can be seen thanks to its integral representation);

\smallskip
\item[--] in Subsection \ref{relationship} we prove Theorem \ref{l:IJ}
(by using in particular a slicing formula for $\widehat E$ which follows from the coarea formula for $H$); 

\smallskip
 
\item[--] in Subsection \ref{proof1} we prove the inequality $\mathcal I \geq \mathcal I ^*$, which is the easiest half of Theorem \ref{t:nogap}.
\smallskip
\end{itemize}

\medskip

\subsection{Construction of the convex extension of the primal energy}\label{construction}

As enlightened by the heuristics given in Section \ref{genesis}, 
the basic idea of our duality method is  to consider  the flux of suitable fields across the graph of functions $u$ admissible in the primal problem; and, along this way, we are naturally led to apply the divergence theorem on subgraphs.

Let us now fix these ideas in a systematic  setting, and develop them into the proposal of a convexification recipe: it consists in extending the non-convex energy introduced in \eqref{f:E} to a {\it convex functional defined in one more space dimension}.

Any element $u $ of $H ^ 1 (\Omega)$ can be identified  with a function in one more dimension, given by the characteristic function $\1 _u$ of its subgraph, defined on $\Omega \times \R$ by
$$ \1 _u (x, t) := \begin{cases}1 & \hbox{ if  }t \leq u (x) \\ 0 &  \hbox{ if  }t > u (x) \,.
 \end{cases}
 $$
Notice that $\1_u$  is not in
$L ^ 1 (\Omega \times \R)$, but merely in $L ^ 1 _{{\rm loc}} (\Omega \times \R)$.
 
Our target is to 
find a convex lower semicontinuous functional $\widehat E:  L ^ 1 _{{\rm loc}} (\Omega \times \R) \to \R \cup \{ + \infty \}$  
and a suitable subclass $\widehat{\mathcal A}$ of $L ^ 1 _{{\rm loc}} (\Omega \times \R)$ such that:

\medskip
\begin{itemize}
\item[--]  for every $u \in \mathcal A$, it holds $\1_u\in \widehat{\mathcal A}$ and $\widehat E(\1 _u)  = E (u)$ ;

\medskip
\item[--] the infimum $\mathcal I$ in  \eqref{defI} can be recast by minimizing $\widehat E$ over the class $\widehat {\mathcal A}$. 
\end{itemize}

\medskip
To that aim we are going to consider separately the {\it bulk part} and the {\it surface part} of the energy $E$. 
 
We start by recalling that, 
for any $u \in W ^ {1,p} (\Omega)$ (and actually more in general for  any $u\in BV (\Omega)$), its subgraph
is a set with finite perimeter \cite[p.~371]{GiMoSo}, or equivalently  $D \1 _u$ belongs to the space $\mathcal M (\Omega\times \R; \R ^ {N+1}) $ of vector valued bounded measures on $\Omega \times \R$.  
However, $\1 _u$ does not belong to $BV (\Omega \times \R)$,  since as already noticed it is not in
$L ^ 1 (\Omega \times \R)$, but merely bounded. We can thus say that $\1 _u$ belongs to the following subspace of $L ^ 1 _{\rm loc} (\Omega \times \R)$:
\begin{equation}\label{BVinfty}
BV _\infty (\Omega \times \R): = \Big \{ v \in L ^ \infty (\Omega \times \R ) \ :\ Dv \in \mathcal M (\Omega \times \R; \R ^ {N+1} ) \Big \}\,.
\end{equation}

For any $v\in BV _\infty (\Omega \times \R)$ and any  $\sigma$ in the space $X _ 1 (\Omega \times \R)$ defined in \eqref{X1}, a pairing 
$\sigma\cdot Dv$ can be defined as the following linear functional, which turns out to be a Radon measure   on $\Omega \times \R$
(see \cite[Thm 1.5 and Corollary 1.6]{A}) 
\begin{equation}\label{lf}
\langle (\sigma \cdot Dv) , \varphi \rangle:= - \int _{\Omega \times \R}v \, (\sigma \cdot \nabla \varphi  +   \varphi \div \sigma)   \, dx \qquad \forall \varphi \in \mathcal C ^ \infty _0 ( \Omega \times \R)\,.
\end{equation}
Moreover this measure is absolutely continuous with respect to $|Dv|$ and satisfies
\begin{equation}\label{lftotal}
\int_{\Omega\times\R}  |(\sigma \cdot Dv)| \le  \|\sigma\|_\infty \int_{\Omega\times\R}\, |Dv| \ .\end{equation}
Notice that definition \eqref{lf} reduces to \eqref{gen_trace} in the special case when $v$ is the characteristic function of a set $A \subset \Omega \times \R$ with finite perimeter.

We are now in a position to define on $L ^ 1 _{\rm loc}(\Omega \times \R)$ the following functional, which will give the required convex extension of the bulk part of the energy $E$: 
\begin{equation}\label{f:H}
H (v) := \begin{cases} \displaystyle
\sup \Big \{ \int _{\Omega \times \R} \sigma \cdot Dv \ :\ \sigma \in \K \Big \} & \text{ if } v \in BV _\infty (\Omega \times \R)
\\ 
+ \infty & \text{ otherwise.} 
\end{cases}
\end{equation}
An integral representation result for $H$ will be proved in Subsection \ref{representation} below. 
In particular, such result will disclose the crucial information that any function $v\in L^1  _{{\rm loc}} (\Omega \times \R)$ lying in the finiteness domain of $H$ satisfies a {\it monotonicity condition}, namely: 
\begin{equation}\label{mon}
H(v) < + \infty \ \ \Rightarrow \ \ \text{for } {\mathcal L} ^ N\text{-a.e.}\ x \in \Omega\, , \text{ the map } t \mapsto v (x, t)\text{ is  decreasing} .
\end{equation} 

We infer that, if $H (v ) < + \infty$, for ${\mathcal L} ^ N$-a.e.\ $x \in \Omega$ and ${\mathcal L} ^ 1$-a.e.\ $s \in (0,1)$, the set $ \{ \tau \in \R  \, :\,v (x, \tau ) \leq s \big \}$ is a nonempty half-line, and we can define for later use the function
\begin{equation}\label{defus} u _ s (x) := \inf \big \{ \tau \in \R  \ :\ v (x, \tau ) \leq s \big \}\,. \end{equation}
Notice that by construction the subgraph of $u _s$ agrees up to a Lebesgue negligible set with 
the level set $\{\tau \in \R\, : \, v (x, \tau) > s\}$, namely 
\begin{equation}\label{eqchar}
\1 _{u _s}(x, t) = \chi _{\{v > s\}}(x, t) \qquad  \hbox{ for ${\mathcal L} ^ {N+1}$-a.e. $(x, t) \in \Omega \times \R$}\,. 
\end{equation}

Next we turn our attention  to extend also the surface part of the energy $E$. To that aim we observe  that, 
though $\1_u \not \in L ^ 1 (\Omega \times \R)$,  it becomes  integrable after a suitable translation. 
Indeed, since $u$ is almost everywhere finite, 
 for a.e.\ $x  \in \Omega$ the map $t \mapsto \1_u (x, t)$ is monotone decreasing, with 
$$ \1_u (x, - \infty) = 1 \qquad {\rm and  }   \qquad \1_u (x, + \infty) = 0 \,.$$ 
We are thus led to introduce the reference function 
\begin{equation}\label{d:v0} v_0 (x, t):= \begin{cases} 1 & \hbox{ if  }t \leq 0 \\ 0 &  \hbox{ if  }t >0 \,.
\end{cases}
\end{equation}  
The equality
$$\int _\R |\1 _u (x, t) - v_0(x, t) | \, dt = |u (x)|\,,$$
 implies that $\1 _u - v _0 \in L ^ 1 (\Omega \times \R)$ as soon as $u \in L ^ 1 (\Omega)$.  
%
%

We infer that the class $\mathcal A$ introduced in \eqref{f:A} can be embedded, through the map $u \mapsto \1 _u$, into the class
\begin{equation}\label{defA} \A  := \Big \{ v \in BV _\infty (\Omega \times \R) \ :\ v -v _0 \in L ^ 1(\Omega \times \R)\, , \ v-\1_{u_0} = 0 \hbox{ on } \Gamma_0 \times \R \Big \}\,,
\end{equation}
where the last equality is intended in the sense of traces.


Notice in particular that, for every $v \in \widehat {\mathcal A}$, the function $v-v_0$ is in $BV (\Omega \times \R)$, so that it has a $L ^1$-trace on $\Gamma _1 \times \R$.  

We are then in a position to define on $L ^ 1 _{\rm loc}(\Omega \times \R)$ the following functional,  which will give the required convex extension of the surface part of the energy $E$: 
\begin{equation}\label{f:ell}
\ell (v) := \begin{cases} \displaystyle
\int_{\Gamma_1 \times \R } \gamma' (t) (v-v_0) \, d \mathcal H ^ {N -1}\, dt & \text{ if } v \in \widehat{\mathcal A}
\\ 
+ \infty & \text{ otherwise.} 
\end{cases}
\end{equation}

Finally, we set 
\begin{equation}\label{f:F}
\widehat E (v) := H ( v) + \ell (v)  \qquad \forall v \in L ^ 1 _{\rm loc} (\Omega \times \R)\, . 
\end{equation}

The next result states that
the functional $\widehat E$ and the class $\widehat{\mathcal A}$ thus defined fit exactly the target conditions demanded at the beginning of this section:

\begin{theorem}{\rm (link between the initial non-convex problem and its convex extension)}  \label{l:IJ} 
There holds 
\begin{equation}\label{FE} \widehat E (\1 _u ) = 
\begin{cases}
E (u) & \text{ if } u \in W ^ {1,p} (\Omega) \\ 
\noalign{\medskip}
+ \infty & \text{ if } u \in BV (\Omega) \setminus W ^ {1,p} (\Omega) \end{cases} 
\end{equation}
\begin{equation}\label{nog}
\inf \Big \{ E (u) \, :\, u \in {\mathcal A}  \Big \}= \inf \Big \{ \widehat E (v) \, :\, v \in \widehat{\mathcal A} \Big \}\,.
\end{equation} 
Moreover, both the infima in \eqref{nog} are finite and attained, and:

\medskip
-- if $u \in {\rm argmin } _{\mathcal A} (E)$, then $\1 _u \in  {\rm argmin } _{\widehat{\mathcal A}} (\widehat E)$; 

\medskip
-- if $v \in  {\rm argmin } _{\widehat{\mathcal A}} (\widehat E)$, then $u _s \in {\rm argmin } _{\mathcal A} (E)$  for $\mathcal L ^ 1$-a.e.\ $s \in (0,1)$ (with $u _s$ as in \eqref{defus}).

\medskip
In particular, if the primal problem $\inf \{ E (u) \, :\, u \in {\mathcal A} \}$ admits a finite number of solutions $\{ u ^1, \dots, u ^ k\}$, 
then 
\begin{equation}\label{cc}
{\rm argmin } _{\widehat{\mathcal A}} (\widehat E) = \sum _{ i=1} ^ k \theta _i \1 _{u ^i }\, , \qquad \theta _i \in [0,1]\, , 
\end{equation}
meaning that $v$ is a piecewise constant function.

\end{theorem}

The proof of  Theorem \ref{l:IJ}  will be given in Subsection \ref{relationship},  after developing the necessary tools 
in Subsections \ref{representation} and \ref{coarea}.

\subsection{Integral representation of $H$}\label{representation}

Let us introduce the one-homogeneous convex integrand $h_f$ which will appear in the integral representation of $H$. 
Such integrand has been already used in several previous works exploiting the classical identification between BV functions and subgraphs of finite perimeter (see for instance \cite{Da}). Its definition reads as follows: 

\begin{definition}\label{d:hf} For $(t, q) \in \R \times \R ^ {N+1}$, we set:
\begin{equation}\label{defhf}
h_f(t,q):=
\begin{cases} -q^t f\left(t,-q^x/q^t\right)& \hbox{ if } q^t<0 \\ 
 +\infty\ \ \ \  & \hbox{ if } q^t >  0 \hbox{ or } q^t=0, q^x\neq 0\\ 
 0\ \ \ \  & \hbox{ if } (q^x , p^t) = (0, 0).
 \end{cases}
 \end{equation}
\end{definition}

The above definition will look more natural recalling that it takes its origins  in  
Convex Analysis, as it corresponds precisely to the support function of the epigraph of the Fenchel conjugate $f ^* _z (t, \cdot) $, namely of the set $K (t)$ introduced in \eqref{f:C} (see \cite[Section 13]{Ro}).   For convenience of the reader, this and the other main properties of $h _f$ are stated below.

 \begin{lemma}[properties of $h _f$] \label{l:proph}

The function $h_f$ is lower semicontinuous in $(t, p)$ and convex, positively $1$-homogeneous in $p$. 

Moreover, $h _f ( t, \cdot)$ is the support function of  the convex set $K(t)$ introduced in \eqref{f:C}, or equivalently  the Fenchel conjugate of the indicatrix function $I _{K(t)}( \cdot )$ (which equals $0$ on $K (t)$ and $+ \infty$ outside): 
\begin{equation}\label{eqh}
h _f(t, q ) =  \sup \big \{ q \cdot \tilde q \ :\ \tilde q \in K (t) \} = I_ {K(t)} ^ * (q) \,.
\end{equation}
In particular, the map $t \mapsto K (t)$ defines a lower semicontinuous multifunction
(meaning that 
$\big \{ t \in \R \, :\, K (t) \cap A \neq \emptyset \}$ is open for every open subset $A$ of $\R ^ {N+1}$). 

\end{lemma}

\proof Since by assumption
$f$ is lower semicontinuous in $(t, z)$,  it is clear that $h_f$ is  l.s.c.\ at any $(t,q)$ with $q^t < 0$.  Let us assume  that $q ^ t \geq 0$, and let $(t_n,q_n)$ be a sequence converging to $(t,q)$, with $\liminf_{n\to\infty} h_f(t_n,q_n) = l \in [0, +\infty)$
(otherwise there is nothing to prove).
Then, possibly passing to a subsequence,  for every $n$ it holds $q_n^t \leq 0$, hence $q ^ t = 0$. Recalling the growth condition from below satisfied by $f$, we infer that 
$$|q^t_n| \Big ( \alpha \Big| \frac{q^x_n }{q^t_n } \Big| ^p  - r ( t_n) \Big ) \leq l\, , $$ 
and therefore also $q^x=0$, so that $h_f(t,q)=0\leq l$.

It is immediate from the definition of $h _f$ that $h _f (t, \cdot)$ is positively $1$-homogeneous. The proof of equality (\ref{eqh}), which in particular implies the convexity of $h _ f (t , \cdot)$, can be found in \cite[Corollary 13.5.1]{Ro}, but for the sake of completeness we sketch it below. 
By definition, it holds
$$I_ {K(t)} ^ * (p) = \sup \big \{ (q^ x \cdot \tilde q ^ x + q ^ t \tilde q ^ t ) \ :\  \tilde q ^ t \geq f ^*_z (t, \tilde q ^x) \big \}\,.$$ 
It is immediately seen the above supremum is $0$ in case $ q^t =|q^x|=0$, and $+ \infty$ in case $q^t >  0$ or  $q^t=0, q^x\neq 0$. 
In case $p ^ t <0$, it holds
$$\begin{array}{ll}
\sup \big \{ (q^ x \cdot \tilde q ^ x + q ^ t \tilde q ^ t ) \ :\ \tilde q ^ t \geq f ^*_z (t, \tilde q ^x) \big \} & = 
\displaystyle \sup \Big \{ (q^ x \cdot \tilde q ^ x + q ^ t  f ^*_z (t, \tilde q ^x)  ) \ :\ \tilde q ^ x \in \R ^N \Big \}  \\ \noalign{\smallskip} & = 
\displaystyle - q ^ t \sup \Big \{ (- \frac{q^ x}{q ^ t} \cdot \tilde q ^ x -   f ^*_z (t, \tilde q ^x)  ) \ :\ \tilde q ^ x \in \R ^N \Big \}   \\ \noalign{\smallskip}& = - q ^ t f ( t , - \frac{q^ x}{q ^ t}) \,. \end{array}$$ 
Finally, the lower semicontinuity of the multifunction $t \mapsto K(t)$ follows from \cite[Theorem 17]{BoVa}.
\qed

\bigskip

As a last ingredient, let us  recall that one-homogeneous convex integrands such as $h _f$ can be integrated in the sense of measures. More precisely,   for any bounded vector-valued measure $\lambda \in \mathcal M (\Omega \times \R; \R ^{N+1})$,  the integral of $h _f (t, \lambda)$ is meant as 
$$
\int _{\Omega \times \R}  h _ f (t, \lambda) := \int _{\Omega \times \R}   h _ f \Big (t , \frac{d \lambda}{d| \lambda|} \Big) \, d |\lambda|\,,
$$
where $|\lambda|$ is the total variation measure of $\lambda$. 

Such convex one-homogeneous functional on measures has been studied in \cite{BoVa}. In particular, it can be characterized 
in terms of the duality $\langle \ , \ \rangle$  between $\mathcal M (\Omega \times \R; \R ^{N+1})$ and $\mathcal C _ 0 (\Omega \times \R; \R ^ {N+1 } )$ according to the next lemma:

\begin{lemma}\label{l:dualite} There holds
\begin{equation}\label{dualite}
\int _{\Omega \times \R}  h _ f (t, \lambda)  = \sup \Big \{ \langle  \lambda, \psi \rangle  \ : \ \psi \in \mathcal C _ 0 (\Omega \times \R; \R ^ {N+1 } ) \, , \ \psi (x, t) \in K (t) \hbox{ on } \Omega \times \R \Big \}
\,.\end{equation}
Moreover, the equality above is still true if the supremum at the right hand side is restricted to functions $\psi \in \mathcal D (\Omega \times \R; \R ^ {N+1})$. 
\end{lemma}

\proof  It is easy to check that the supremum at the right hand side of \eqref{dualite} is not larger than $\int_{\Omega \times \R} h _ f (t, \lambda)$.
This follows by applying 
the inequality $h _f (t, q) \geq q \cdot \tilde q$,  holding  for every $\tilde q \in K (t)$,  with $q = \frac{ d \lambda}{d |\lambda|}$ and $\tilde q = \psi$. 

Therefore, the proof of the lemma is concluded if we show that
\begin{equation}\label{dualite2}
\int _{\Omega \times \R}  h _ f (t, \lambda)  = \sup \Big \{ \langle  \lambda, \psi \rangle  \ : \ \psi \in \Sigma \Big \} \, , 
\end{equation}
with
$$\Sigma := \Big \{ \psi \in \mathcal D(\Omega \times \R; \R ^ {N+1 } ) \, , \ \psi (x, t) \in K (t) \hbox{ on } \Omega \times \R  \Big \}\,.$$
Clearly, in \eqref{dualite2} we can replace $\Sigma$ by its closure $\overline \Sigma$ (in the uniform norm of $\mathcal C _ 0 (\Omega \times \R; \R ^ {N+1 } )$). Then,  according to \cite[Theorem 5]{BoVa}, in order to prove \eqref{dualite2} it is enough to establish that
$$\overline \Sigma = \Big \{ \psi \in \mathcal C_0 (\Omega \times \R; \R ^ {N+1 } ) \, , \ \psi (x, t) \in K (t) \hbox{ on } \Omega \times \R  \Big \}\,.$$
As $\Sigma$ is $\mathcal C^ \infty$-convex, we may apply \cite[Proposition 10]{BoVa}, yielding
$$\overline \Sigma = \Big \{ \psi \in \mathcal C_0 (\Omega \times \R; \R ^ {N+1 } ) \, , \ \psi (x, t) \in \Gamma (x, t) \hbox{ on } \Omega \times \R  \Big \}\,,$$
 with 
 $$\Gamma (x, t)  := \overline{ \big \{ \psi ( x, t) \ :\ \psi \in \Sigma \big \}  } \,.$$
 Thus we are reduced to prove the equality $K ( t_0) = \Gamma (x_0, t_0)$ for every $(x_0, t_0) \in \Omega \times \R$. 
Since $K (t_0)$ is closed, it is immediate that $\Gamma (x_0, t_0) \subseteq K (t_0)$. Conversely, let $z \in {\rm int} ( K (t_0))$. 
There exists $\delta >0$ such that, for $|t- t _0| < \delta$, we have $z \in K (t)$ (see \cite[Lemma 15]{BoVa}), 
and consequently the whole interval $[0, z]$ lies in $K (t)$ for $|t- t _0| < \delta$. 
Then we define $\psi (x, t) := z \alpha (x) \beta _\delta (t)$, being $\alpha \in \mathcal D (\Omega; [0, 1])$, $\beta _\delta \in \mathcal D (\R; [0, 1])$ with ${\rm spt} (\beta _\delta) \subset [ t _0 - \delta, t _0 + \delta] $, and $\alpha ( x_0) = \beta _0 ( t _0) = 1$. It is easy to check that the function $\psi$ belongs to $\Sigma$, and hence $z \in \Gamma (x_0, t_0)$.  Since $\Gamma (x_0, t_0)$ is closed, and $K ( t_0)$ coincides with the closure of its interior, we have proved that
$\Gamma (x_0, t _0) \supseteq K ( t_0)$. \qed \bigskip


\medskip

We are now ready for the announced integral representation result. 

\begin{proposition}[integral representation of $H$]\label{repH} For every $v \in BV _\infty (\Omega \times \R)$,
the functional $H$ defined in \eqref{f:H} satisfies the equality
 $$H (v) 
 = \int_{\Omega \times \R} h _ f (t, Dv)\,.$$
\end{proposition}

\proof In view of the  definition \eqref{f:H} of the functional $H$ and of Lemma \ref{l:dualite}, for every $v \in BV _\infty (\Omega \times \R)$ there holds
$$H (v) \geq  \sup \Big \{\int_{\Omega \times \R} \sigma \cdot Dv \ :\ \sigma \in \K \cap \mathcal D ( \Omega \times \R; \R ^ {N+1})   \Big \}=  \int_{\Omega \times \R} h _f (t, Dv)\,.$$
To obtain also the converse inequality we have to show that, for every $\sigma \in \K$ and every $v \in BV _\infty (\Omega \times \R)$, there holds $\int_{\Omega \times \R } \sigma \cdot Dv  \leq \int_{\Omega \times \R} h _ f (t, Dv)$. This is established in the lemma below which completes our proof. \qed

 \begin{lemma}[lower bound for $H$] \label{l:julygen}
For every  $\sigma \in \K$ and every $v \in BV _\infty (\Omega \times \R)$, there holds $$\int_{\Omega \times \R } \sigma \cdot Dv  \leq \int_{\Omega \times \R} h _ f (t, Dv)\,.$$
\end{lemma}

\proof The lemma will be obtained by showing separately the following two inequalities:
\begin{eqnarray}
& \displaystyle \int _{ \Omega \times (\R \setminus D ) } \sigma \cdot Dv  \leq  \int _{  \Omega \times (\R \setminus D )  } h _f (t, Dv) & \label{f:intermedia2}
\\ \noalign{\smallskip}
& \displaystyle \int _{ \Omega \times D  } \sigma  \cdot Dv    \leq  \int _{\Omega \times D  }
h _f (t, Dv) \,. & \label{f:intermedia1}
\end{eqnarray}

In order to prove \eqref{f:intermedia2}, we need to exploit some facts established in \cite{Abis, A} (see also 
\cite{BoDa}). 
Recall that, for every  $\sigma \in \K$ and  $v \in BV _\infty (\Omega \times \R)$, the measure $\sigma \cdot Dv$ defined in \eqref{lf} is absolutely continuous with respect to $|Dv|$ ({\it cf.} \eqref{lftotal}).  Moreover, setting $\nu _ v:= \frac{(\partial Dv) }{\partial |Dv|}$, the Radon-Nikodym derivative
of $\sigma \cdot Dv$ with respect to $|Dv|$ is given by
\begin{equation}\label{anz}
\frac{d (\sigma \cdot Dv )}{|Dv|} = q _\sigma (x, \nu _v) \qquad |Dv|\text{-a.e. in }  \Omega \times \R \,,
\end{equation}
where $q _\sigma: \Omega \times \R \times S ^N \to \R$ is the Borel function given by
$$q _{\sigma} ((x, t), \zeta) := \limsup _{\rho \to 0 ^+} \limsup _ {r \to 0 ^+}  \frac{1}{ \mathcal L ^ {N+1} \big ( C_{r, \rho} ((x, t), \zeta) \big ) }
{\int_{C_{r, \rho} ((x, t), \zeta) } \sigma (y, s) \cdot \zeta \, d \mathcal L ^ {N+1}}
\,,$$
being 
$$C_{r, \rho} ((x, t), \zeta):= \Big \{ (y, s) \in \R ^ {N+1} \ :\  |(y-x, s-t) \cdot \zeta| \leq r \, , \, \big | (y-x, s-t) - \big ( (y-x, s-t) \cdot \zeta \big ) \zeta \big |  \leq \rho   \Big \} \,.$$ 

%
%
%
In view of \eqref{anz}, we can rewrite \eqref{f:intermedia2} as
\begin{equation}\label{f:intermedia3}\int _{\Omega \times (\R \setminus D) } q _\sigma \big ( (x, t), \nu _v \big )  \, d |Dv|  \leq  \int _{\Omega \times (\R \setminus D) } h _f (t, \nu _v) \, d |Dv| \,.
\end{equation}
We observe that 
\begin{equation}\label{f:claim}
q _\sigma (( x_0, t_0), \zeta) \leq 
h _ f ^ + ( t_0, \zeta):= \limsup _{t \to t_0} h _ f (t, \zeta) \qquad \forall (x_0, t_0 ) \in \Omega \times \R \,,\ \forall \zeta \in S ^N\,.
\end{equation}
Namely,  
since $\sigma$ satisfies condition \eqref{s3} (or equivalently \eqref{s3'}), by Lemma \ref{l:proph} it holds
$$\sigma (x, t) \cdot \zeta \leq h _ f (t, \zeta) \quad  \ \text{ for } \mathcal L ^ {N+1}\text{-a.e.}\ (x, t) \in \Omega \times \R\,.$$
By taking the mean value over the cylinder $C_{r, \rho} ((x_0, t_0), \zeta)$, and passing to the limsup as $\rho$ and $r$ converge to zero, we obtain
$$
\begin{array}{ll} q _\sigma (( x_0, t_0), \zeta) & \leq  \displaystyle \limsup _{\rho \to 0 ^+ } \limsup _{r \to 0 ^+} 
 \frac{1}{ \mathcal L ^ {N+1} \big ( C_{r, \rho} ((x_0, t_0), \zeta) \big ) }
{\int_{C_{r, \rho} ((x_0, t_0), \zeta) }\!\!h _ f (t, \zeta) \, d \mathcal L ^ {N+1}}
\\ \noalign{\medskip}
& \leq h _ f ^+ ( t _0, \zeta) \,.
\end{array}$$
Now we notice that,  thanks to our hypothesis \eqref{H1bis}, we have
\begin{equation}\label{f:intermedia4}
h _f ^+(t, \zeta) = h _f (t, \zeta)  \qquad \forall t \in \R \setminus D\, , \ \forall \zeta \in  S ^N \,.\end{equation}
The required inequality \eqref{f:intermedia3} follows from
\eqref{f:claim} and \eqref{f:intermedia4}.

%

Let us now prove inequality \eqref{f:intermedia1}.  
To that aim it is enough to show that 
$\frac{d (D_x v) } {d|Dv| }= 0$  $|Dv|\text{-a.e.\ on }  \Omega \times D$, or equivalently that
\begin{equation}\label{nuvert}
\nu _v = - e _{N+1} \qquad \mathcal |Dv|\text{-a.e.\ on }  \Omega \times D \,.
\end{equation}
Indeed in this case, by exploiting condition \eqref{s4} (or equivalently \eqref{s4'}), we obtain 
$$ \sigma (x, t) \cdot (- e _{N+1} ) = - \sigma ^ t (x, t) \leq  f (t, 0) = h _ f (t, - e _{N+1}) 
\  \forall\, t \in  D  \text{ and for } \mathcal L ^ {N}\text{-a.e.}\ x \in \Omega. $$

To prove \eqref{nuvert} we simply observe that 
$$0 = \int _ D \Big ( \int _\Omega |D_x v (\cdot, t) | \Big ) \, dt = \int _{\Omega \times D} |D_x v| \,,$$
where the first equality follows from the assumption $\mathcal L ^ 1 ( D) = 0 $, and the second one from the slicing formula for $BV$ functions (see \cite[Section 3.11]{AmFuPa}). 
\qed

\bigskip
\begin{remark}\label{r:easyine2} As one can easily check by inspection of the proof of Lemma \ref{l:julygen},  the inequality \eqref{bflux} can be strengthened into
$$\int_{\omega}  f(u, \nabla u) \,dx = \int_{G_u \cap ( \omega \times \R) }  h _f (t,  \nu_u) \, d\mathcal H^N   \geq \int_{G_u \cap ( \omega\times \R) }  \sigma \cdot \nu_u \, d\mathcal H^N   \quad  \forall \, \omega \text{ Borel set } \subset \Omega \,.$$
By the arbitrariness of the Borel set $\omega$  we infer that,  for all $u \in \mathcal A$ and $\sigma \in \X$, there holds 
\begin{equation}\label{f:fun}
h _ f (t, \nu _u) \geq \sigma \cdot \nu _u \qquad \mathcal H ^N\hbox{-a.e.\ on } G_u\,.
\end{equation}
Consequently, we see that the validity of  inequality \eqref{s4}  is extended for free also to values $t \in \R \setminus D$. Indeed, 
by taking locally constant functions $u$ in \eqref{f:fun} we obtain that,  for all $\sigma \in \X$,  
there holds $f (t, 0 ) \geq - \sigma ^ t(x, t) $ for every  $t \in \R$ and $\mathcal L ^N$-a.e.\  $x\in \Omega$. 
\end{remark}

\bigskip

\subsection{Generalized coarea formula}\label{coarea}

Let $A$ be an open subset of $\R^d$.  For every function $v \in L^1_{\rm loc}(A) $ and every $s \in \R$, let 
$\chi _{\{v>s\}}$ denote the characteristic function of the set $\{v>s \}$, {\it i.e.}
$$\chi_{\{v>s\}} (x) := \begin{cases}1 & \hbox{ if  } v (x) >s \\ 0 &  \hbox{ if  } v (x) \leq s \,.
 \end{cases}
 $$
Following a terminology introduced in \cite{Vi}, we give the following

\begin{definition}

We say that a functional $J: L^1_{\rm loc}(A) \to [0,+\infty]$ 
satisfies the {\it generalized coarea formula} if
for every $u\in L^1_{\rm loc}(A)$ the function $t \mapsto J(\chi _{\{v>s\}})$
is Lebesgue-measurable on $\R$ and there holds
\begin{equation}\label{gcf}
J(v) = \int_{-\infty}^{+\infty} J(\chi_{\{v>s\}}) \, ds \qquad \forall v \in  L ^ 1_{\rm loc} (A) \,.
\end{equation}
\end{definition}

\begin{remark}{\rm It is readily seen that the following conditions are necessary in order that a functional $J : L^1_{\rm loc}(A) \to  [0,+\infty] $ 
satisfies the generalized coarea formula:

\smallskip 
 -- $J$ is positively $1$-homogeneous ({\it i.e.}\ $J(\lambda v)
= \lambda J(v)$ for all $v \in L ^ 1_{\rm loc} (A)$ and $\lambda \geq 0$)

\smallskip
-- $J (\chi _A) = 0$. 

\smallskip
 
Indeed, the $1$-homogeneity is immediately obtained via a change of variable in (\ref{gcf}), whereas the second property follows by applying (\ref{gcf}) to $v = \chi _A$, which gives 
$J(\chi _A) =
\int^1_{-\infty} J(\chi _A) dt$.}
\end{remark}

\medskip The next result establishes sufficient conditions in order that a functional $J$ 
satisfies the generalized coarea formula. Its proof is postponed to Section \ref{proofcoarea}.

\begin{theorem}{\rm (generalized coarea formula)} \label{theo}
Let $J : L^1_{\rm loc}(A) \to [0,+\infty]$ be positively $1$-homogeneous
 and such that $J(\chi _A) = 0$.
Assume in addition that $J$ is convex, lower semicontinuous,  and satisfies the following property: if  $\{\alpha _i\} _{1 \leq i \leq k}$ is a family of functions in ${\mathcal{C}}^\infty(\R;[0,1])$ with $\sum^k_{i=1} \alpha_i \equiv 1$, setting $\beta_i(t) := \int_0^t \alpha_i(s) ds$, it holds
\begin{equation}\label{sub} \sum^k_{i=1} J(\beta _i \circ v) \leq J(v) \qquad \forall v \in L ^ 1_{\rm loc} (A)\,.\end{equation}
Then $J$ satisfies the generalized coarea formula.
\end{theorem}

\begin{remark}
{\rm It is easy to check that the functional $J: L^1_{\rm loc}(A) \to [0,+\infty]$ defined by  $\int_A |Dv|$ if $u \in BV(A)$ and $+ \infty$ otherwise 
fulfills all the hypotheses of Theorem \ref{theo}.
Hence $J$ satisfies the generalized coarea formula, which allows to recover the classical coarea formula $
\int_A |Dv| =
\int_{-\infty}^{+\infty} {\rm Per} ( \{ v > t \} ) \,  dt 
$ holding for every function  $u \in BV(A)$ (see {\it e.g.}\ \cite[p.145]{AmFuPa}).
}
\end{remark}

Theorem \ref{theo} applies in particular to the functional $H$,  as stated below. 

\begin{proposition}{\rm (coarea formula for $H$)} 
\label{l:coarea} 
The functional $H$  satisfies  the generalized coarea formula. 
\end{proposition}

\proof  Let us  check  that $H$ satisfies all the assumptions of Theorem \ref{theo}. 
It is clear from  definition  \eqref{f:H} 
that $H$ is positively $1$-homogeneous, convex, lower semicontinuous, and satisfies $H (\chi _{\Omega \times \R} ) = 0$.
It remains to check that, if  $\{\alpha _i\} _{1 \leq i \leq k}$ is a family of functions in $\mathcal C ^\infty(\R; [0,1])$ with $\sum^k_{i=1} \alpha_i \equiv 1$, and
$\beta_i(t) := \int_0^t \alpha_i(s) ds$, 
the inequality (\ref{sub}) holds. 
To that aim we may assume without loss of generality that $H (v) < + \infty$, namely that $v \in BV _{\infty} (\Omega \times \R)$. We observe that
$v \in BV _{\infty} (\Omega \times \R)$ implies $\chi _{\{v> s \} } \in BV _\infty (\Omega \times \R)$ for $\mathcal L ^1$-a.e.\ $s \in \R$. 
%
%
%
%
Then, according to 
Proposition \ref{repH}, we have to prove that
\begin{equation}\label{tesic}\int _{\Omega \times \R} \!\!\!
h _ f (t, Dv) \geq  \sum _{i=1}^ k \int _{\Omega \times \R} \!\!\!
h _ f (t, D(\beta _i \circ v))    \,.
\end{equation}
Denoting by $\widetilde D v$ the diffuse part of the measure $Dv$ (namely the sum of the absolutely continuous part plus the Cantor part), by the chain rule formula  \cite[Theorem 3.96]{AmFuPa}, for every $i = 1, \dots, k$, the function $\beta _i (v)$ belongs to $BV _\infty (\Omega \times \R)$ and it holds
$$
D(\beta_i(v))=
\alpha_i(v)\widetilde D v \ +\ 
\left( \int_{v^-(x,t)}^{v^+(x,t)}\alpha_i(s)\,ds \right)
\nu_v(x,t)d {\mathcal H}^{N} \res J_v \,.
$$
Then, since $\alpha_i$ are nonnegative functions and $h_f(t,\cdot)$ is positively
$1$--homogeneous, we have
$$
\int _{\Omega \times \R} \!\!\!
h _ f (t, D(\beta _i \circ v))  =  \int_{\Omega\times \R} \alpha_i(v)h_f(t,\widetilde D v)
\ +\ 
\int_{J_v} \left( \int_{v^-(x,t)}^{v^+(x,t)}\alpha_i(s)\,ds \right)
h_f(t,\nu_v(x,t)) d{\mathcal H}^{N}\,.
$$
Summing over $i$, and recalling that $\sum_{i=1}^k \alpha_i \equiv 1$, we deduce that \eqref{tesic} is satisfied with equality sign. Thus (\ref{sub}) holds, we are in a position to apply Theorem \ref{theo}, and we obtain that $H$ satisfies the generalized coarea formula. 

\qed

\subsection{Proof of Theorem \ref{l:IJ}}\label{relationship} 
We are going to prove the theorem in  two steps.  
In the first step we prove the equality \eqref{FE}, and 
in the second one we prove the equality \eqref{nog} and the subsequent part of the statement.  
For the second step we need a slicing formula for the functional $\widehat E$ (stated in Proposition \ref{p:slicing} below), 
which is obtained thanks to the the coarea formula for $H$ proved in the previous subsection.

\medskip

{\it Step 1 (Proof of $\eqref{FE}$)}.  Let us show separately the two equalities
\begin{eqnarray}
& \displaystyle{\int _{\Omega \times \R} \!\!\!
h _ f (t, D \1_u)} = \begin{cases}
\displaystyle \int _\Omega f(u, \nabla u ) \, dx & \text{ if } u \in W ^ {1,p} (\Omega) \\
+ \infty &  \text{ if } u \in BV (\Omega) \setminus W ^ {1,p} (\Omega) \end {cases}
& \label{FE1} 
\\ \noalign{\medskip}
 &\displaystyle  \int_{\Gamma_1 \times \R} \!\!\! (\1_u -v_0 ) \gamma' (t) \, d \mathcal H ^ {N-1} (x) = \int _{\Gamma _1} \gamma (u) d \mathcal H ^ {N-1} (x)
 \qquad \forall u \in BV (\Omega)
  \,. 
& \label{FE2}
\end{eqnarray}

In order to show \eqref{FE1} it will 
be useful to recall few basic facts about subgraphs of $BV$ functions. 
For any $u \in BV (\Omega)$, the singular set of $\1 _u$, or equivalently the measure theoretic boundary of the subgraph of $u$, is called
the {\it complete graph} of $u$, and is denoted by $\Gamma _u$.  Moreover, we set $\nu _{\Gamma _u}$ the inward unit normal to 
$\Gamma _u$.  In particular, we have
$$D \1 _u = \nu _{\Gamma _u} \,d (\mathcal H ^ N \res \Gamma _u )\,,$$
and 
$$ \int _{\Omega  \times \R}
h _f (t , D\1_ u) = \int  h _ f (t, \nu _{\Gamma _u}) \, d ({\mathcal H } ^N \res {\Gamma _u})\,.$$

By writing  $D \1 _{u}$ as the sum of the two measures
$$D \1 _{u} \res ( J _u \times \R) \qquad \hbox{ and } \qquad D \1 _{u} \res ( (\Omega \setminus J _u) \times \R)\,,$$
where $J _u$ denotes the jump set of $u$, one obtains a decomposition of $\Gamma _u$  into a  ``vertical part'' plus an ``approximately continuous part''.
On the vertical part, $\nu _{\Gamma _u} $ is horizontal, and precisely it is given by
\begin{equation}\label{nuhor}
\nu _{\Gamma _u}(x, t) = \big ( \nu _{J _u}(x), 0 \big ) \,.
\end{equation}
On the approximately continuous part, denoting by  $\displaystyle{ u _ + (x) = {\rm aplimsup} _ {y \to x} u (y)}$, 
$\nu _{\Gamma _u}$  is given by
\begin{equation}\label{nucont}
\nu _{\Gamma _u}(x, u _+(x)) = \frac{(\nabla u (x), - 1) }{\sqrt { 1 + |\nabla u(x)| ^ 2 }} 
\end{equation}
if $u$ is approximately differentiable at $x$ (with approximate gradient $\nabla u (x)$), 
and it is horizontal otherwise (namely at points corresponding to the Cantor part of $Du$). 
We refer to \cite[Section 4.1.5]{GiMoSo} for a detailed account of these properties. 

In particular, when dealing with functions $u \in W ^ {1,p} (\Omega)$, the complete graph $\Gamma _u$ agrees with the usual graph 
$G_u$, and $ \nu _u (x, u (x))  = \nu _{\Gamma _u}(x, u _+ (x))$. 

\medskip

Then, from the explicit expression  (\ref{defhf}) of $h _f$  and the fact that $\nu _{\Gamma_u}$ is horizontal except at point $(x, u (x))$ where $u$ is approximately differentiable, we see that
$ \int _{\Omega  \times \R}
h _f (t , D\1_ u)$  is finite only if $u \in W ^ {1, 1} (\Omega)$.
In this case, 
the measure $D \1 _u$ is given by $\nu _u \mathcal H ^N \res G _u$, and 
we have
$$ \int _{\Omega  \times \R}
h _f (t , D\1_ u) = \int  _ {\Gamma _u} h _ f (t, \nu _{ u}) \, d{\mathcal H } ^N \,.$$
Since the Jacobian of the mapping $\Omega \ni x \mapsto (x, u (x)) \in \Gamma _u$ is given by $\sqrt {1 + |\nabla u| ^ 2}$ and since $h _f (t, \cdot)$ is  positively $1$-homogeneous, via change of variable we get $$  \int_{G_u} h _ f ( t, \nu _u) \, d \mathcal H ^N = \int _\Omega h _f (u (x) , ( \nabla u (x), -1)) \, dx\,.$$ 
Now, by using the definition (\ref{defhf}) of $h _f$, it is immediate to check that the r.h.s.\ of the above equality agrees with 
$\int _\Omega f (u, \nabla u)\,dx$, which yields \eqref{FE1}. 

The identity
$$(\1_u - v_0 )   
= \begin{cases}
\1 _{[0, u (x)]} & \text{ if } u (x) >0
\\ 
- \1 _{[u (x), 0]} & \text{ if } u (x) <0\,.
\end{cases}
 $$
together with $\gamma (0 )= 0$, yields
\begin{equation}\label{iduseful}
\int _\R (\1_u - v_0 ) \gamma' (t) \, dt   = \gamma (u)\,. 
\end{equation}
We obtain \eqref{FE2} after an integration over $\Gamma _1$. 
The identity \eqref{FE} follows by adding \eqref{FE1} and \eqref{FE2}. \qed

\bigskip
\begin{proposition}\label{p:slicing}
For every $v \in \widehat {\mathcal A}$ such that $\widehat E (v) < + \infty$, there holds
\begin{equation}\label{super}
\widehat E (v) = \int _0 ^ 1 \widehat E(\chi _{\{v > s\}})  \, ds 
= \int _0 ^ 1 \widehat E(\1 _{u _s})  \, ds  = \int_0 ^ 1 E (u _ s ) \, ds\,.
\end{equation}
\end{proposition}

\proof  Let $v \in \widehat {\mathcal A}$ be such that $\widehat E (v) < + \infty$. We claim that it holds
\begin{equation}\label{f:slicel}
\ell (v) = \int _0 ^ 1 \ell (\chi _{\{v > s\}})\, ds = \int_0 ^ 1 \ell ( \1 _{u_s} ) \, ds = \int_0 ^ 1 \int _{\Gamma_1} \gamma (u _s) \, d \mathcal H ^ {N-1} \, ds \,. 
\end{equation} 
Notice that the second and the third equalities in \eqref{f:slicel} are satisfied in view of \eqref{eqchar} and \eqref{FE}. 
Thus we have just to prove the first equality, which can be rewritten as
$$\int _{\Gamma _ 1 \times \R}  \gamma ' (t) (v- v_0)  \, d \mathcal H ^ {N-1} \, dt  = \int_0 ^ 1 \Big \{  
\int _{\Gamma _ 1 \times \R}  \gamma ' (t) (\chi _{\{v>s\}} - v_0)  \, d \mathcal H ^ {N-1} \, dt \Big \} \, ds\,.$$
We write
$$
\begin{array}{ll}
& \displaystyle \ell (v) = \int_{\Gamma_1 \times \R_+ } v (x, t) \gamma ' (t) \, d \mathcal H ^ {N-1} \, dt -
 \int_{\Gamma_1 \times \R_- } [1-v (x, t)] \gamma ' (t) \, d \mathcal H ^ {N-1}  \, dt
\\ 
\noalign{\medskip}
& \displaystyle \ell (\chi _{\{v>s\}}) = \int_{\Gamma_1 \times \R_+ } \chi _{\{v>s\} }(x, t) \gamma ' (t) \, d \mathcal H ^ {N-1} \, dt -
 \int_{\Gamma_1 \times \R_- } [1- \chi _{\{v>s\} } (x, t)] \gamma ' (t) \, d \mathcal H ^ {N-1} \, dt\,.
\end{array}
$$
Now we observe that, since $v - v_0 \in L ^ 1 (\Omega \times \R)$, and $v(x, \cdot)$ is nonincreasing, 
$v$ takes values into $[0, 1]$. 
Any function $w$ with values in $[0,1]$ can be written as 
$w (x) = \int_0 ^ 1 \chi _{ \{w > s\} } \, ds$ (which is commonly called {\it layer cake representation formula}).  
Then, by applying  Fubini Theorem separately to the integrals over $\Gamma_1 \times \R _+$ and 
over $\Gamma _ 1 \times \R _-$, we have:
$$\int_{\Gamma_1 \times \R_+ } v (x, t) \gamma ' (t) \, d \mathcal H ^ {N-1}  \, dt  = \displaystyle  \int _0 ^ 1 \, ds 
\int_{\Gamma_1 \times \R _+} \chi _{\{v>s\}}  (x, t) \gamma ' (t) \, d \mathcal H ^ { N-1}  \, dt $$

\medskip
$$
\begin{array}{ll}
\displaystyle \int_{\Gamma_1 \times \R_- } [1-v (x, t)] \gamma ' (t) \, d \mathcal H ^ {N-1}  \, dt & \displaystyle = \int _0 ^ 1 \, d\tau 
\int_{\Gamma_1 \times \R _-} \chi _{\{1-v>\tau \}}  (x, t) \gamma ' (t) \, d \mathcal H ^ { N-1} \, dt 
\\ \noalign{\medskip}
&\displaystyle = \int_0 ^ 1 \, ds \int_{\Gamma_1 \times \R _-}   [1- \chi _{\{v>s\} } (x, t)] \gamma ' (t) \, d \mathcal H ^ {N-1} \, dt\,,
\end{array}
$$
and we obtain \eqref{f:slicel} by addition. 

We are now ready to conclude. By using Proposition \ref{l:coarea},  the equality \eqref{f:slicel}, 
the fact that (as noticed above) $v$ takes values into $[0, 1]$, and the equalities
$H (0) = H ( \chi _{\Omega \times \R} ) = 0$, we obtain
$$\begin{array}{ll} \widehat E (v) = H ( v ) + \ell (v) & \displaystyle = \int _{- \infty} ^ {+ \infty}
H (\chi _{\{v>s\}}) \, ds + \int_0 ^ 1 \ell (\chi _{\{v>s\}}) \, ds 
\\  \noalign{\medskip}
& \displaystyle = \int _{ 0} ^ {1}
H (\chi _{\{v>s\}}) \, ds + \int_0 ^ 1 \ell (\chi _{\{v>s\}}) \, ds  
\\  \noalign{\medskip}
& \displaystyle = \int _{ 0} ^ {1}
\widehat E (\chi _{\{v>s\}}) \, ds  
 \,.\end{array} $$
Finally, recalling the equalities \eqref{eqchar} and \eqref{FE}, we obtain
$$\widehat E (v) = \int _0 ^ 1 \widehat E(\1 _{u _s})  \, ds = \int_0 ^ 1 E (u _ s) \, ds\,.$$

\qed

\bigskip

{\it Step 2  (Proof of \eqref{nog} and of last part in Theorem \ref{l:IJ}).} Set for brevity 
$$\mathcal I = \inf \Big \{ E (u) \, :\, u \in {\mathcal A}  \Big \} \qquad \hbox{ and } \qquad  \mathcal J:= \inf \Big \{ \widehat E (v) \, :\, v \in \widehat{\mathcal A} \Big \}\,.$$

For every $u \in \mathcal A$, the function $v = \1 _u$ belongs to $\A $. Therefore, in view of the equality \eqref{FE}, we immediately see that the inequality $\mathcal  I \geq \mathcal J$ is satisfied. 

Conversely, let $v \in \A $ be such that $\widehat E (v) < + \infty$. For such a function $v$,  the slicing formula 
 \eqref{super} holds. Such equality  implies in particular
 that,   for $\mathcal L ^1$-a.e.\ $s \in (0, 1)$, $u_s$ lies in $W ^ {1,p} (\Omega)$; moreover,  since $v = \1_{u _0}$ on $\Gamma _0 \times \R$, it holds $u_s= u_0$ on $\Gamma _0$. Therefore,  for $\mathcal L ^1$-a.e.\ $s \in (0, 1)$, we have $u _ s \in \mathcal A$, which implies $E (u _s) \geq \mathcal I$.  After an integration over $(0, 1)$, by \eqref{super}, we obtain $\widehat E (v) \geq \mathcal I$. 
By the arbitrariness of $v \in \A$, we conclude that $\mathcal J \geq \mathcal I$. 

The equalities \eqref{FE} and \eqref{nog} imply immediately that, if $u \in {\rm argmin } _{\mathcal A} (E)$, then $\1 _u \in  {\rm argmin } _{\widehat{\mathcal A}} (\widehat E)$. Since we know from Proposition \ref{l:existence} that the infimum $\mathcal I$ is finite and attained, we deduce that the same holds true for the infimum $\mathcal J$. 

Finally, if $v \in  {\rm argmin } _{\widehat{\mathcal A}} (\widehat E)$, \eqref{super} and \eqref{nog} imply that $u _s \in {\rm argmin } _{\mathcal A} (E)$  for $\mathcal L ^ 1$-a.e.\ $s \in (0,1)$. In particular this assertion implies  that, in case the primal problem has a finite number of solutions, $v$ must be a convex combination of them as stated in \eqref{cc}. 

\qed

\bigskip

\subsection{Proof of the inequality $\mathcal I \geq \mathcal I ^*$ in Theorem \ref{t:nogap}}\label{proof1}

We are going to prove that,  for every $u \in \mathcal A$  and  $\sigma \in \X$, there  holds  
\begin{equation}\label{disug}
\begin{array}{ll} \displaystyle
\int_\Omega  f(u, \nabla u) \,dx + \int _{\Gamma _1} \gamma (u) \, d \mathcal H ^{N-1} 
& \displaystyle \geq   \int_{G_u}  \sigma \cdot \nu_u \, d\mathcal H^N  + \int _{\Gamma _1} \gamma (u) \, d \mathcal H ^{N-1} \\ \noalign{\medskip}
& \displaystyle =  \int_{G_{u_0}}  \sigma \cdot \nu_{u_0} \, d\mathcal H^N  + \int _{\Gamma _1} \gamma (u_0) \, d \mathcal H ^{N-1}\,.
\end{array}
\end{equation}
Once proved \eqref{disug}, by passing to the infimum over $u \in \mathcal A$ and to the supremum over $\sigma \in \X$  respectively at the left hand side and at the right hand side, we obtain the inequality $\mathcal I \geq \S$. 

\medskip
Let us prove separately the inequality in the first line of \eqref{disug}, and the equality in the second line. 

\medskip

The inequality in the first line of \eqref{disug}  follows simply by recalling \eqref{FE1} and applying 
 Lemma \ref{l:julygen} with $v = \1 _u$:
\begin{equation}\label{bflux} 
\int_\Omega  f(u, \nabla u) \,dx  = \int_{G_u} h _ f ( t, \nu _u) \, d \mathcal H ^N \geq \int_{G_u}  \sigma \cdot \nu_u \, d\mathcal H^N    \,.
\end{equation}

The equality in the second line of \eqref{disug} follows via an integration by parts formula 
that we state separately in the next lemma, since it will be useful again in the sequel.  
It is obtained as an
application of the following generalized divergence theorem, that we recall  from \cite{A} (see also \cite{BoDa}): 
for every $\sigma \in X_1(\Omega \times \R)$ and every  $v \in BV_\infty (\Omega \times \R) \cap  L ^ 1 (\Omega \times \R)$, there holds
\begin{equation}\label{gen_div}
\int_{\Omega \times \R} \sigma \cdot D v + \int _{\Omega \times \R} v \div \sigma \, dx = \int _{\partial \Omega \times \R} ( \sigma\cdot \nu_\Omega )\, v \, d \mathcal H ^ {N}\,.  
\end{equation}
Notice that the boundary integral at the r.h.s.\ is well-defined since the normal trace $\sigma  \cdot \nu _\Omega$ is in $L ^ \infty (\partial \Omega \times \R)$, and the function $v$ is in $L ^ 1 (\partial \Omega \times \R)$ because $v \in BV (\Omega \times \R)$.

\begin{lemma}{\rm (an integration by parts formula)} \label{l:intbyp} For every $\sigma$  in $X_1(\Omega \times \R)$ satisfying \eqref{f:bc} and \eqref{f:div}, and every $v$ in the class $\A$ defined in \eqref{defA}, there holds
\begin{equation}\label{casev}
\int_{\Omega \times \R} \!\! \sigma \cdot Dv \, + \int _{\Gamma _ 1 \times \R} \!\! \gamma '(t)\, ( v-v_0)\, d \mathcal H ^{N-1} \, dt = 
\int_{G_{u_0}} \!\! \sigma \cdot \nu _{u _0} \, d \mathcal H ^{N} + \int _{\Gamma _1}\!\! \gamma (u _0) \, d \mathcal H ^{N-1}\!\!.
\end{equation}
In particular, if $v$ is of the form $v = \1 _u$ for some $u\in \mathcal A$, we obtain
\begin{equation} \int_{G_u} \!\! \sigma \cdot \nu _{u } \, d \mathcal H ^{N}\! + \! \int _{\Gamma _1} \!\!\gamma (u ) \, d \mathcal H ^{N-1} 
=\! \int_{G_{u_0}} \!\!\sigma \cdot \nu _{u _0} \, d \mathcal H ^{N} \! + \! \int _{\Gamma _1} \!\! \gamma (u _0) \, d \mathcal H ^{N-1}. 
\end{equation}
 \end{lemma}

\proof

For every $\sigma$ and $v$  as in the assumptions, we have that the function $v - \1 _{u_0} = (v - v_0) + (v_0 - \1 _{u_0})$ is in $BV_\infty (\Omega \times \R) \cap L ^ 1 (\Omega \times \R)$, and $\sigma$ is in $X_1(\Omega \times \R)$. 
Therefore, we are in a position to apply the generalized Gauss-Green formula \eqref{gen_div}. 
Exploiting also the condition $\div \sigma = 0$ in $\Omega \times \R$, we obtain 
$$ \begin{array}{ll} \displaystyle\int_{\Omega \times \R} \sigma \cdot (Dv - D \1_{u _0}) &\displaystyle = 
\int_{\partial \Omega \times \R} (\sigma ^x \cdot \nu _\Omega)\, (v- \1_{u _0}) 
= -
 \int _{\Gamma _ 1 \times \R} \gamma '(t)\, ( v-\1_{u _0})\, d \mathcal H ^{N-1} \, dt \\
 \noalign{\medskip}
& \displaystyle = 
 - \int _{\Gamma _ 1 \times \R} \gamma '(t)\,  [( v-v_0)  -  (\1_{u _0}-v_0)]\, d \mathcal H ^{N-1} \, dt \,.
 \end{array} 
$$
Hence, 
$$\begin{array}{ll}& \displaystyle
\int_{\Omega \times \R} \sigma \cdot Dv + \int _{\Gamma _ 1 \times \R} \gamma '(t)\, ( v-v_0)\, d \mathcal H ^{N-1} \, dt 
\\ \noalign{\medskip}
= 
& \displaystyle \int_{\Omega \times \R} \sigma \cdot D\1_{u _0} + \int _{\Gamma _ 1 \times \R} \gamma '(t)\, ( \1_{u _0}-v_0)\, d \mathcal H ^{N-1} \, dt 
\\ 
\noalign{\medskip}
=   &\displaystyle 
\int_{G_{u_0}} \sigma \cdot \nu _{u _0} \, d \mathcal H ^ {n-1}+ \int _{\Gamma _1} \gamma (u _0) \, d \mathcal H ^{N-1} \,.
\end{array}$$

Notice that in the last equality we have used the identity 
$\int _\R \gamma' (t) (\1_{u_0} - v_0 )  \, dt   = \gamma (u_0)$, already shown in the proof of Theorem \ref{l:IJ} ({\it cf.}\ equation \eqref{iduseful}). 
\qed

\bigskip
We have thus completed the proof of \eqref{disug} and hence of the inequality $\mathcal I \geq \mathcal I ^*$ in Theorem \ref{t:nogap}. 

\section{Optimality conditions and min-max formulation}\label{secminmax}
Out next goal is to provide necessary and sufficient conditions for optimality:

\begin{theorem}{\rm (geometric optimality condition)} \label{t:calib} Let $u \in \mathcal A$ and $\sigma \in \X$.  
Then $u$ is a solution to the primal problem problem $({\mathcal P})$ in \eqref{defI}  and $\sigma$ is a solution to the dual problem $({\mathcal P}^*)$ in \eqref{defI*}  if and only if 
\begin{equation}\label{f:fun=}
h _ f (t, \nu _u) = \sigma \cdot \nu _u \qquad \mathcal H ^N\hbox{-a.e.\ on } G_u\,.
\end{equation}
In this case, we say that $\sigma$ is a {\rm calibration} for $u$. 
\end{theorem}

{\it Proof}. Assume that $u \in \mathcal A$ and $\sigma \in \X$ satisfy \eqref{f:fun=}. 
By using in the order the definition of $\S$, Lemma \ref{l:intbyp}, condition \eqref{f:fun=}, the equality \eqref{FE1}, and the definition of $\mathcal I$, we obtain
$$\begin{array}{ll}
\S & \displaystyle \geq \int_{G_{u_0}} \!\!\sigma \cdot \nu _{u _0} \, d \mathcal H ^{N} \! + \! \int _{\Gamma _1} \!\! \gamma (u _0) \, d \mathcal H ^{N-1}  =  \int_{G_u} \!\! \sigma \cdot \nu _{u } \, d \mathcal H ^{N}\! + \! \int _{\Gamma _1} \!\!\gamma (u ) \, d \mathcal H ^{N-1} 
\\ \noalign{\medskip}  
& \displaystyle =  \int_{G_u}h _ f (t, \nu _u) \, d \mathcal H ^N + \! \int _{\Gamma _1} \!\!\gamma (u ) \, d \mathcal H ^{N-1} =  \int_{\Omega} \!\!f (u, \nabla u)  \, dx \! + \! \int _{\Gamma _1} \!\!\gamma (u ) \, d \mathcal H ^{N-1} \geq \mathcal I \,.
\end{array}$$
Since we know from Theorem \ref{t:nogap} that $\mathcal I = \S$, we infer that all the inequalities above hold as equalities, which means in particular that $u$ and $\sigma$ are optimal respectively for the primal and the dual problem. 

Assume that $u \in \mathcal A$ and $\sigma \in \X$ are optimal respectively for the primal and the dual problem. 
By using in the order Lemma \ref{l:intbyp}, the optimality of $\sigma$, Theorem \ref{t:nogap}, the optimality of $u$, and the equality \eqref{FE1}, we obtain
$$\begin{array}{ll}
\displaystyle \int_{G_u} \!\! \sigma \cdot \nu _{u } \, d \mathcal H ^{N}\! + \! \int _{\Gamma _1} \!\!\gamma (u ) \, d \mathcal H ^{N-1} 
&\displaystyle = \int_{G_{u_0}} \!\!\sigma \cdot \nu _{u _0} \, d \mathcal H ^{N} \! + \! \int _{\Gamma _1} \!\! \gamma (u _0) \, d \mathcal H ^{N-1}
 = \S  
\\ \noalign{\medskip}  
& \displaystyle  = \mathcal I  =  \int_{\Omega} \!\!f (u, \nabla u)  \, dx \! + \! \int _{\Gamma _1} \!\!\gamma (u ) \, d \mathcal H ^{N-1}
\\ \noalign{\medskip}  
& \displaystyle  =  \int_{G_u} h _f ( t, \nu _u ) \, d \mathcal H ^N \! + \! \int _{\Gamma _1} \!\!\gamma (u ) \, d \mathcal H ^{N-1}
 \,.
\end{array}$$
We infer that $\int_{G_u} \!\! \sigma \cdot \nu _{u } \, d \mathcal H ^{N} = \int_{G_u} h _f ( t, \nu _u ) \, d \mathcal H ^N $. In turn, recalling the inequality \eqref{f:fun} in Remark \ref{r:easyine2}, this implies \eqref{f:fun=}. 
\qed
\bigskip

From a practical point of view, in order to construct a calibration, it is useful to rephrase condition  \eqref{f:fun=} more explicitly as done in the next result. 

\begin{corollary}[user's form of optimality conditions] \label{p:september}
Let $u \in \mathcal A$ and $\sigma \in \X$, with $\sigma$ continuous on $\Omega \times (\R \setminus D)$. Then condition \eqref{f:fun=} is satisfied if and only if there holds
\begin{eqnarray}
& \sigma ^ x (x, u (x)) \in  \partial_z f (u (x), \nabla u (x)) \ \text{ for } \mathcal L ^ {N}\text{-a.e.}\ x \in u ^ {-1} (\R \setminus D) ;  &\label{o1}
\\
\noalign{\smallskip}
& \sigma ^ t (x, u (x)) = f ^*_z (u (x), \sigma ^x(x, u (x))) \ \text{ for } \mathcal L ^ {N}\text{-a.e.}\ x \in u ^ {-1} (\R \setminus D)  ;  &\label{o2}
\\
\noalign{\smallskip}
& \sigma ^ t (x, t) = - f (t, 0) \quad \forall\, t \in \R  \text{ and for } \mathcal L ^ {N}\text{-a.e.}\ x \in \{u = t \}\, . & \label{o3}
\end{eqnarray}
(Note that the set of values $t \in \R$ such that $\mathcal L ^N ( \{ u = t \} )>0$ is at most countable.) 
%
\end{corollary}

\proof  By \eqref{o1}-\eqref{o2}, we infer that the following equality  is satisfied $\mathcal H ^N$-a.e.\ on $G_u \cap [\Omega \times (\R \setminus D)]$:
$$ \begin{array}{ll} \sigma  \cdot \nu _u & \displaystyle = \nu _u ^ x \cdot \sigma ^x + \nu _u ^t \, f ^* _z (u , \sigma ^ x) \\ \noalign{\medskip} & = - \nu _u ^ t \big [ - \frac{\nu _u ^ x}{\nu _u ^ t } \cdot \sigma ^ x - f ^* _z (u, \sigma ^x) \big ]
\\ \noalign{\medskip} &= - \nu _u ^ t \big [\nabla u \cdot \sigma ^x - f ^* _z (u, \sigma ^x) \big ] \\ \noalign{\medskip} & = - \nu _u ^ t f (u, \nabla u) \\ \noalign{\medskip} &\displaystyle= 
- \nu _u ^ t \, f (u,  - \frac{\nu _u ^ x}{\nu _u ^ t } ) = h _f (t, \nu _u)\,. \end{array} $$ 
On the other hand, by \eqref{o3}, $\mathcal H ^N$-a.e.\ on $G_u \cap [\Omega \times D]$ we have
$$ \sigma \cdot (- e _{N+1} ) = - \sigma ^ t (x, t) = f (t, 0) = h _ f (t, - e _{N+1}) \,.
 $$
Recalling \eqref{nuvert} with $v = \1 _u$, we conclude that \eqref{f:fun=} is fulfilled. 

Conversely, assume that \eqref{f:fun=} holds true. 

Since $\sigma$ satisfies \eqref{s3} and is assumed to be continuous on $\Omega \times (\R \setminus D)$, the following chain of inequalities is satisfied 
$\mathcal H ^N$-a.e.\ on $G_u \cap [\Omega \times (\R \setminus D)]$:
$$ \begin{array}{ll} h _ f (t, \nu _u ) & = \sigma  \cdot \nu _u  \displaystyle \leq \nu _u ^ x \cdot \sigma ^x + \nu _u ^t \, f ^* _z (u , \sigma ^ x)  \\ \noalign{\medskip} & \displaystyle = - \nu _u ^ t \big [ - \frac{\nu _u ^ x}{\nu _u ^ t } \cdot \sigma ^ x - f ^* _z (u, \sigma ^x) \big ]
\\ \noalign{\medskip} &\displaystyle = - \nu _u ^ t \big [\nabla u \cdot \sigma ^x - f ^* _z (u, \sigma ^x) \big ]\\ \noalign{\medskip} &\displaystyle \leq - \nu _u ^ t f (u, \nabla u) \\ \noalign{\medskip} & = 
- \nu _u ^ t \, f (u,  - \frac{\nu _u ^ x}{\nu _u ^ t } ) = h _f (t, \nu _u)\,. \end{array} $$ 
We deduce that the two inequalities appearing in the chain are actually equalities, which yields \eqref{o1}-\eqref{o2}. 

On the other hand, 
since $\sigma$ satisfies \eqref{s4} on $\Omega \times \R$ ({cf.}\ Remark \ref{r:june}), $\mathcal H ^N$-a.e.\ on $G_u \cap [\Omega \times \R]$ we have 
$$ h _ f (t, - e _{N+1})  =  \sigma \cdot (- e _{N+1} ) = - \sigma ^ t (x, t) \leq f (t, 0) = h _ f (t, - e _{N+1}) \,.
 $$
We conclude that the inequality appearing in the line above holds with equality sign, which yields \eqref{o3}.  \qed
\bigskip

\begin{remark}\label{p:convexcali}
In the case when $f$ is differentiable and convex in $(t, z)$ and $\gamma ' \equiv c$, it is easy to construct an explicit calibration for a given solution $\o u $ to problem $(\mathcal P)$. Indeed, denoting by  $\overline \sigma$   a solution to the classical dual problem  ({\it cf.}\ Remark \ref{r:convexdual}),  we claim that
the field $\sigma$ defined on $\Omega \times \R$ by 
 \begin{equation}\label{convexcali}
 \begin{cases}
 \sigma ^ x (x, t ) = \o \sigma (x) & \\ \noalign{\medskip}
 \sigma ^ t (x, t) = f ^ * _z (\o u, \o \sigma  ) 
 -  (\div \o \sigma ) \big ( t - \o u (x) \big )\,, & 
 \end{cases}
\end{equation} 
is a calibration for $\o u$, provided it is continuous on $\Omega \times \R$.

%
%
%
%
%
Namely, 
by classical duality, $\o u$ and $\o \sigma$ satisfy the optimality conditions
\begin{eqnarray}
& \o \sigma = \partial _z f (\o u, \nabla \o u)\, , \quad \div \o \sigma = \partial _t f (\o u, \nabla \o u) \qquad 
\mathcal L ^ N\hbox{-a.e. in } \Omega
& \label{optcond}
\\ 
& \o \sigma \cdot \nu _\Omega = - \gamma' (\overline u) \qquad \mathcal H ^ {N-1}\hbox{-a.e. on } \Gamma_1\,. & \label{optcondbord}
\end{eqnarray}
In view of \eqref{optcond} and of the continuity assumption made on $\sigma$, Corollary \ref{p:september} (applied with $D = \emptyset$) ensures that $\sigma$ is a calibration for $\o u$, provided we show that $\sigma \in \X$. 

It is immediate to verify that $\sigma$ satisfies \eqref{f:div}. By \eqref{optcondbord} and the assumption $\gamma ' \equiv c$, it satisfies also \eqref{f:bc}. It only remains to check  \eqref{s3}, namely
$$ f ^ * _z (\o u, \o \sigma  ) 
 - ( \div \o \sigma ) \big ( t - \o u (x) \big ) \geq f _ z ^* ( t, \o \sigma) \qquad \mathcal L ^ {N+1}\hbox{-a.e.\ on } \Omega \times \R\,,$$
or equivalently
$$ f ^ * _z (\o u, \o \sigma  )  \geq  \sup _{z \in \R ^N} \big [\o \sigma \cdot z - f ( t , z) \big ]
+   ( \div \o \sigma ) \big ( t - \o u (x) \big )
  \qquad \mathcal L ^ {N+1}\hbox{-a.e.\ on } \Omega \times \R\, .$$
  In turn the latter inequality is satisfied provided
 \begin{equation}\label{penultima}
\begin{array}{ll} f ^ * _z (\o u, \o \sigma  )  & \geq  \displaystyle{\sup _{(t,z) \in \R ^{N+1}} \big [(\div \o \sigma, \o \sigma ) \cdot (t, z) - f ( t , z) \big ]
- \o u(x) \div \o \sigma } 
 \\ \noalign{\medskip} & =  f ^* (\div \o \sigma, \o \sigma ) 
- \o u(x) \div \o \sigma  
 \qquad \mathcal L ^N\hbox{-a.e.\ on } \Omega \, , \end{array}
  \end{equation}
where $f ^*$ denotes the global Fenchel conjugate of $f$ with respect to the pair $(t, z)$.

 Now, by the two equations in (\ref{optcond}), we have that $(\div \o \sigma, \o \sigma)$ satisfy the Fenchel equality 
 $$f ^* (\div \o \sigma, \o \sigma ) + f (\o u,  \nabla \o u) = \o u (x) \div \o \sigma + \nabla \o u (x) \cdot \o \sigma\qquad \mathcal L ^N\hbox{-a.e.\ on } \Omega \,.$$
 Inserting this identity into (\ref{penultima}), we are reduced to
 $$f ^ * _z (\o u, \o \sigma  ) \geq  \nabla \o u (x) \cdot \o \sigma- f (\o u,  \nabla \o u) \qquad \mathcal L ^N\hbox{-a.e.\ on } \Omega\,,$$
 which is satisfied by definition of $f ^ * _z$ (and actually holds as an equality since $\o \sigma = \partial _z f (\o u, \nabla \o u )$).  
\end{remark}

\bigskip
Hereafter we give a min-max formulation of our duality result. 
For every pair $(v, \sigma)$, with $v \in BV _\infty (\Omega \times \R)$ and $\sigma \in X _ 1 (\Omega \times \R)$, we introduce the Lagrangian
\begin{equation}\label{f:defL}
L (v, \sigma ):=  \int _{\Omega \times \R} \sigma \cdot Dv + \int _{\Gamma _1 \times \R }  \gamma ' (t) (v-v_0) \, d \mathcal H ^ {N-1}\,.
\end{equation}


\begin{theorem}[saddle point]\label{minmax}
There holds
$$\mathcal I = \inf _{v \in \widehat{\mathcal A}} \ \sup _{\sigma \in \mathcal K } \ L (v, \sigma) =
\sup _{\sigma \in \mathcal K } \  \inf _{v \in \widehat{\mathcal A}} \ L (v, \sigma) =  \mathcal I ^* \,.$$
Moreover, a pair  $(\overline v, \overline \sigma)$ is optimal for the convexified infimum problem $\inf \big \{ \widehat E (v) \, :\, v \in \widehat{\mathcal A} \}$ 
and for the dual problem $( \mathcal P ^*)$ in \eqref{defI*} if and only if it is a saddle point for $L$, namely
$$L ( \overline v, \sigma) \leq L (\overline v, \overline \sigma ) \leq L (v, \overline \sigma ) \qquad \forall (v, \sigma) \in  
\widehat{\mathcal A } \times \mathcal K  \,.$$
\end{theorem}

\begin{remark}
(i) Notice that, since the class $\mathcal A$ is not weakly compact, the equality $\mathcal I = \mathcal I ^*$ already established in Theorem \ref{t:nogap} 
cannot be deduced by applying an inf-sup commutation argument to the bivariate Lagrangian $L$ over the product space $\mathcal A \times \mathcal K$. 

(ii) We emphasize that  the class $\mathcal K$ appearing in the saddle point problem does not include the divergence free condition. In fact, such condition is handled by duality, through the use of the variable $v$ seen as a Lagrange multiplier. (In analogy with fluid dynamic, one may think of $\sigma$ as the speed of an incompressible fluid, and of $v$ as its pressure).  

\end{remark}

\bigskip

{\it Proof  of Theorem \ref{minmax}}.  Thanks to the equality \eqref{nog} in Theorem \ref{l:IJ}, 
recalling the definitions \eqref{f:F}, \eqref{f:H}, \eqref{f:ell}, and \eqref{f:defL} of $\widehat E$, $H$, $\ell$,  and $L$, 
we obtain
$$\begin{array}{ll}
\mathcal I  & \displaystyle =  \inf \big \{ \widehat E (v) \, :\, v \in \widehat{\mathcal A} \} = \inf \big \{ H (v) + \ell (v)  \, :\, v \in \widehat{\mathcal A} \}   
\\ \noalign{\medskip}
& \displaystyle  =  \inf _{v \in \widehat{\mathcal A}} \ \sup _{\sigma \in \mathcal K } \Big \{  \int_{\Omega \times \R} \sigma \cdot Dv + \int _{\Gamma_1 \times \R } \gamma' (t) \, (v- v_0) \, d \mathcal H ^ {N-1} \, dt  \Big \} =  \inf _{v \in \widehat{\mathcal A}} \ \sup _{\sigma \in \mathcal K } L (v, \sigma)\,.
\end{array}$$

Since we know from Theorem \ref{t:nogap} that $\mathcal I = \mathcal I ^*$, in order to complete the proof it remains to show that $\mathcal I ^* = 
\sup _{\sigma \in \mathcal K } \  \inf _{v \in \widehat{\mathcal A}} \ L (v, \sigma)$.  To that aim let us show that, for every $\sigma \in \mathcal K$, it holds
\begin{equation}\label{f:2mai}\inf _{v \in \widehat{\mathcal A}} \ L (v, \sigma) = \begin{cases}
- \infty & \text{ if } \sigma \not \in \mathcal B
\\ \noalign{\medskip}
\displaystyle \int_{G_{u_0}} \sigma \cdot \nu _{u _0} \, d \mathcal H ^N + \int _{\Gamma _1} \gamma (u_0) \, d \mathcal H ^ {N-1} & \text{ if } \sigma \in \mathcal B\,.
\end{cases}
\end{equation}
Indeed, the Lagrangian $L(v, \sigma)$ can be rewritten as $$\begin{array}{ll}
L ( v, \sigma) & \displaystyle =  \int_{\Omega \times \R} \sigma \cdot (Dv - D \1 _{ u _0} ) +  \int_{G_{u_0}} \sigma \cdot \nu _{u _0} \, d \mathcal H ^N  +  \int _{\Gamma _1 \times \R }  \gamma ' (t) (v-v_0) \, d \mathcal H ^ {N-1}\,.\end{array}$$
Then, by exploiting the generalized Gauss-Green formula \eqref{gen_div}, we get
$$\begin{array}{ll}
L ( v, \sigma) & \displaystyle
 =  - \int_{\Omega \times \R} \div \sigma \cdot (v -  \1 _{ u _0} ) +  \int_{\Gamma_0 \times \R} (v- \1 _{u_0}) (\sigma ^ x \cdot \nu_ \Omega) \, d \mathcal H ^ {N-1} \, dt \\ \noalign{\bigskip}
&  \displaystyle
+ \int_{\Gamma_1 \times \R} (v- v_0) (\gamma ' (t) + \sigma ^ x \cdot \nu _\Omega)  \, d \mathcal H ^ {N-1} \, dt \\ \noalign{\bigskip}
&  \displaystyle
+    \int_{\Gamma_1 \times \R} (v_0- \1_{u_0}) 
(\sigma ^ x \cdot \nu _\Omega )\, d \mathcal H ^ {N-1} \, dt 
 +  
 \int_{G_{u_0}} \sigma \cdot \nu _{u _0} \, d \mathcal H ^N  
\end{array}
$$
Now, by taking $v \in \widehat{\mathcal A }$ of the form $v = \1 _{u _0} + \varphi$, with  $\varphi \in \mathcal D (\Omega \times \R)$, we obtain that $\inf _{v \in \widehat{\mathcal A}} \ L (v, \sigma)$ cannot be finite unless $\div \sigma = 0$ in $\Omega \times \R$. Next, by taking
$v = v_0 + \varphi$, with 
$\varphi \in \mathcal D (\overline \Omega \times \R)$ such that $\varphi = 0$ on $\Gamma _0 \times \R$, we see that the normal trace of $\sigma$ must 
agree with $- \gamma ' (t)$ on $\Gamma _ 1 \times \R$.  
We conclude that \eqref{f:2mai} is true by recalling \eqref{iduseful}. 

The last part of the statement is a standard equivalence in min-max theory (see for instance  \cite{EkTe}).

\qed

As mentioned in Remark \ref{r:bounded}, whenever the solutions to the primal problem are bounded, we can settle our duality theory on a bounded set of the form $\Omega  \times [m, M]$.    

For a given $u _0 \in W ^ {1, p} (\Omega; [m, M]$, we denote by  $\mathcal I (m, M)$ and $\mathcal I ^* (m, M)$ respectively the infimum of the primal problem $(\mathcal P)$   and the supremum of the dual problem $(\mathcal P^*)$  over the classes $\mathcal A (m, M)$  and $\mathcal B (m, M)$ introduced in Remark \ref{r:bounded}. Then we set
$$\begin{array} {ll}
& \widehat {\mathcal A} (m, M) : = \Big \{ v \in \mathcal A \ :\  v (x, t) = 1 \text{ for } t <m \, , \ v (x, t) = 0 \text{ for } t >M \Big\} 
\\ \noalign{\smallskip}
& {\mathcal K } (m, M) := \Big \{ \sigma \in X _1 (\Omega \times (m, M)) \text{ satisfying \eqref{s31}-\eqref{s41}} \Big \} \, , 
\end{array}
$$
Accordingly, the Lagrangian $L$ must be now intended as
\begin{equation}\label{f:defL}
L (v, \sigma ):=  \int _{\Omega \times [m, M]} \sigma \cdot Dv + \int _{\Gamma _1 \times [m, M] }  \gamma ' (t) (v-v_0) \, d \mathcal H ^ {N-1}\,.
\end{equation}

\begin{remark}
Note that in \eqref{f:defL} the first integral may have a non vanishing contribution on the horizontal part of the boundary (namely the set $\Omega \times \{m, M\}$), in case the function $v$ has a jump on such interfaces. More precisely, we have:
$$\int_{\Omega\times [m, M]} \sigma \cdot Dv = \int _{\Omega\times (m, M)} \sigma \cdot Dv + \int _\Omega \big [ \sigma ^ t (x, M) (0- v (x, M ^-)) + \sigma ^ t (x, m) (v (x, m ^+) - 1 ) \big ] \,,$$
 being $v (x, m^+)$ and $v (x, M^-)$ respectively the traces of $v$ on $\Omega \times \{m\}$ and $\Omega \times \{M \}$.  
\end{remark}

\bigskip
We can now reformulate the following variant of Theorems \ref{t:nogap} and \ref{minmax}.  

\begin{proposition}\label{c:minmax} With the above notation, there holds:
$$\mathcal I (m , M) = \inf _{v \in \widehat{\mathcal A}(m, M)}  \sup _{\sigma \in \mathcal K (m, M)}  L (v, \sigma) =
\sup _{\sigma \in \mathcal K(m , M) }   \inf _{v \in \widehat{\mathcal A}(m, M)}  L (v, \sigma) =  \mathcal I ^* (m, M) \,.$$
Moreover, a pair  $(\overline v, \overline \sigma)$ is optimal for the infimum problem $\inf \big \{ \widehat E (v) \, :\, v \in \widehat{\mathcal A}(m, M) \}$ 
and for the dual problem $( \mathcal P ^*)$ settled over $\mathcal B(m, M)$ if and only if 
$$L ( \overline v, \sigma) \leq L (\overline v, \overline \sigma ) \leq L (v, \overline \sigma ) \qquad \forall (v, \sigma) \in  
\widehat{\mathcal A } (m, M) \times \mathcal K(m, M)  \,.$$
\end{proposition} 

\proof  The statement can be proved in the analogous way as done for Theorem \ref{minmax}, taking into account that, by following the same proof as in Theorem \ref{t:nogap}, one can check that  $\mathcal I (m , M) = \mathcal I ^* (m, M)$.  \qed

 
\section {Application to a free boundary problem}\label{secexamples}


\subsection{Description of the problem.} In this section we illustrate the application of our method to the free boundary problem
\begin{equation}\label{altcaff}
\mathcal I (\Omega, \lambda):= \inf \Big \{ \int _\Omega \frac{1}{2} |\nabla u| ^ 2 + \lambda  \big |\{ u >0 \} \big | \ :\ u \in W ^ {1, 2} (\Omega) \, ,\ u = 1 \hbox{ on } \partial \Omega \Big \}\,,
\end{equation}
which has been firstly considered in the pioneering paper \cite{AlCa}. 

\smallskip
The free boundary in the minimization problem \eqref{altcaff} is the frontier of the zero level set $E:= \{ u =0\}$. Actually, the infimum $\mathcal I(\Omega, \lambda)$ can be recast by solving the shape optimization problem
$$\inf _E  \Big \{ \int _\Omega \frac{1}{2} |\nabla u_E| ^ 2 + \lambda  \big |\Omega \setminus E | \Big \} \, , $$ 
being $u _E$ the solution to 
$$\begin{cases}
\Delta u = 0 & \text{ in  } \Omega \setminus E
\\
u = 0 & \text{ in } E
\\
u = 1 & \text{ on  } \partial \Omega\,.
\end{cases} 
$$
Such problem falls in our setting by choosing
$$f (t,  z) = \frac{1 }{2} |z| ^ 2 + \lambda \chi _{(0, + \infty)} ( t) \ , \qquad \ (\Gamma _0, \Gamma _ 1)  = (\partial \Omega, \emptyset)\, , \qquad u _0 \equiv 1 \,.$$ 
Notice that the function $f$ satisfies the standing assumptions, and in particular the discontinuity set $D$ appearing in \eqref{H1bis} is given by
$\{t= 0 \}$.

Then, according to Theorem \ref{t:nogap}, we have $\mathcal I (\Omega, \lambda) = \mathcal I ^* (\Omega, \lambda)$. As disclosed in the Introduction, the dual problem reads:
\begin{equation}\label{altcaffdual}
\mathcal I^* (\Omega, \lambda):= \sup \Big \{ - \int _\Omega \sigma ^ t (x, 1) \, dx  \ :\ \sigma \in \mathcal B \Big \}\,, 
\end{equation} 
where $\mathcal B$ is the class of bounded divergence free vector field on $\Omega \times \R$ satisfying the constraints
\begin{equation}\label{consigma}
\sigma ^ t (x, t) + \lambda \geq \frac{1}{2} |\sigma ^ x ( x, t) | ^ 2   \quad \text{a.e. on } \Omega \times \R \, , \qquad \sigma ^ t (x, 0) \geq 0  \text{ a.e. on } \Omega\,.
\end{equation}

It is easy to check that any solution $u\in W ^ {1, 2} (\Omega)$ to problem \eqref{altcaff}  takes values in $[0, 1]$. Therefore, according Remark \ref{r:bounded}, we can work on the bounded subset $\Omega \times [0, 1]$. Then,  in virtue of Proposition \ref{c:minmax}, searching for an optimal pair $(\overline u, \overline \sigma)$ amounts to find a saddle point for the bivariate functional 
 \begin{equation}\label{Lclosed}
 \inf_{v \in \widehat {\mathcal A} } \ \sup _{ \sigma \in \mathcal K } \ 
\int_{\Omega \times[0, 1]}   \sigma \cdot Dv \, , 
 \end{equation}
 with
$$\begin{array} {ll}
& \widehat {\mathcal A} = \Big \{ v \in BV _\infty (\Omega \times\R ) \ :\  v  = 1 \text{ for } t <0\, , \ v  = 0 \text{ for } t >1\, , \ v = 1  \text{ on } \partial \Omega \times [0, 1]\Big\} 
\\ \noalign{\medskip}
& {\mathcal K } = \Big \{ \sigma \in X _1 (\Omega \times (0, 1)) 
\ :\   \sigma ^ t + \lambda \geq \frac{1}{2} |\sigma ^ x| ^ 2 \ \text{a.e. on } \Omega \times (0, 1) \, ,  \ \sigma ^ t (\cdot, 0) \geq 0 \text{ a.e. on } \Omega \Big \} \, . 
\end{array}
$$  
Notice carefully that the integration domain in \eqref{Lclosed} is the product of $\Omega$ times the {\it closed} interval $[0, 1]$. Actually,  minimizing over $\widehat {\mathcal A}$ the functional $v \mapsto \int_{\Omega \times[0, 1]}   \sigma \cdot Dv$ appearing in \eqref{Lclosed} is equivalent  to minimizing over the space of functions $v \in BV (\Omega \times (0, 1))$  satisfying the boundary condition $v = 1$  on  $\partial \Omega \times [0, 1]$ the functional
 $$v \mapsto \int_{\Omega \times(0, 1)}   \sigma \cdot Dv  + \int _{\Omega} \big  [ \sigma ^ t (x, 0)(v (x, 0^+) -1 )- \sigma ^ t (x, 1) v (x, 1^-)  \big ] \, dx $$
   being $v (x, 0^+)$ and $v (x, 1^-)$ respectively the traces of $v$ on $\Omega \times \{0\}$ and $\Omega \times \{1 \}$.  

Before proceeding  to solve the min-max problem \eqref{Lclosed} let us recall that,
if $(\overline v, \overline \sigma)$ is an optimal pair,  
 the function $\overline v$ should be a step function. Indeed, 
we expect that the primal problem \eqref{altcaff} admits only one or at most a finite number of solutions. Then, by virtue of \eqref{cc}, the function $\overline v$ will  take only the values $0$ and $1$ in case of a unique solution, or a finite number of values in $[0 , 1]$  in case of multiple solutions.

\subsection{Numerical algorithms.}
In order to solve the saddle point problem \eqref{Lclosed}, we adopt two different numerical schemes.

\smallskip
The first one  is a primal-dual algorithm which generalizes a classical method of Arrow-Hurwicz \cite{arrow}, which we took from \cite{PCBC1} (see also \cite{PCBC2}).  We choose an initial point $(v_0, \sigma _0) \in \widehat {\mathcal A} \times \mathcal K$ and two positive time steps $\alpha$, $\beta$. Then, for each $n \in \N$, denoting by $h$ the size parameter of a cartesian grid in $\R ^ {N+1}$,  we let
\begin{equation}\label{syst1}
\begin{cases}
\sigma_{n+1}^h = \Pi^h _K(\sigma_n^h + \alpha \nabla^h \overline{v}_n^h) 
\\ \noalign{\smallskip} 
v_{n+1}^h = v_n^h + \beta \div^h (\sigma^h_{n+1})\\ \noalign{\smallskip} 
\overline{v}^h_{n+1} = 2 v^h_{n+1} - v^h_n\, , 
\end{cases}
\end{equation}
where $\Pi ^h _K$ is  a suitable projection operator associated with the convex constraint $K (t)$ . 
The convergence for system \eqref{syst1} requires that the stringent condition $\alpha \beta \leq c_h^2 \leq 1$ is satisfied, where $c _h$  equals ${2 \sqrt N}/{h}$ (namely the norm of the discretized gradient operator).

The computational cost  in terms of the mesh size $h$ can be shown to be of order  $\frac{1}{h ^ {N+2}}$. 
\smallskip
The second scheme is inspired from the projection method for Navier-Stokes system, 
in which a $L ^2$-orthogonal projection is performed on the space of divergence free field
 (in this analogy, $\sigma$ and $v$ represent respectively the speed and the pressure of the fluid). 
 Roughly, in our case we start from the reformulation of  problem \eqref{Lclosed} as
  \begin{equation}\label{Lclosed2}
 \inf_{p \in {\mathcal C} } \ \sup _{ \sigma \in \mathcal K } \ 
\int_{\Omega \times[0, 1]}   \sigma \cdot p \, \qquad \text{ with } \mathcal C:= \{ Dv \ :\ v \in \widehat {\mathcal A}\}.
 \end{equation}
Then we replace the second equation in \eqref{syst1} by  
$$p_{n+1}^h =  \Pi ^h _ {\mathcal C}  (p_n^h - \beta  \sigma^h_{n+1})\, , $$
where $\Pi ^h_ {\mathcal C}$  is the $L ^2$-orthogonal projector on the convex set $\mathcal C$. 
Denoting by $(\Delta ^h) ^ {-1}$ the discretization of the inverse Dirichlet-Neumann Laplacian operator 
which associates to a function $\varphi$ the solution $w$ to
$$\Delta w = \varphi\, , \qquad w = 0 \text{ on } \partial \Omega \times (0, 1) \, , \qquad \frac{\partial w}{\partial n} = 0 
\text{ on } \Omega \times \{0, 1\} \, , $$
we are led to the following semi-implicit algorithm \begin{equation}\label{syst2}
\begin{cases}
\sigma_{n+1}^h = \Pi^h _K(\sigma_n^h + \alpha \nabla^h \overline{v}_n^h) 
\\ \noalign{\smallskip} 
v_{n+1}^h = v_n^h - \beta (\Delta ^h) ^ {-1} (\div ^h(\sigma^h_{n+1}))\\ \noalign{\smallskip} 
\overline{v}^h_{n+1} = 2 v^h_{n+1} - v^h_n\,.
\end{cases}
\end{equation}
Notice that  \eqref{syst2} differs from  \eqref{syst1} just in the term $ -(\Delta ^h) ^ {-1} (\div ^ h (\sigma^h_{n+1}))$ which replaces 
${\div}^h( \sigma^h_{n+1})$.

The theoretical convergence of this second algorithm can be proved under the condition $\alpha\beta  \leq 1$, which is independent from both the mesh side and the space dimension.  
Moreover, in this case the convergence  occurs after a relatively small number of iterations. 
In fact, 
the inverse Laplacian computation is the most costly (in particular for $\Omega \subset\R ^ 2$ when one works in $\R ^ 3$), and the computational cost depends highly on the solver used for the inverse Laplace operator;  if one uses
a multigrid or a FFT solver, it can be of order  $\frac{1}{h ^ {N+1}\log h}$.

\subsection{Some simulations in case $N=1$.} 

When the open set $\Omega$ is an interval $(0, a)$ of the real line, we can solve explicitly the primal problem, which reads
\begin{equation}\label{1Dpb}
\mathcal I (a, \lambda) := \inf \Big \{ \int _0 ^ a  \frac{| u'| ^ 2}{2}  + \lambda \big | \{ u \neq 0 \} \big |  \, dt \ :\ u \in W ^ {1, 2} (0, a)\, , \ u (0) = u (a) = 1 \Big \}\,.
\end{equation}
The Euler-Lagrange equation  written in the integrated conservation law form reads
\begin{equation}\label{EL}\frac{1}{2} |u'| ^ 2 - \lambda \chi _{\{ u \neq 0 \} } = C\,.\end{equation}
Two cases may occur, according to whether the measure of the level set  $\{ u = 0 \}$ is null or strictly positive. 
In the first case, the solution is the constant function equal to $1$ on $(0, a)$, with cost equal to $\lambda a$.  
In the second case, the constant $C$ in   \eqref{EL} equals zero, so that $u' \in \{ 0 , \pm \sqrt{2 \lambda }  \} $. 
Setting $E ^ {\pm}  \:= \{ x \in (0, a) \, :\, u' = \pm  \sqrt{2 \lambda }  \} $, since $\int_0 ^ a u' = 0$, there holds $|E^+ | = |E^-|$, and the cost is $4 \lambda |E^-|$. On the other hand, since $u (0)= 1$ and $u$ reaches the level zero, we have the lower bound $|E^-| \geq {1}/{\sqrt{2 \lambda}}$. Therefore, such a function $u$ can be a minimizer only if $a \geq 2 \sqrt{{2}/{\lambda}}$, and in this case the minimal cost is larger than or equal to $2 \sqrt{2 \lambda}$, with equality if $E ^- = (0, {1}/{\sqrt{2 \lambda}} )$, 
$E^+ = (h - {1}/{\sqrt {2 \lambda}}, a )$.  To summarize, we have
$\mathcal I (\lambda, a) = \min \{ \lambda a, 2 \sqrt { 2 \lambda}\}$,
and

(i)
 for $a \in (0,
 2 \sqrt{ \frac{2}{\lambda}} ]$, the unique solution is $\overline u_1 \equiv 1$; 

\smallskip
(ii) for $a > 2 \sqrt{ \frac{2}{\lambda}}  $, the unique solution is
 $$\overline u _2 (x) = \begin{cases}
 - \sqrt { 2 \lambda} x + 1 & \text{ if } x \in [0, \frac{1}{\sqrt{2 \lambda}}]
 \\
 0 &\text{ if } x \in [ \frac{1}{\sqrt{2 \lambda}},a - \frac{1}{\sqrt {2 \lambda}} ]
\\ 
\sqrt { 2 \lambda} x + 1 - \sqrt 2 a & \text{ if }  x \in [a - \frac{1}{\sqrt {2 \lambda}}, a ]\,.
 \end{cases}
 $$
 
(iii) for $a =2 \sqrt{ \frac{2}{\lambda}} $ there are two solutions, given by the two functions $\overline u _1$ and $\overline u _2$.

Contrarily to the primal problem, the dual problem does not admit easy explicit bounded solutions. 
In particular, the one obtained through the value function ({\it cf.} Remark \ref{r:value}) blows up near the lateral boundary of the cylinder (see Remark \ref{r:blow} for more details).  

Below we give some numerical results obtained,  for $a = 2$,  by using the algorithm \eqref{syst2}.

\begin{figure}
\includegraphics[scale=0.3]{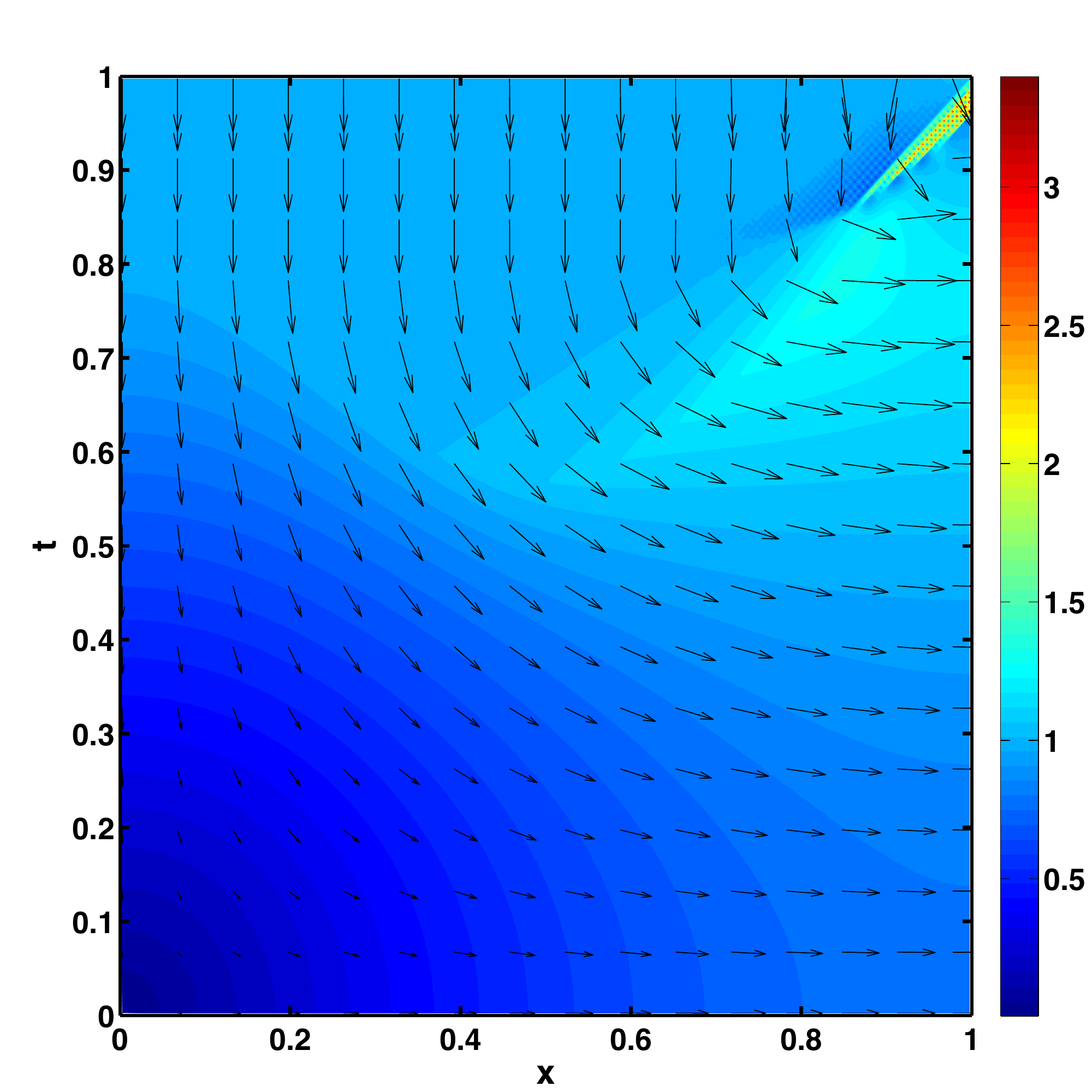}
\qquad \includegraphics[scale=0.3]{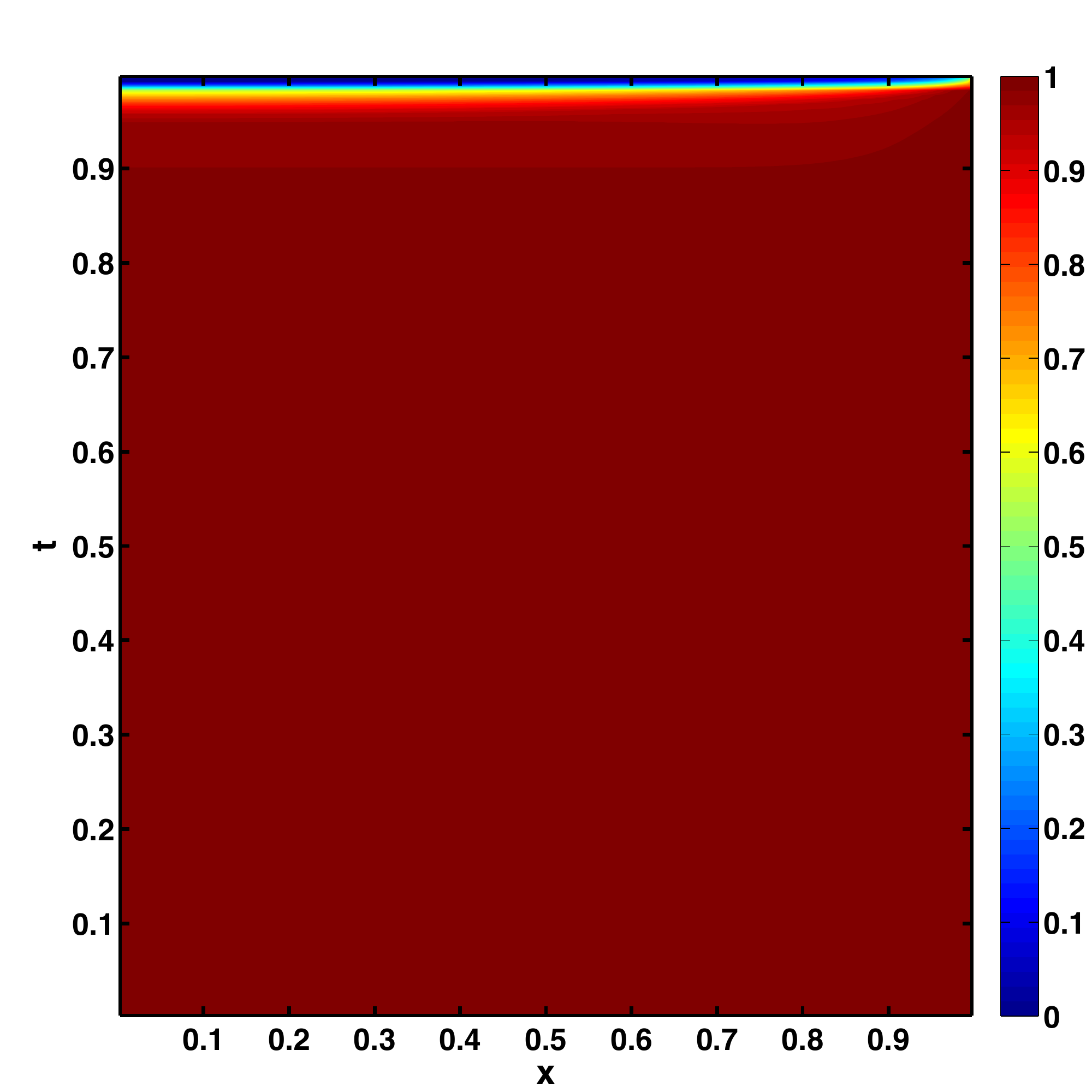}
\caption{Streamlines of $\overline \sigma$ and level sets of $\overline v$ in the case $\lambda = 1$} 
\label{l1}
\end{figure}

\begin{figure}
\includegraphics[scale=0.3]{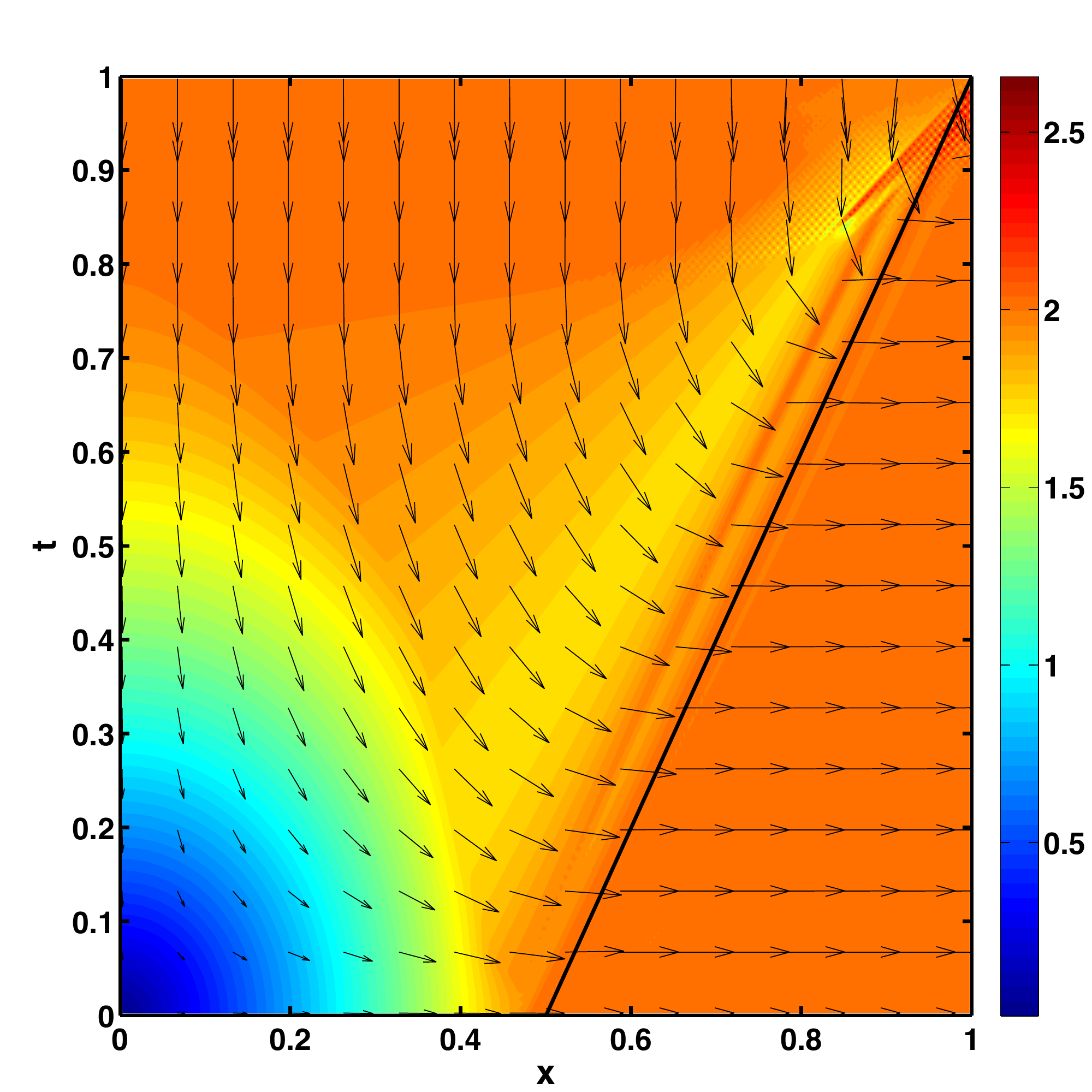}
\qquad \includegraphics[scale=0.3]{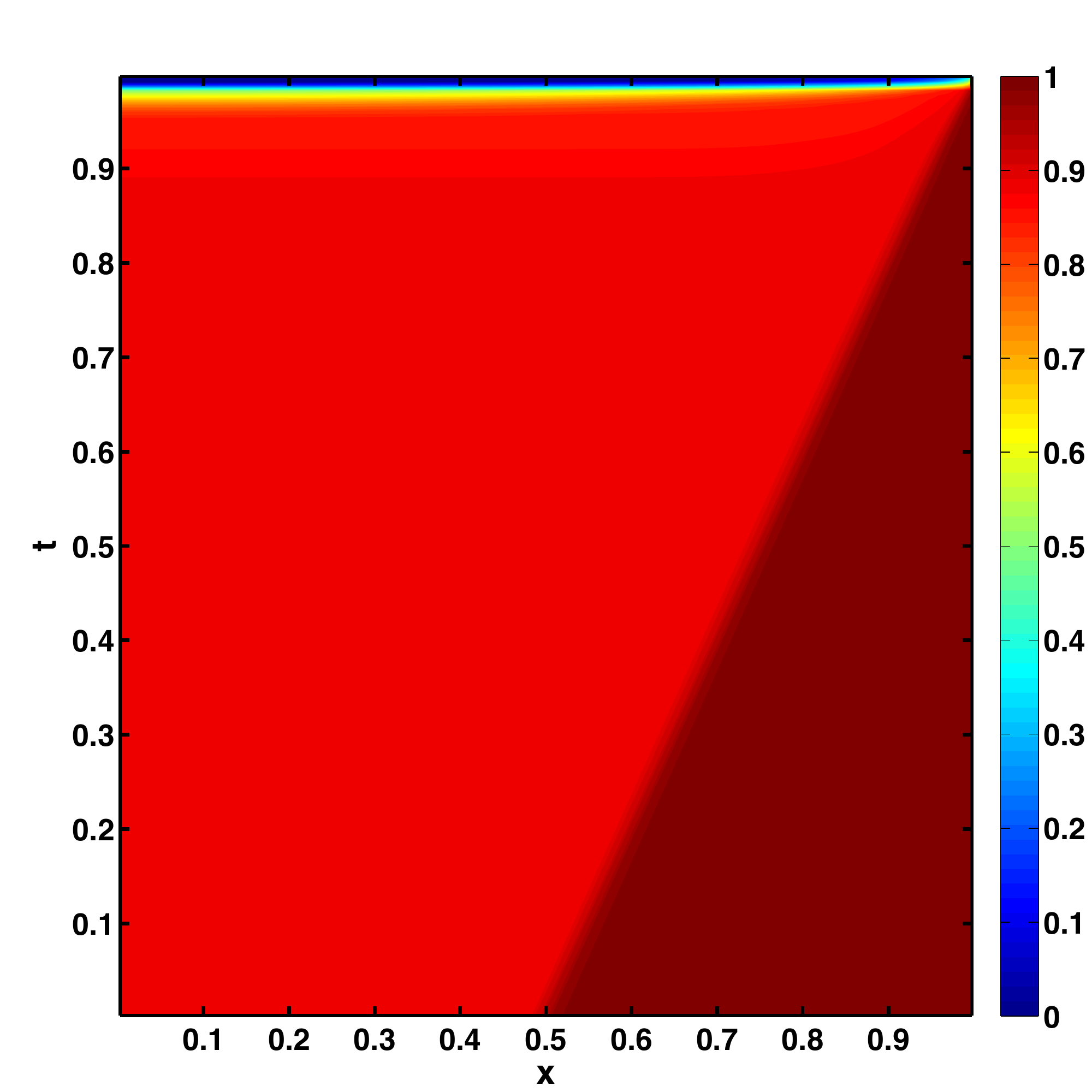}
\caption{Streamlines of $\overline \sigma$ and level sets of $\overline v$ in the case $\lambda = 2$} 
\label{l2}
\end{figure}

\begin{figure}
\includegraphics[scale=0.3]{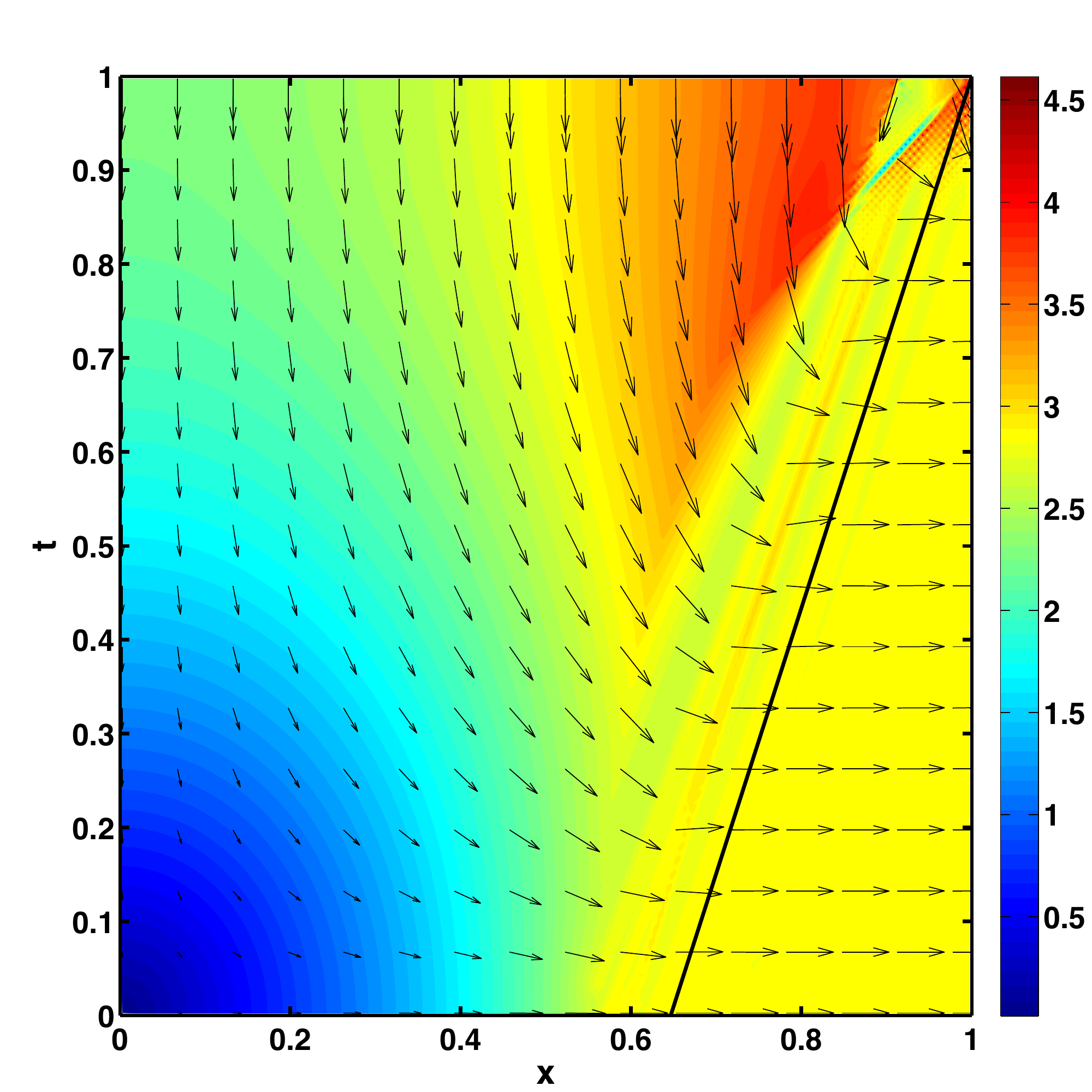}
\qquad \includegraphics[scale=0.3]{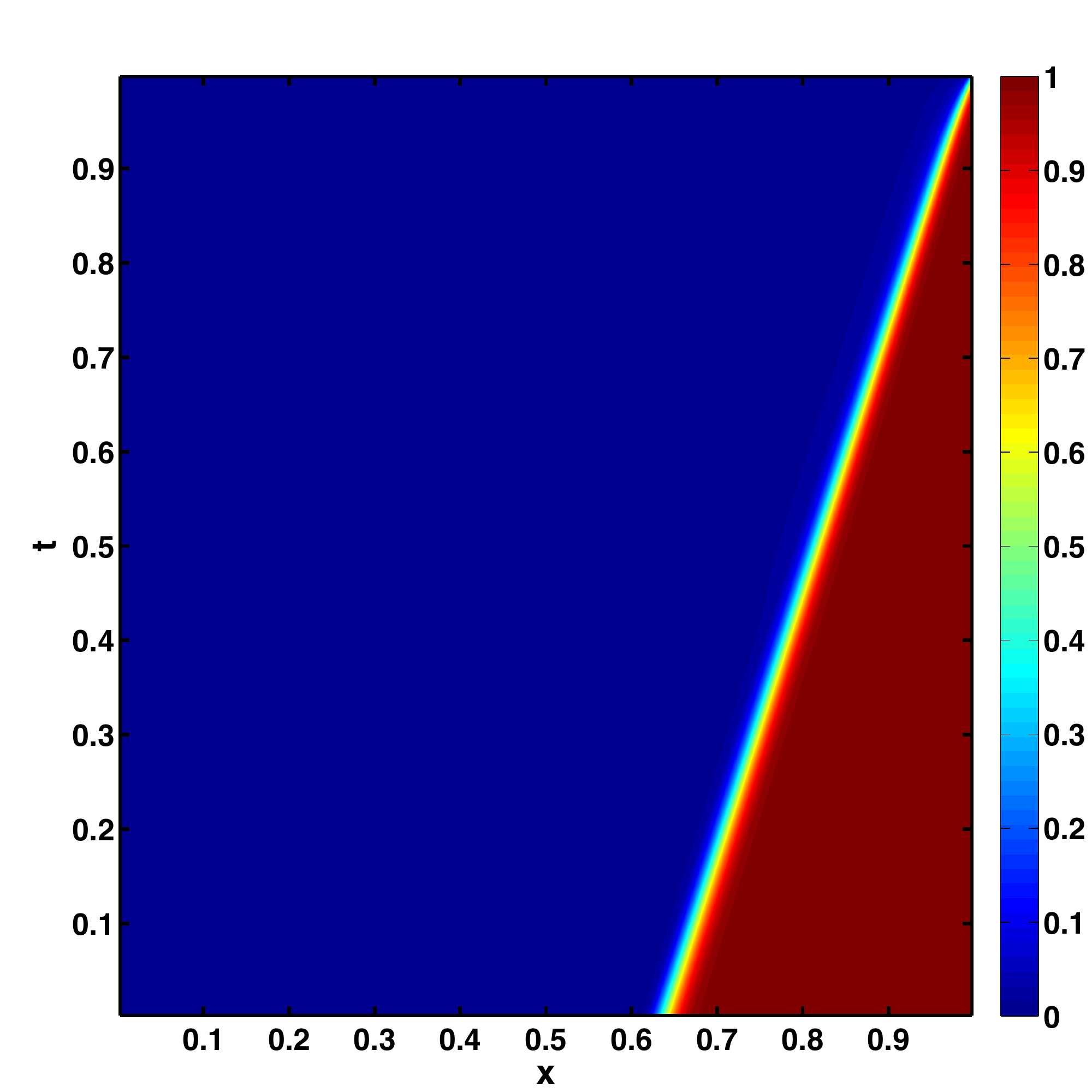}
\caption{Streamlines of $\overline \sigma$ and level sets of $\overline v$ in the case $\lambda = 4$} 
\label{l4}
\end{figure}

Figures \ref{l1}, \ref{l2}, \ref{l4} correspond to three cases $\lambda  =  1, 2, 4$. 
They represent the behaviour of the optimal $\overline \sigma$ and $\overline v$ in each case. Up to a translation of the interval $\Omega = (0, 2)$ into $(-1, 1)$, we can work on the cylinder $(-1, 1) \times (0, 1)$; then, for symmetry reasons, we limit ourselves to plot  our functions on the right part $(0 ,1) \times (0,1)$ of the cylinder. Notice that the most important issue is the location of the discontinuity set of $\overline v$, as the free boundary is given by  the intersection of this set with the horizontal axis. 

For $\lambda = 2$, we recover the two solutions $\overline u _1$ and $\overline u _2$ since the optimal function $\overline v$ exhibits three values (see Figure \ref{l2}, were the regions in blue, red, and brown correspond respectively to the level sets $\{\overline v= 0\}$, $\{\overline v = 0.8886\}$, and $\{\overline v = 1\}$). 

In constrast, for $\lambda =1$ or $\lambda = 4$, when the primal problem admits a unique solution, 
the function $\overline v$ exhibits only two values (see the regions in blue and brown in Figures \ref{l1} and \ref{l4}).

\begin{remark}\label{r:blow} {\rm Let us compute the candidate calibration obtained for problem \eqref{1Dpb} through  the method described in Remark \ref{r:value}. Through some straightforward computations it is easy to obtain that the value function  introduced in \eqref{vf} is given by
$$\begin{array}{ll}
V(x, t)& \displaystyle = \inf \Big \{\int_0 ^ x  f (u, u') \, dt   \ :\ u \in W ^ {1, 2} (0, h)  \, , \ u (0) = 1\, , \ u (x) = t \Big \} 
\\ \noalign{\bigskip}
&\displaystyle=  \min \Big \{ \frac{1}{2} \frac{(t-1) ^ 2 }{x} + \lambda x \, , \, \sqrt { 2 \lambda} ( 1 + |t| ) \Big \}  
\end{array}
$$
Accordingly, the explicit expression of the vector field $\sigma (x, t) := ( \partial _ t  V , - \partial _ x  V)$  reads
$$\sigma (x, t)=
\begin{cases}
\Big ( \frac{t-1}{x} , \frac{1}{2} \frac{(t-1)^2}{x^2} - \lambda \Big )  & \text { if } x \leq \frac{1}{\sqrt {2 \lambda}} ( 1 + \sqrt t ) ^ 2, \ 
t >0 , \text{  or }  x \leq \frac{1}{\sqrt {2 \lambda}} ( 1 +|t| ) , \ 
t < 0
\\ \noalign{\medskip}
(\sqrt {2 \lambda }, 0 )  & \text { if } x > \frac{1}{\sqrt {2 \lambda}} ( 1 + \sqrt t ) ^ 2, \ 
t >0
\\ \noalign{\medskip}
(-\sqrt {2 \lambda }, 0 )  & \text { if } x > \frac{1}{\sqrt {2 \lambda}} ( 1 + |t| ) , \ 
t <0\,,
\end{cases} 
$$
It is easy to check that $\sigma$ satisfies conditions \eqref{s3}-\eqref{s4}. 
However, $V$ is not optimal for the formulation \eqref{MK} of the dual problem because is not Lipschitz; indeed, 
it turns out that $\sigma$ blows up near $x = 0$, see Figure  \ref{fig:sing}
for a plot representing in case $\lambda = 2$ the symmetrized field $\widetilde \sigma (x, t):=\big (\frac{1}{2} [ \sigma ^ x (x, t) + \sigma ^ x ( 2-x, t )], \frac{1}{2} [ \sigma ^ t (x, t) - \sigma ^ t ( 2-x, t )] \big )$ (which also satisfies conditions \eqref{s3}-\eqref{s4}). Again, for symmetry reasons, the plot is restricted to the right half of the cylinder.  
  \begin{figure}
\includegraphics[scale=0.28]{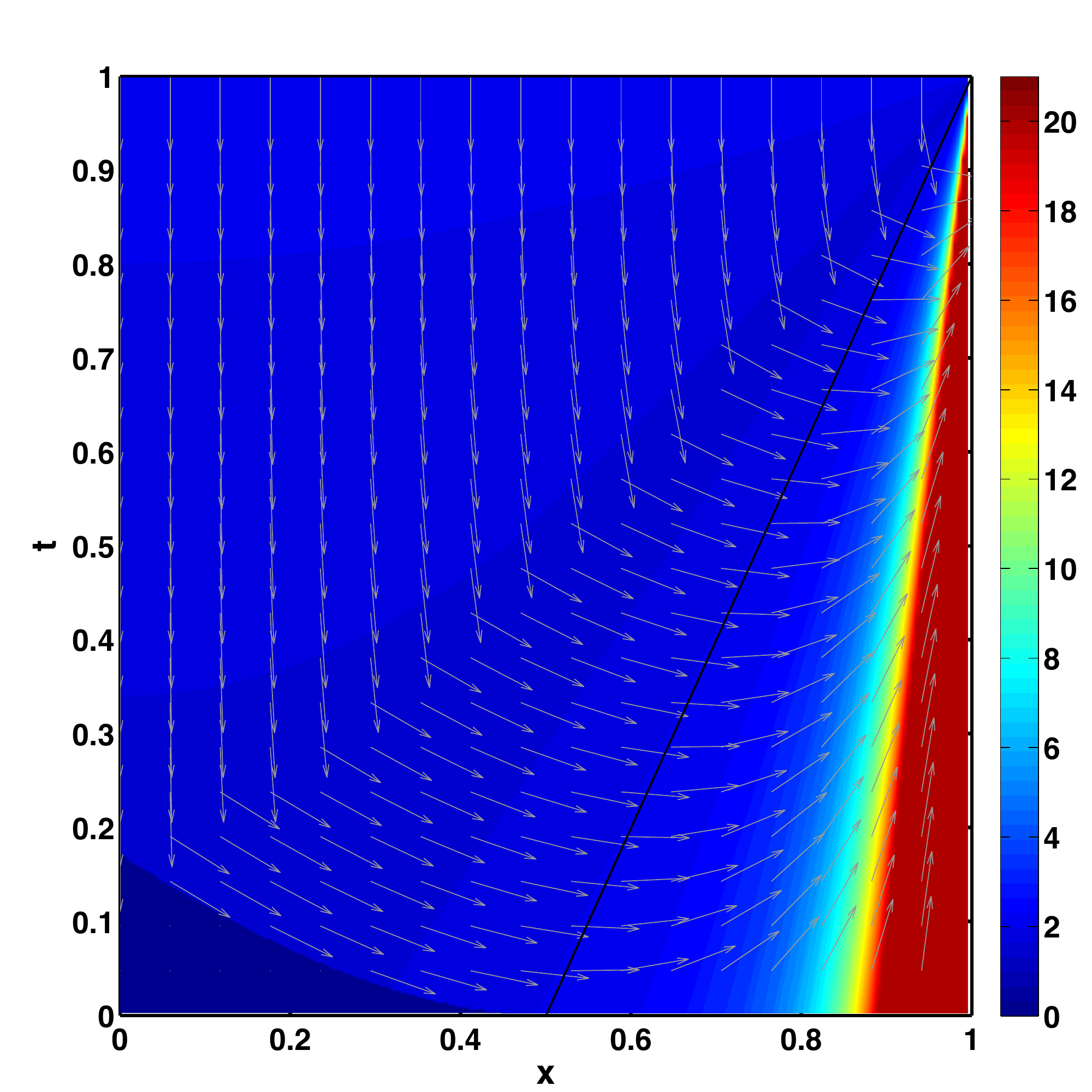}
\caption{Streamlines of the field $\widetilde \sigma$ given by the value function }
\label{fig:sing}
\end{figure}
}
\end{remark}

\subsection{Some simulations in case $N=2$.}

 By using the concavity of the map $\lambda \mapsto \mathcal I (\Omega, \lambda)$ one can check that, 
similarly to the one dimensional case, there exists a critical value $\lambda ^ *= \lambda ^* (\Omega)$ below which the unique solution of the primal problem is $\overline u _1 \equiv 1$, corresponding to the function $\overline v_1 \in \widehat {\mathcal A}$ which vanishes identically in $\Omega \times (0, 1)$.  
For $\lambda = \lambda ^ *(\Omega)$ this solution may coexist with a non constant solution $\overline u _2$, exhibiting a free boundary $E$.

Moreover, the function $\Omega \mapsto \lambda ^* (\Omega)$ turns out to be monotone decreasing  with respect to domain inclusions. In the special case when $\Omega = B _ R := \{ |x| < R \}$, we find the explicit value $\lambda ^ * (B _R) = \frac{2e}{R ^ 2}$. 

We now present some numerical simulations obtained for  $\Omega = (-1, 1) ^ 2$. 
Noticing that $B _ 1 \subset \Omega \subset B _{\sqrt 2}$, we can predict a critical value $\lambda ^ *(\Omega)$ in the interval $(e, 2e)$. In fact, 
by using the second algorithm described above with a mesh size $10 ^ {-2}$ and  by tuning the value of $\lambda$, we obtained the estimate $\lambda ^*(\Omega) \sim 4.7$. 

\begin{figure}
\includegraphics[scale=0.2]{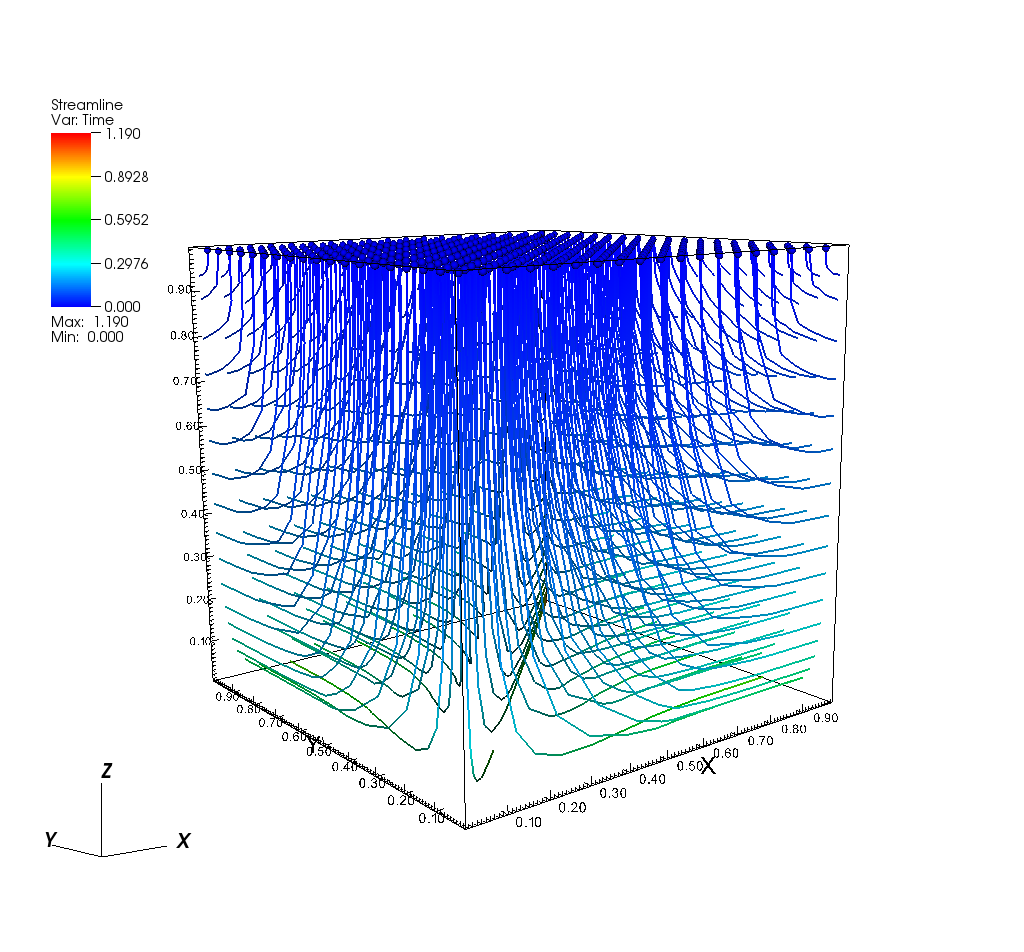}
\caption{Streamlines of $\overline \sigma$ in the case $\lambda = 2e$} 
\label{f:31}
\end{figure}

\begin{figure}
\includegraphics[scale=0.15]{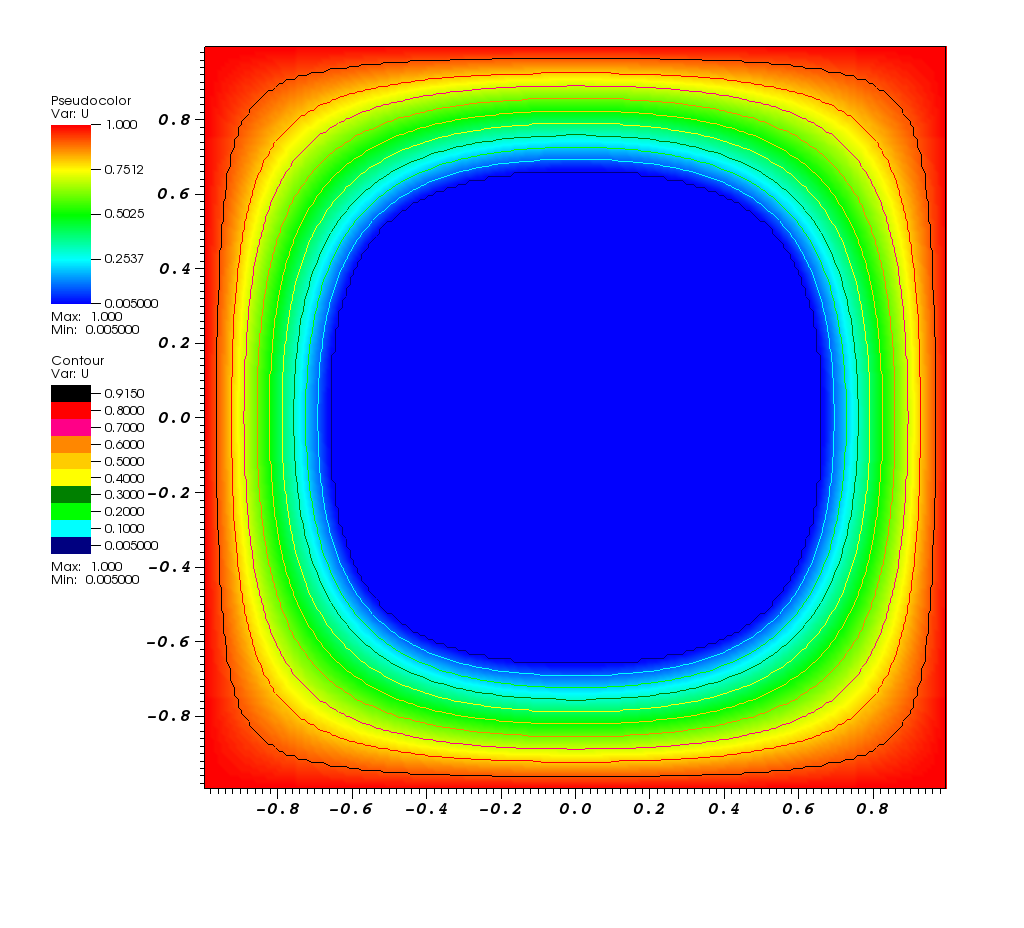} \qquad \includegraphics[scale=0.18]{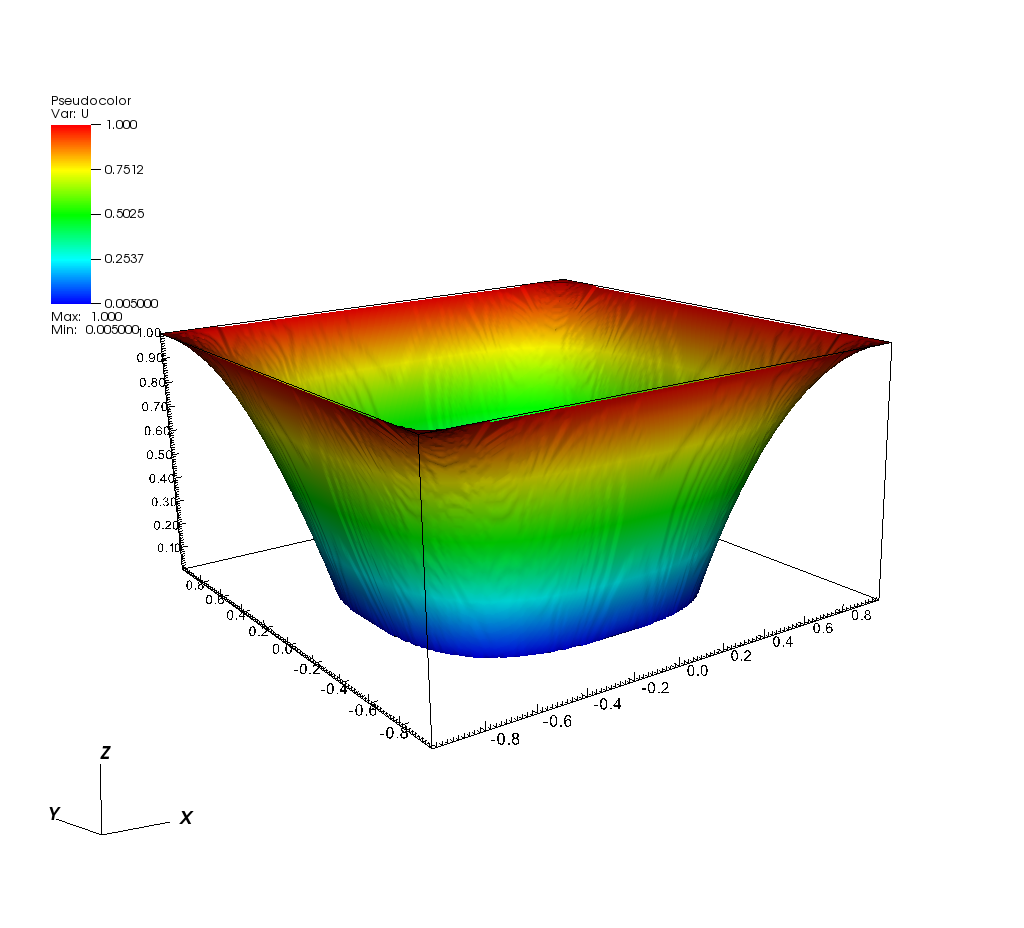}
 \caption{Level sets and plots of $\overline u$ in the case $\lambda = 2e$} 
 \label{f:32}
\end{figure}

In Figures \ref{f:31} and \ref{f:32} we represent respectively the behaviour of the optimal field $\overline \sigma$ and of the optimal function $\overline u$ for $\lambda = 2e$ (for symmetry reasons, Figure \ref{f:31} is referred just to a quarter of $\Omega$, namely to the set $(0, 1)^ 2$). 
Notice that the free boundary is given by the frontier of the region in dark blue.  

\section{Completion of the proofs}\label{secproofs}

In this section  we prove the duality principle stated in Theorem \ref{t:nogap} and the coarea formula stated in Theorem \ref{theo}.

Before starting with the proof of Theorem \ref{t:nogap}, we give some preliminary lemmas.

\begin{lemma}\label{l:B0}  
\begin{itemize}
\item[(i)] If $\Gamma _1 \neq \emptyset$, for every compact neighbourhood $U$ of $\Gamma _1$, there exists  $\sigma_U \in \mathcal B$ such that
$${\rm spt} (\sigma_U ) \subset  U \times \R \,;$$

\medskip

\item[(ii)]  There exists $\sigma _0\in \mathcal B $ such that, for $\delta >0$ sufficiently small, it holds
$$  \| \eta \| _{L^\infty(\Omega \times \R)} \leq \delta\ \Rightarrow \ \sigma _0 + \eta \in \mathcal K\,. $$

\end{itemize}

\end{lemma}

{\it Proof}. (i) Let $A:= U \cap\Omega$,  that we can assume to be Lipschitz. 
Thanks to \eqref{gen_claim}, we know there exists a field $q \in X _ 1 (\Omega)$, with ${\rm spt} (q) \subseteq A$, such that
\begin{equation}\label{f:q} \div q = \frac{|\Gamma_1|}{|A|}\,   \chi _A \hbox  { in } \Omega \quad \, \quad q \cdot \nu _ A = 1\  \hbox{ on } \Gamma_1  
\quad \, \quad q \cdot \nu _ A = 0\  \hbox{ on } \partial A\setminus \Gamma _1  \,.\end{equation}

We define the vector field $\sigma$ by
$$\sigma (x, t) = \big ( - \gamma' ( t) q (x) \, , \, \frac{|\Gamma_1|}{|A|} \, \chi _ A (x)\,  [\gamma (t) - \inf _{\R} \gamma ] + \lambda \big )\, , $$
being $q$ as in \eqref{f:q}, and $\lambda \in \R$ to be chosen later. 
By the choice of $q$, it is immediate that $\div \sigma = 0$ in $\Omega \times \R$ and $\sigma \cdot \nu _\Omega = - \gamma'$ on $\Gamma _ 1 \times \R$. Let us show that  it is possible to choose $\lambda$ such that $\sigma$ belongs to $\mathcal K$. 
By the growth condition from below in \eqref{H3} satisfied by $f$, setting $k := \sup_{\R}  r $, it holds 
$$ f ^ * _ z (t, z^*) \leq  b |z^*| ^ {p'} + k\, , \qquad \hbox{ with }  b = b (\alpha, p) := \frac{1}{p'} \frac{1}{ (\alpha p) ^ {p'-1}}\, , $$ and $-f (t, 0) \leq  r (t) \leq k$.  Therefore, in order that $\sigma$ 
 satisfies   (\ref{s3}) and \eqref{s4}, it is enough to choose $\lambda$ such that 
$$\lambda \geq  b \| \gamma ' \|^{p'} _\infty \| q \| ^{p'} _\infty + k  \,. $$ 

\smallskip
(ii)  Let us consider separately the cases $\Gamma _ 1 \neq \emptyset$ and $\Gamma _ 1 = \emptyset$. 

\smallskip
{\it Case $\Gamma _ 1 \neq \emptyset$.}  We define the vector field $\sigma_0$ by 
$$\sigma _0 (x, t):= \big ( - \gamma' (t) \psi (x) , c_{\Omega} \gamma (t) + \lambda \big ) \,,$$
being $\lambda \in \R$ to be chosen later, $c_\Omega := |\partial \Omega| / |\Omega|$, and $\psi := \nabla w$, with $w$ the unique solution to the boundary value problem $\Delta w = c _\Omega$ in $\Omega$,  $w _\nu= 1$ on $\partial \Omega$. Clearly 
$\sigma _0\in L ^ \infty (\Omega \times \R;  \R ^{N+1} )$ and by construction it holds
$ \div \sigma _0= 0$  in  $\Omega \times  \R$   and $(\sigma _0)^ x \cdot \nu _\Omega =  - \gamma' $  on  $\partial \Omega \times \R $ (thus in particular on $\Gamma_1 \times \R$). 
Let us  check that it is possible to choose $\lambda$ so that 
 $\sigma _0+ \eta$
belongs to $\mathcal K$  if  $\| \eta \| _\infty \leq \delta $. 

We recall that, by our hypothesis \eqref{f:condinf} (in its weaker version asked for $\Gamma _0 \neq \emptyset$), 
and thanks to the boundedness of $r$, 
there exists a constant $m \in \R$ such that
\begin{equation}\label{usem}
c _\Omega \gamma (t) - r (t) \geq m \qquad \text{ for $\mathcal L^1$-a.e. } t \in \R \,.
\end{equation}

In order that $\sigma _0+ \eta$
 satisfies   (\ref{s3}), we need to choose $\lambda$ such that 
$$ 
|q| < \delta \ \Rightarrow \ c_{\Omega} \gamma (t) + \lambda + q ^ t \geq f ^ * _ z (t, - \gamma' (t) \psi (x ) + q ^ x) \quad \text{ for } \mathcal L ^ {N+1} \text{--a.e.}\ (x, t) \in \Omega \times \R \,.
$$
Since, by the growth condition from below in \eqref{H3} satisfied by $f$, it holds $ f ^ * _ z (t, z^*) \leq  b  |z^*| ^ {p'} + r (t)$, 
it is enough to have 
$$|q| < \delta \ \Rightarrow \ c_{\Omega} \gamma (t) + \lambda + q ^ t \geq b  \big | - \gamma' (t) \psi (x ) + q ^ x| ^ {p'} + r (t) \quad \text{ for } \mathcal L ^ {N+1} \text{--a.e.}\ (x, t) \in \Omega \times \R \,.
$$
In  turn, in view of \eqref{usem},  we are reduced to choose $\lambda$ so that
$$|q| < \delta \ \Rightarrow \ m+ \lambda \geq b   \big | - \gamma' (t) \psi (x ) + q ^ x| ^ {p'} - q ^t \qquad \text{ for } \mathcal L ^ {N+1} \text{--a.e.}\ (x, t) \in \Omega \times \R \,,
$$
which is clearly possible since $\psi$ is bounded and $\gamma$ is Lipschitz.

In order that $\sigma _0+ \eta$
 satisfies   (\ref{s4}), we need to choose $\lambda$ such that
$$
|q| < \delta \ \Rightarrow \ c_{\Omega} \gamma (t) + \lambda + q ^ t \geq - f (t, 0) \qquad \forall t \in D \, , \ \text{ for } \mathcal L ^ {N} \text{--a.e.}\ x \in \Omega  \,.
$$
This is possible because, by the growth assumption \eqref{H3}, we have $- f (t, 0) \leq r (t)$, and hence, in view of \eqref{usem}, it is enough to choose $\lambda$ so that $$
|q| < \delta \ \Rightarrow \  m+ \lambda  \geq - q ^ t  \,.
$$

\smallskip
{\it Case $\Gamma _ 1 = \emptyset$.}
We define the vector field $\sigma_0$ simply by 
$$\sigma _0 (x, t):= \big ( 0,  \lambda \big ) \,.$$
Clearly, it holds  
$\sigma _0\in L ^ \infty (\Omega \times \R;  \R ^{N+1} )$ and  $ \div \sigma _0= 0$  in  $\Omega \times  \R$. We have just to choose $\lambda$ so that   $\sigma _0+ \eta$
belongs to $\mathcal K$  if  $\| \eta \| _\infty \leq \delta $. 

In order that $\sigma _0+ \eta$
 satisfies   (\ref{s3}), in view of the inequality $ f ^ * _ z (t, z^*) \leq  b  |z^*| ^ {p'} + r (t)$, 
it is enough to have 
$$|q| < \delta \ \Rightarrow \ \lambda +  q ^ t \geq b  \big | q^x | ^ {p'} + r (t) \quad \text{ for } \mathcal L ^ {N+1} \text{--a.e.}\ (x, t) \in \Omega \times \R \,.
$$
This is clearly possible since we assumed $r (t)$ bounded. 

In order that $\sigma _0+ \eta$
 satisfies   (\ref{s4}), in view of the inequality $- f (t, 0) \leq r (t)$, it is enough  to choose $\lambda$ such that
$$
|q| < \delta \ \Rightarrow \  \lambda + q ^ t \geq r (t) \qquad \forall t \in D \, , \ \text{ for } \mathcal L ^ {N} \text{--a.e.}\ x \in \Omega  \,.
$$
Once again, this is possible thanks to the boundedness of $r$. 
 \qed 

\bigskip

\begin{lemma}\label{l:3mai} For every $\sigma \in X_1 (\Omega \times \R)$ and every $v \in BV _\infty (\Omega \times \R)$, it holds
\begin{equation}\label{f:eta} H ( v) = \sup \Big \{ \int_{\Omega \times \R} (\sigma + \eta ) \cdot Dv \ :\ \eta \in \mathcal D (\Omega \times \R ; \R ^ {N+1}) \, , \ \sigma + \eta \in \mathcal K \Big \}\,.
\end{equation}
\end{lemma}

\proof Let $G (\sigma)$ denote the right hand side of  \eqref{f:eta}. 
 The map $\rho:t\in\R \mapsto G ( t \sigma)$ is convex
(as it is the supremum of affine functions).   By Lemma \ref{l:dualite}, it holds  $\rho( 0 ) = H (v)$, whereas $\rho(t) \leq H (v)$ for every $t$ by Lemma \ref{l:julygen}. It follows that  $\rho(t)$ is constant. We deduce in  particular that $G(\sigma)= \rho(1)=\rho(0) = H(v)$. 
\qed

\bigskip
%

{\bf -- Proof of Theorem \ref{t:nogap}}.

Thanks to the equality \eqref{nog} established in Theorem \ref{l:IJ}, the thesis of Theorem \ref{t:nogap} (namely the equality $\S = \mathcal I$) can be reformulated as 
\begin{equation}\label{ineq} \S =  \inf \Big \{ \widehat E (v) \ :\ v \in \widehat{ \mathcal A}  \Big \} \,.\end{equation}

In order to prove \eqref{ineq}, we introduce on $\mathcal C _ 0 (\Omega \times \R; \R ^ {N+1})$ the perturbation function 
$$\begin{array}{ll} \Phi (\eta) :=  
\inf \Big \{\displaystyle & \displaystyle - \int _{G_{u_0}} \sigma \, \cdot \,  \nu _{u_0} \, d \mathcal H ^N - \int_{\Gamma _1} \gamma (u_0) \, d \mathcal H ^{N-1} 
 \ :\  \displaystyle \sigma \in   X_1  (\Omega \times \R; \R ^ {N+1})  \, , \\ \noalign{\medskip}
&  \displaystyle \ \div \sigma = 0\, , \ \sigma^x \cdot \nu _\Omega= - \gamma ' \hbox{ on } \Gamma _1 \times \R \, , \ \sigma + \eta \in \mathcal K \Big \}\,.$$
\end{array} 
$$

It is easy to check that the map $\eta \mapsto \Phi(\eta)$ is convex. Moreover, in view of the 
choice of admissible fields $\sigma$ in the definition of $\Phi (\eta)$, it holds   
\begin{equation}\label{chaineq}
\S = - \Phi (0) \,. 
\end{equation}

Let us compute $\Phi (0)$. Observe  that $\Phi$ is continuous at $0$: namely,  
for any $\eta$ with $\| \eta \| _\infty \leq \delta$, thanks to Lemma \ref{l:B0} (ii)  it holds  
$$\Phi (\eta) \leq -\int _{G_{u_0}} \sigma _0\, \cdot \,  \nu _{u_0} \, d \mathcal H ^N - \int_{\Gamma _1} \gamma (u_0) \, d \mathcal H ^{N-1} \, .$$
Hence we have \begin{equation}\label{chaineq2}
 - \Phi (0) = - \Phi ^{**} (0) = \min (\Phi ^ *) \,,
\end{equation}
where $\Phi ^*= \Phi ^ * (\lambda)$ denotes the Fenchel conjugate of $\Phi$ in the duality between continuous functions and bounded measures. 
Let us compute $\Phi ^ *$, and let us show that it satisfies
\begin{equation}\label{conjugate}
\Phi ^ * (\lambda) = \begin{cases}
\widehat E (v) & \text{ if } \lambda = Dv \, , \hbox{ with } v \in \A
\\ 
+ \infty & \text{ otherwise. }  
\end{cases}
\end{equation}
Once proved \eqref{conjugate}, our proof will be achieved. Indeed, \eqref{conjugate} implies in particular that $\min (\Phi ^* ) = \min \big \{ \widehat E(v)  :  v \in \widehat {\mathcal A} \big \} $. 
Taking into account \eqref{chaineq} and \eqref{chaineq2}, we deduce that the required equality \eqref{ineq} is satisfied.

In order to establish \eqref{conjugate}, we fix now a bounded vector measure $\lambda$ such that $\Phi ^* (\lambda ) < + \infty$
and we proceed  in three steps. 

\bigskip
{\it Step 1. If $\Gamma _ 1 \neq \emptyset$,  for every compact neighborhood $U$ of $\Gamma_1$ and every bounded continuous  $\psi: \overline \Omega \times \R \to  \R ^ {N+1}$, it holds
\begin{eqnarray}& \langle \lambda - D \1_{u_0}, \psi \rangle = 0 \ \text{ whenever }   \  \div \psi = 0 \hbox{ in } \Omega \times \R \text{ and }  \  \psi = 0 \hbox{ on }  U  \times \R & \label{star} 
\\  \noalign{\medskip}
& \int_{(\Omega \setminus U) \times \R} h _ f (t, \lambda ) < + \infty \,.
& \label{star2}
\end{eqnarray} 
If $\Gamma _ 1 =\emptyset$, conditions \eqref{star} and \eqref{star2} hold true with $U = \emptyset$. 
}

\smallskip
\noindent Assume first that $\Gamma _1 \neq \emptyset$.  Given a compact neighborhood $U$ of $\Gamma_1$  and a function $\psi$ as in \eqref{star},  we consider the vector field $\sigma= \sigma_U + \psi$ with $\sigma_U$ chosen according to Lemma \ref{l:B0}  (i).
Since such $\sigma$ is divergence free and satisfies $\sigma^x \cdot \nu _\Omega= - \gamma ' \hbox{ on } \Gamma _1 \times \R$,  in view of the definition of $\Phi$, one has:  
$$\Phi(\eta) \le  -\int _{G_{u_0}} (\sigma_U + \psi) \, \cdot \,  \nu _{u_0} \, d \mathcal H ^N - \int_{\Gamma _1} \gamma (u_0) \, d \mathcal H ^{N-1}$$
for every smooth field $\eta$ with compact support in $\Omega\setminus U$ such that $\psi + \eta \in \mathcal K$. 

This implies 
$$
\Phi ^ * (\lambda)  \ge \ \langle \lambda, \eta+ \psi \rangle +  \langle  D \1_{u_0}-\lambda, \psi \rangle +
 \int _{G_{u_0}} \sigma_U \, \cdot \,  \nu _{u_0} \, d \mathcal H ^N + \int_{\Gamma _1} \gamma (u_0) \, d \mathcal H ^{N-1} \,,
$$
where we have used the identity $\int _{G_{u_0}} \psi \, \cdot \,  \nu _{u_0}= \langle D \1_{u_0}, \psi \rangle $.
Now by fixing $\psi$ and  taking the supremum  with respect to $\eta$ satisfying the   conditions above, by exploiting Lemma \ref{l:3mai} applied on $\Omega\setminus U$, we deduce that, for a suitable constant $C$, there holds:
$$ \Phi ^ * (\lambda)  \ge \ \int_{(\Omega\setminus U)\times \R} h_f(t,\lambda) + \langle  D \1_{u_0}-\lambda, \psi \rangle + C\ .$$
Thus, since by assumption $\Phi ^* (\lambda )$ is finite, \eqref{star} and \eqref{star2} follow. 

In case $\Gamma _ 1 = \emptyset$, we can repeat the same proof above with $\sigma _U \equiv 0$.

\medskip
 {\it Step 2.  There exists a scalar function $v \in L ^ 1 _{\rm loc} (\Omega \times \R)$, 
with $v(x,\cdot)$ monotone non-increasing, such that $\lambda=Dv $.
Moreover, up to adding a constant to $v$, we have $v \in \widehat {\mathcal A}$,  as it holds:
\begin{eqnarray}
& v \in BV_\infty(\Omega\times\R;[0,1])\, , \quad  v(x,-\infty)= 1\ ,\quad  v(x,\infty)= 0 & \label{i}
\\  \noalign{\smallskip}
& v-v_0 \in L^1(\Omega \times \R) & \label{ii}
\\  \noalign{\smallskip}
& v = \1_{u_0} \hbox{ on } \Gamma_0\times\R\,. & \label{iii} 
\end{eqnarray} 
}

From \eqref{star}, 
since $U$ is arbitrarily small (and empty in case $\Gamma _ 1 = \emptyset$), we infer that the bounded measure $\lambda- D \1 _{u_0}$ is orthogonal to all smooth vector fields $\psi$ which are divergence free and compactly supported in $\Omega\times\R$.
 As $\Omega\times\R$ is simply connected, this implies the existence of a scalar function $v \in L^1 _{\rm loc} (\Omega \times \R)$ such that $\lambda =  Dv$.   
Then, since $\int_{K\times\R} h_f(t,\lambda) < + \infty$ for every compact set $K\subset\Omega$, 
 we infer that $- D_t v$ is a non-negative measure on $\Omega\times \R$, which yields the desired monotonicity property of $v(x,\cdot)$ for a.e. $x\in\Omega$. 

 Let us now prove that $v$ satisfies \eqref{i}, \eqref{ii}, and \eqref{iii}.

To prove \eqref{i},  we choose $\varphi \in \mathcal{D}(\Omega; \R^+)$ and  we set $\psi=(0, \varphi(x))$. Integrating by parts over $\Omega\times(-R,+R)$ and taking into acount that
for a.e. $x\in\Omega$, $v(x,-R+0) -v(x,R-0)$ is non negative and converges increasingly to $var(v(x,\cdot))$ as $R\to+\infty$, we obtain 
$$ \begin{array}{ll} \displaystyle \langle Dv - D v_0, \psi \rangle 
& \displaystyle = \lim_{R\to +\infty} \int_\Omega  \varphi(x) [v(x,R-0) -(v(x,-R+0)-1)]\, dx 
\\ \noalign{\medskip}
& \displaystyle = \ \int_\Omega  \varphi(x) \, (1 -var(v(x,\cdot))\, dx\ .
\end{array} $$
By the arbitrariness of $\varphi$, if we combine the above equality  with \eqref{star} and with the identity $\langle D \1_{u_0} - D v_0, \psi\rangle =0$, we get
\begin{equation}\label{varv}
var(v(x,\cdot))=1 \quad \text{ for a.e.\ $x\in\Omega$}\, .
\end{equation}
Next, we consider a function $\varphi\in \mathcal{D}(\Omega)$ such that $\int_\Omega \varphi\, dx=0$, to which we associate a vector field $q\in L^\infty(\Omega;\R^N)$ such that
$ - \div_x q =\varphi$ in $\Omega$ and  $q\cdot\nu_\Omega= 0$ on $\partial\Omega$. Set:
\begin{equation}\label{dpsi} \psi (x,t) := (\alpha'(t) \, q(x), \alpha(t) \varphi(x))\, ,\quad \text{ with } \ \alpha(t) : = H(t) (1- e^{-t}) 
\end{equation}
(being $H$ the Heavyside function). 
Then, integrating  once more by parts over $\Omega\times(-R,+R)$  and letting $R$ tend to $+\infty$, we obtain
 $$ \langle Dv - D v_0, \psi \rangle = \lim_{R\to +\infty} \int_\Omega  \varphi(x)\,  \alpha(R) \, v(x,R-0)\, dx = \ \int_\Omega  \varphi(x) \, v(x,+\infty)\, dx\ ,$$
 where in the second equality we use dominated convergence taking into account that $|v(\cdot,R-0)|\le 1+ |v(\cdot,t_0)|$ for a suitable $t_0>0$ such that $v(\cdot,t_0)=v(\cdot,t_0\pm 0)$ belongs to $L^1(\Omega)$.
 Then, by applying \eqref{star} to the function $\psi$ introduced in \eqref{dpsi}, and recalling the arbitrariness of the smooth function $\varphi$ with vanishing average, we deduce that $v(x,+\infty)$ is a constant
 that we may fix to be zero. 
  Thus, with the help of \eqref{varv}, we conclude the proof of \eqref{i}. 
  
To prove \eqref{ii}, 
 we fix $\sigma_0 \in \mathcal{B}$  (for instance, we can take the one given by Lemma \ref{l:B0} (ii)). 
Similarly as above, we integrate by parts over $\Omega \times (-R \times R)$, and we obtain:
 $$ \begin{array}{ll}  \displaystyle \langle Dv - D v_0, \sigma_0 \rangle 
 \displaystyle &  \displaystyle = \lim_{R\to +\infty} \bigg(\int_\Omega \left[ \sigma_0^t(x,R) v(x,R-0) + \sigma_0^t(x,-R) (1- v(x,-R+0)\right]\, dx 
\\ \noalign{\medskip}
& \displaystyle \qquad \qquad-  \ \int_{\Gamma_1\times(-R,R)}   \gamma'(t) (v-v_0)\, d \mathcal H^{N-1}(x) dt \bigg)\\
\\ 
& \displaystyle = - \lim_{R\to +\infty} \bigg(\int_{\Gamma_1\times[0,R)}   \gamma'(t) v\, d \mathcal H^{N-1}(x) dt 
\\ \noalign{\medskip}
& \displaystyle \qquad \qquad
- \int_{\Gamma_1\times(-R,0))}   \gamma'(t) (1-v)\, d \mathcal  H^{N-1}(x) dt    \bigg)
\end{array} $$
 where 
 in the second equality we used the fact that $\sigma_0^t$ is bounded together with the convergence of $v(\cdot,R-0)$ and of $1-v(\cdot,-R+0)$ to $0$ in  $L^1(\Omega)$.
 
Now, recalling that $u _s$ are defined as in \eqref{defus}, and using the slicing property \eqref{f:slicel} proved in Proposition \ref{p:slicing}, we can rewrite the above equality as
 
 $$ \begin{array}{ll} \displaystyle \langle Dv - D v_0, \sigma_0 \rangle 
& \displaystyle = - \lim_{R\to +\infty} \int_0^1 ds \bigg(\int_{\Gamma_1\cap\{u_s\ge0\}}   \gamma(u_s\wedge R)\, d \mathcal  H^{N-1}(x) 
\\ \noalign{\medskip  } & \displaystyle \qquad \qquad +  \int_{\Gamma_1\cap\{u_s<0\}}   \gamma(u_s\vee -R)\,  d \mathcal H^{N-1}(x)     \bigg)
\\ \noalign{\medskip}
& \displaystyle = - \lim_{R\to +\infty} \int_0^1 ds\int_{\Gamma_1}  \gamma(u_s^R) \, d \mathcal H^{N-1}(x)    \,
\end{array} $$
 
 where $u_s^R := (u_s \wedge R) \vee -R$. 
Clearly $u_s^R\to u_s$ as $R\to+\infty$. Then, since  $\gamma$ is assumed to be bounded from below, by applying Fatou's Lemma we get 
$$ \langle Dv_0 - D v, \sigma_0 \rangle =  \liminf_{R\to +\infty}  \int_0 ^ 1 \, ds \int_{\Gamma_1}  \gamma(u_s^R) \, d \mathcal H^{N-1}(x)  \ge  \int_0 ^ 1 \, ds \int_{\Gamma_1}  \gamma(u_s) \, d \mathcal H^{N-1}(x)\,.$$
Recalling that $ \langle Dv - D v_0, \sigma_0\rangle$ is finite, we infer that 
\begin{equation}\label{info1}
\int_0 ^1 \, ds  \int_{\Gamma_1}  \gamma(u_s) d \mathcal  H^{N-1}(x) <+\infty\,.
\end{equation}

Now, by Proposition \ref{l:coarea} and Step 1, we know that  
\begin{equation}\label{info2}
\int_0^1  \, ds \int_\Omega f(u_s,\nabla u_s)\,dx = \int_{\Omega\times\R} h_f (t,Dv) <+\infty \,.
\end{equation}
Notice in particular that, in case $\Gamma _ 1 = \emptyset$, 
the last inequality follows from \eqref{star2}, applied with $U = \emptyset$. 
In case $\Gamma _ 1 = \emptyset$, we can still apply \eqref{star2} by letting $\Omega \setminus U$ increase to $\Omega$; this is possible thanks to the fact that 
$h _f(t, \lambda)$ is bounded below by a multiple of the total variation of $\lambda$.  The last assertion is easily checked, since, for all $(t, q)$ with $q \neq 0$ such that $h _ f (t, q) < + \infty$, it holds
$$h _ f(t, q) = - q ^ t f \Big ( t, - \frac{q ^x}{q ^t} \Big )  \geq - q ^ t  \Big | \frac{q ^x}{q ^t} \Big | ^ p + q ^ t r (t) \geq q ^ t r (t) \, , $$
and our assumption \eqref{hypr} ensures that $r (t)$ is bounded. 

Combining \eqref{info1} and \eqref{info2}, we deduce that 
$\int_0^1 E(u_s)\, ds <+\infty$. 
%
In view of the estimate \eqref{f:coerc} obtained in the proof of Proposition \ref{l:existence},  
we deduce that 
$\int_0 ^ 1  \|u _s \| _{ W^ {1,p} (\Omega)} \, ds < + \infty$. 
This implies \eqref{ii}  since
$\int_{\Omega \times \R} |v- v_0| \, dx \, dt = \int _0 ^ 1 \,ds \int _\Omega |u _s| \, dx$. 

To conclude the proof of Step 2, it remains to show \eqref{iii}. To that aim, it is enough to apply 
\eqref{star}. Indeed integrating by parts we obtain  $\int _{\Gamma _0 \times \R} ( v - \1 _{ u _0}) \psi \cdot \nu _\Omega = 0$
 for every bounded continuous function $\psi$ as in \eqref{star}, and the conclusion follows recalling \eqref{gen_claim}.

\medskip
 
{\it Step 3. There holds $\Phi ^ * (\lambda) = \widehat E (v)$}.

\smallskip
Let $\sigma \in X_1(\Omega\times\R;\R^{N+1})$. We observe that,  by Step 2, the duality  bracket $ \sigma \cdot Dv$ is well defined 
({\it cf.}\ \eqref{lf}). Moreover, by Lemma \ref{l:3mai}
it holds
$$ H ( v) = \sup \Big \{ \int_{\Omega \times \R} (\sigma + \eta ) \cdot Dv \ :\ \eta \in \mathcal D (\Omega \times \R ; \R ^ {N+1}) \, , \ \sigma + \eta \in \mathcal K \Big \}\,.
$$

We are now ready to compute the Fenchel conjugate of $\Phi$. We have :
$$\begin {array}{ll}
\Phi ^ * (\lambda ) & \displaystyle = \sup \Big \{ \int_{\Omega \times \R} \eta \cdot  Dv  - \Phi (\eta) \ :\  \eta \in \mathcal C _ 0 (\Omega \times \R; \R ^ { N+1} ) 
\Big \} 
\\ \noalign{\medskip}
& \displaystyle = \sup \Big \{ \int_{\Omega \times \R} \eta \cdot  Dv  - \Phi (\eta) \ :\  \eta \in \mathcal D  (\Omega \times \R; \R ^ { N+1} ) 
\Big \} 
\\ \noalign{\medskip}
& = \displaystyle \sup \Big \{  \int_{\Omega \times \R} (\eta + \sigma) \cdot  Dv  + \int _{G_{u_0}} \sigma \, \cdot \,  \nu _{u_0} \, d \mathcal H ^N + \int_{\Gamma _1} \gamma (u_0) \, d \mathcal H ^{N-1} 
 - \langle Dv, \sigma \rangle \ :\  
 \\ \noalign{\medskip}
 & \qquad 
\eta \in \mathcal D (\Omega \times \R; \R ^ { N+1} ) \, , \ 
\sigma \in  X_1  (\Omega \times \R ; \R ^ {N+1} ) \, , \\ \noalign{\medskip}
 & \qquad \displaystyle \div \sigma = 0\, , \ \sigma^x \cdot \nu _\Omega= - \gamma ' \hbox{ on } \Gamma _1 \times \R \, , \ \ \sigma + \eta \in \mathcal K  \Big \} \\ \noalign{\medskip}
& = \displaystyle \int_{\Omega \times \R} h _ f (t, Dv)   + \int_{\Gamma _1} \gamma (u_0) \, d \mathcal H ^{N-1} 
 + \sup \Big \{  \langle D \1 _{u_0}- Dv, \sigma \rangle \ :\  
 \\ \noalign{\medskip}
 & \qquad 
  \sigma \in  X_1  (\Omega \times \R ; \R ^ {N+1} ) \, , \  \div \sigma = 0\, , \ \sigma^x \cdot \nu _\Omega= - \gamma ' \hbox{ on } \Gamma _1 \times \R \Big \}\,,
\end{array}
$$
 where :
 \begin{itemize}
 \item[--] the first equality is just the definition of $\Phi^*$;
 \item[--] the second equality follows from the density of $\mathcal D (\Omega \times \R; \R ^ { N+1} )$ in $\mathcal C _ 0 (\Omega \times \R; \R ^ { N+1} )$ and from the  continuity of the convex function $\Phi$ at $0$;
 \item[--] the third equality is just the definition of $\Phi$; 
 \item[--] the fourth equality holds by Lemma \ref{l:3mai}.
 \end{itemize} 

Finally we observe that,  thanks to \eqref{casev},
the expression of $\Phi ^ *(\lambda)$ appearing in the fourth equality above coincides with
$$\int_{\Omega \times \R} h _ f (t, Dv)   +  \int _{\Gamma _ 1 \times \R} (v - v_0) \gamma ' (t) \, d \mathcal H ^N = \widehat E ( v)\ .$$
Since from Step  2 we already know  that $v \in \widehat{\mathcal A}$, 
the proof of \eqref{conjugate} is complete. 

\qed

{\bf -- Proof of Theorem \ref{theo}} \label{proofcoarea}

Throughout the proof we set for brevity 
$$u_t(x) := {\chi}_{\{u>t\}}(x)\,.$$

Let us first show that the map $t \mapsto J (u_t)$ is Lebesgue measurable. 

For every fixed open set $V \subset \subset A$, consider the function of a real variable  defined by
 $$\psi _V( t):= 
\int_V u_t(x) dx \,.$$ 
Clearly $\psi_V$ is monotone decreasing, non negative and bounded; in particular, it turns out to be
continuous on $\R \setminus D_V$,  where $D_V$ is a countable subset of $\R$ (depending on $V$). 
Moreover, since
$$
\forall t\, ,  \quad \forall \delta >0\, ,  \hspace{4mm} \int_V |u_t - u_{t+\delta}| dx =
\psi_V(t) - \psi_V(t+\delta)\, , 
$$
the map $t \mapsto u_t$ is continuous from $\R \setminus D_V$ to
$L ^ 1 (V)$. Then, by considering an increasing sequence of open sets $V_h \uparrow A$, and exploiting the assumption that  $J$ is  lower semicontinuous on $L ^ 1 _{\rm loc} (A)$, we obtain that 
the map $t \mapsto J (u _t)$ is lower semicontinuous on $\R \setminus D$,  with $D = \cup _h D_{V_h}$ countable. 
Consequently,  the map $t \mapsto J (u _t)$ is Lebesgue-measurable on $\R$.

\medskip
We now prove separately the inequality $J (u) \leq \int _\R J (u _ t ) \, dt$ and its converse. 

\medskip
{\it -- Proof of the inequality $J (u) \leq \int _\R J (u _ t ) \, dt$ } 

Since $J$ is convex, lower semicontinuous, and proper (recall that by assumption $J(\chi_A) = 0$),  we have  $J^{**}=J$, where $J ^{**}$  is the Fenchel biconjugate in the duality between $L ^ 1 _{\rm loc}(A)$ and the space 
$L ^ \infty _c (A)$ of bounded functions with compact support. 
Namely, 
\begin{equation}\label{f1}
J(u)  = \, J^{**}(u)  =  \sup \Big \{ \int_A u w dx - J^*(w) \ :\ w \in L^\infty_c(A) \Big \}\,.
\end{equation}

Let us compute $J ^*$. We claim that \begin{equation}\label{F*}
J ^* (w) = \begin{cases}
0 & \hbox{ if } w \in X \\
+ \infty  & \hbox{ otherwise\, }
\end{cases}
\end{equation}
for some nonempty closed convex set $X \subseteq \big \{ w \in L ^\infty _c (A) \, :\, \int _A w \, dx = 0\big \}$

We begin by showing that $J ^*$ takes only the values $0$ and $+ \infty$.  By definition, 
there holds
$$
J^*(w) = \sup \Big \{ \int_A u w dx - J(u)\ :\ u \in L ^ 1 _{\rm loc}(A) \Big \}  \qquad \forall w \in L ^\infty _c (A)\,.
$$
Let $w \in L ^\infty _c (A)$ be fixed. If $J ^ * (w) \neq 0$, necessarily there exists some $u \in L ^ 1 _{\rm loc}(A)$ such that $\displaystyle \int_A u w dx - J(u) =:
r \not = 0$.  Since for every $\lambda \geq 0$ we have
$\displaystyle \int_A (\lambda u) w dx - J(\lambda u) =  \lambda r$, we infer that 
 $J^*(w) = +\infty$ if $r$ is positive
(by letting $\lambda$ tend to $+ \infty$) and
$J^*(w) \geq 0$ if $r$ is negative (by letting $\lambda$ tend to $0$); moreover, we see that that $J ^ * (w)$ cannot be strictly positive unless it is $+ \infty$ (because if $J ^* (w)>0$ there exists some $u \in L ^ 1 _{\rm loc}(A)$ such that $\displaystyle \int_A u w dx - J(u) >0$, and arguing as above we see that $J ^* (w) = + \infty$). We deduce that $J ^*$ is of the form (\ref{F*}) for some subset $X$ of $L^ \infty _c(A)$. 
Since $J^*$ is convex, lower semicontinuous, and proper, $X$ is a nonempty closed
convex subset of  $L^ \infty _c(A)$. Moreover, if $w \in X$, taking into account that by assumption
 $J(\chi_A) = 0$, we have
 $$0 = J ^ * (w) \geq  \sup_{\lambda \in \R} \Big [  \lambda \, \int_A  w dx\Big ]\, , $$
hence all  functions in $X$ have zero mean on $A$, which concludes the proof of the claim. 

\smallskip
We infer from \eqref{f1} and \eqref{F*} that 
\begin{equation}\label{f2}
J(u) = \sup _{w \in X} \int _A u w \, dx \,.
\end{equation}

As a next step let us show that, for every $w \in X$, setting
$\displaystyle j_w ( t):= \int_A u_t w \, dx$, 
there holds \begin{equation} \label{claim1}
 \int_A u w \,dx = \int_{-\infty}^{\infty}  j_w ( t) \, dt \,.
\end{equation}

To that aim, we apply Fubini's theorem to compute the
following two integrals:
\begin{eqnarray*}
\int_{u\geq 0} u w dx & = & \int_{u \geq 0} \int_0^{u(x)}
w(x) dt dx  = \int_{u \geq 0} w(x) \int_0^{u(x)} dt dx \\
& = & \int_{u \geq 0} w(x) \int_0^{+\infty} u_t(x) dt dx =
\int_0^{+\infty} \int_A u_t w dx dt = \int_0^{+\infty} j_w(t) dt\, , 
\end{eqnarray*}
and\begin{eqnarray*}
\int_{u \leq 0} u w dx & = & - \int_{u\leq 0} w \int^0_{u(x)}
dt dx  = - \int_{u \leq 0} w \int^0_{-\infty} (1-u_t(x)) dt dx \\ 
& = & \int^0_{-\infty} \int_{u\leq 0} (u_t w - w) dx dt =
\int^0_{-\infty} \left[ \int_{u \leq 0} u_t w dx + \int_{u > 0} w dx \right] dt \\
& = & \int^0_{-\infty} \left[\int_{u \leq 0} u_t w dx +
\int_{u > 0} u_t w dx \right] dt = \int^0_{-\infty} j_w (t) dt.
\end{eqnarray*}
Notice that, in the computation of the second integral (fourth equality), we used the fact that $w$ has zero mean on $A$. 

\smallskip
By \eqref{f2} and \eqref{claim1}, we have
\begin{equation}\label{f3} J(u) = \sup_{w \in X} \int_{-\infty}^\infty j_w(t) dt 
 \leq  \int^\infty_{-\infty}
{\mathcal{L}}^1  \esssup _{w \in X} (j_w) \, dt\,.
\end{equation}
Since we know from the first part of the proof that the map $t \mapsto u_t$ is continuous from $\R \setminus D$  to $L ^ 1 _{\rm loc}(A)$  (with $D$ countable), taking into account that $w \in L ^ \infty _c (A)$ we see that $j_w(t)$ is continuous on $\R \setminus D$. Therefore, 
$$
\forall t \in \R \setminus D\, ,\qquad {\mathcal{L}}^1
\esssup_{w \in X}  j_w(t)  =  \sup_{w \in X} j_w(t) \\
 =  \sup_{w \in X} \int_A u_t w dx   = J^{**}(u_t) =
J(u_t)\,,$$ 

so that
\begin{equation}\label{uso2}
\int_{-\infty}^{+\infty} {\mathcal{L}}^1 \esssup_{w
\in X} j_w(t)  \,dt =
\int_{-\infty}^{+\infty} J(u _t) dt\,.
\end{equation} 

By \eqref{f3} and \eqref{uso2}, the proof of the inequality $J(u) \leq \int_{-\infty}^{+\infty} J(u _t) dt$ is achieved.

\bigskip 
{\it -- Proof of the inequality $J (u) \geq \int _\R J (u _ t ) \, dt$}.  Let us start by showing that,  for every $w \in X$,   if $\alpha$ is any function  in
${\mathcal{C}}^\infty(\R,[0,1])$ and $\displaystyle \beta(t) :=
\int_0^t \alpha(s) ds$, there holds:
\begin{equation}\label{f4}
\int_{-\infty}^{+\infty} \alpha(t) j_w(t) dt \leq J(\beta \circ u).
\end{equation}
Indeed, by applying Fubini's theorem we get
$$\begin{array}{ll}\displaystyle
\int_0^{+\infty} \alpha(t)  j_w (t)dt  &\displaystyle =  \int_0^{+\infty}
\alpha(t) \int_A u_t(x) w(x) dx dt 
 =  \int_A w(x) \int_0^{+\infty} \alpha(t) u_t(x) dt dx \\ \noalign{\medskip}
&\displaystyle =  \int_{u \geq 0} w(x) \int_0^{+\infty} \alpha(t) u_t(x) dt dx 
 =  \int_{u\geq 0} w(x) \int_0^{u(x)} \alpha(t) dt dx \\ \noalign{\medskip}
& \displaystyle=  \int_{u \geq 0} \beta \circ u(x) w(x) dx
\end{array}
$$
and
$$
\begin{array}{ll} \displaystyle
\int^0_{-\infty} \alpha(t)  j_w(t) dt &\displaystyle = \int^0_{-\infty}
\alpha(t) \int_A u_t(x) w(x) dx dt 
 =  \int^0_{-\infty} \alpha(t) \int_A
(u_t(x)-1) w(x) dx dt  \\ \noalign{\medskip}
& \displaystyle =  \int_A w(x) \int^0_{-\infty} \alpha(t) (u_t(x)-1) dt dx   =  
 \int_{u < 0} w(x) \int^0_{-\infty} \alpha(t) (u_t(x)-1) dt dx  \\  \noalign{\medskip}
 &\displaystyle =  \int_{u <0} w(x) \int^0_{u(x)} (-\alpha(t)) dt dx =  \int_{u<0} \beta \circ u(x) w(x) dt\,.
\end{array}
$$
Let us remark that, similarly as above,  in the computation of the second integral (second equality), we exploited the fact that $w$ has zero integral mean on $A$. The validity of \eqref{f4} readily follows, since \begin{eqnarray*}
\int^{+\infty}_{-\infty} \alpha(t) j_w(t) dt & = & \int_A
\beta \circ u(x) w(x) dx \\
& \leq & \sup_{w \in X} \int_A
\beta \circ u(x) w(x) dx = J^{**} (\beta \circ u) = J(\beta \circ u)\,.
\end{eqnarray*}
 
We are now ready to   prove the inequality $J (u) \geq \int _\R J (u _ t ) \, dt$. 
 We consider the  $\mathcal C^\infty$-convex subset of
$L ^ 1 _{\rm loc}(A)$ defined by
$$\mathcal H := \left\{\sum^k_{i=1} \alpha_i j_{w_i} \ : \ \alpha_i \in 
\mathcal C^\infty(\R;[0,1]), \ \sum^k_{i=1}
\alpha_i \equiv 1, \ w_i \in X\right\}
$$
For every $u \in L ^ 1 _{\rm loc} (A)$,  if  $v = \sum^k_{i=1} \alpha_i j_{w_i}$ is any function in $\mathcal H$, we have
$$J(u) \geq \sum^k_{i=1} J(\beta_i \circ u) \geq  \sum^k_{i=1} \int_{-\infty}^{+\infty}
\alpha_i(t) j_{w_i}(t) dt  = \int_{-\infty}^{+\infty} v(t) dt  \, $$
where the first inequality holds by assumption (\ref{sub}), and the second one by \eqref{f4}. 

By the arbitrariness of $v \in \mathcal H$, by applying the commutation argument between supremum and integral proved in \cite[Theorem 1]{BoVa}, and recalling the equality (\ref{uso2}), we eventually get
$$\begin{array}{ll}\displaystyle J(u)  \geq \sup_{v \in \mathcal H} \int_{-\infty}^{+\infty} v dt & = \displaystyle
\int_{-\infty}^{+\infty} {\mathcal L}^1
\esssup _{v\in \mathcal H}  v(t) dt  \\ \noalign{\medskip}
 &\displaystyle \geq  \int_{-\infty}^{+\infty} {\mathcal L}^1 \esssup _{w \in X} j_w(t) dt 
 =  \int_{-\infty}^{+\infty} J(u_t) dt\,.
 \end{array}
$$ 
\qed

\bigskip


\bibliographystyle{mybst}
\bibliography{References}

\end{document}